\documentclass[11pt,a4paper]{amsart}
\usepackage{geometry}           
\usepackage[utf8]{inputenc}
\usepackage[english]{babel}
\usepackage{mathrsfs} 
\usepackage{amsfonts} 
\usepackage{enumitem}
\usepackage{amsthm}
\usepackage{amssymb}
\usepackage{amsmath}
\usepackage{hyperref}
\usepackage{mathtools}
\usepackage{overpic}
\usepackage{tikz}
\usetikzlibrary{matrix}
\usetikzlibrary{positioning}
\usepackage{tikz-cd}
\usepackage{mathabx}
\usepackage{graphicx}
\usetikzlibrary{decorations.pathreplacing}
\newcommand{\upartial}{\rotatebox[origin=c]{10}{$\partial$}}

\usepackage{cleveref}

\numberwithin{equation}{section}

\newtheorem{theorem}[equation]{Theorem}
\newtheorem{maintheorem}{Theorem}

\newtheorem{maincorollary}[maintheorem]{Corollary}

\newtheorem{lemma}[equation]{Lemma}
\newtheorem{proposition}[equation]{Proposition}

\theoremstyle{definition}
\newtheorem{definition}[equation]{Definition}
\newtheorem{question}{Question}
\newtheorem*{conjecture}{Conjecture}
\newtheorem{assumption}{Assumption}
\newtheorem{notation}[equation]{Notation}

\newtheorem*{convention}{Convention}

\newtheorem*{observation}{Observation}
\newtheorem{step}{Step}
\newtheorem{stepA}{Step}

\newtheorem{remark}[equation]{Remark}
\newtheorem{claim}[equation]{Claim}

\newcommand{\nc}{\newcommand}
\nc{\dmo}{\DeclareMathOperator}
\dmo{\ra}{\rightarrow}
\dmo{\N}{\mathbb{N}}
\dmo{\F}{\mathbb{F}}
\dmo{\Z}{\mathbb{Z}}
\dmo{\w}{\mathbf{w}}
\dmo{\R}{\mathbb{R}}
\dmo{\C}{\mathcal{C}}
\dmo{\Stab}{Stab}
\dmo{\B}{Big}
\dmo{\Sat}{Sat}
\dmo{\diam}{diam}
\dmo{\Aut}{Aut}
\dmo{\MCG}{MCG}
\dmo{\Min}{Min}
\dmo{\Isom}{Isom}
\dmo{\syl}{syl}
\dmo{\dist}{\mathsf{d}}

\DeclareMathOperator{\lcm}{lcm}

\renewcommand\bar\overline

\title[Embeddability of RAAGs into HHGs]{Embeddability of right-angled Artin groups into hierarchically hyperbolic groups}
\author{Sangrok Oh}
\address{Innovation Center for MathScience Research \& Education, Pusan National University, Busan, Korea}
\email{SangrokOh.math@gmail.com}
\author{Jihoon Park}
\address{Department of Mathematics Education, Kyungpook National University, Daegu, Korea}
\email{ttsiug@knu.ac.kr}
\begin{document}
\begin{abstract}
For a hierarchically hyperbolic group, we give sufficient conditions ensuring that suitable powers of a finite collection of elements generate a right-angled Artin subgroup; under an additional condition, this subgroup is undistorted.
We verify these hypotheses in two natural situations: one modeled on mapping class groups, using structural assumptions on the ambient HHG, and one modeled on RAAGs, requiring the chosen elements to be rigidly fully supported.
Our results recover known embedding theorems for mapping class groups and extend the non-annular undistortion theorem of Clay--Leininger--Mangahas, while also recovering the Kim--Koberda extension graph theorem for RAAGs.
\end{abstract}


\maketitle
\tableofcontents

\section{Introduction}\label{Sec:Introduction}
Given a simplicial graph $\Gamma$, the associated \emph{right-angled Artin group} (RAAG) $A(\Gamma)$ is defined by the presentation:
\[A(\Gamma):= \langle\, \mathcal{V}(\Gamma)\,\mid\, [v_i,v_j]=1\quad \text{whenever }\{v_i,v_j\}\in \mathcal{E}(\Gamma)\,\rangle,\]
where $\mathcal{V}(\Gamma)$ and $\mathcal{E}(\Gamma)$ are the vertex and edge sets of $\Gamma$, respectively.
For a finitely generated group $G$, a fundamental problem in understanding the subgroup structure of $G$ is to determine which RAAGs embed into $G$, as RAAGs interpolate between free groups and free abelian groups.
The RAAG embedding problem has been studied in a number of settings, including mapping class groups \cite{Kob12, CLM12, Seo21, Run21}, RAAGs themselves \cite{DYER90, KK13, CR15, DL24}, and other classes of groups \cite{CP01, JS22, DL24, Sal07, BBM20}. 

\emph{Hierarchically hyperbolic spaces and groups} (HHSs and HHGs), introduced in \cite{BHS17I}, provide a common framework for the geometry of mapping class groups, RAAGs, and many related groups; for background on HHGs, see Section~\ref{Sec:IntroofHHG}. 
Building on the Masur--Minsky machinery for mapping class groups, HHS theory shows that many algebraic and geometric features of these groups follow from a shared collection of axioms. 
Examples include the structure of top-dimensional quasi-flats \cite{BHS21, Bow19}, semihyperbolicity \cite{DMS23, HHP23}, asymptotic dimension \cite{BHS17smallcancellation}, uniform exponential growth \cite{ANSGP24,WY25}, conjugator length bounds \cite{AB23}, and the behavior of stable subgroups \cite{ABD21}.

This paper develops a hierarchically hyperbolic framework for constructing RAAG subgroups in HHGs.
In an HHG, orthogonality of domains provides a natural geometric analogue of commutation.
In general, however, this geometric relation need not translate directly into algebraic commutation in the group.
We identify hypotheses under which axial elements with localized hierarchical dynamics generate, after passing to sufficiently large powers, a RAAG whose defining graph is determined by the orthogonality relations among their supporting domains.

Let $(G,\mathfrak{S})$ be an HHG. 
We consider axial elements \(f_1,\dots,f_m\in G\) such that, for each \(i\), the element \(f_i\) is \emph{fully supported} on an unbounded domain \(U_i\in\mathfrak{S}\), in the sense of Definition~\ref{Def:FullySupportedElement}.  
Suppose moreover that \(\{f_1,\dots,f_m\}\) is \emph{geometrically irredundant}, in the sense of Definition~\ref{Def:irredundant}. 
The notion of being fully supported is modeled on pure mapping classes with connected support, including Dehn twists about essential curves.
Geometric irredundancy is a hierarchical analogue of the usual irredundancy condition: it prevents two distinct generators from having the same large-scale dynamics. 
In the mapping class group case, for pure mapping classes with connected support, it agrees with the condition that no two distinct elements share a common power.


\begin{definition}[Orthogonality graph]\label{Def:OrthogonalityGraph}
Let $(G,\mathfrak{S})$ be an HHG, and let $\Sigma$ be a collection of domains in $\mathfrak{S}$, allowing multiple distinct copies of the same domain.
The \emph{orthogonality graph} $\mathcal{O}^{\Sigma}$ of $\Sigma$ (with respect to $\mathfrak{S}$) is the graph whose vertex set is $\Sigma$, with two distinct vertices joined by an edge if and only if they are orthogonal in $\mathfrak{S}$.
\end{definition}

For the collection of supporting domains \(\Sigma=\{U_i\}_{i=1}^m\), the graph \(\mathcal{O}^{\Sigma}\) records the orthogonality relations among the $U_i$. It is therefore the natural candidate for the defining graph of the RAAG generated by sufficiently large powers of the \(f_i\). The relevance of this graph depends on whether orthogonality of supporting domains is reflected algebraically.
We refer to the corresponding virtual commutation phenomenon as \emph{weak commutativity}; it is defined precisely in Definition~\ref{Def:WeakCommutativity}.
When it holds, after replacing the \(f_i\) by suitable powers, there is a natural homomorphism $A(\mathcal{O}^{\Sigma})\rightarrow G$ sending the standard generator corresponding to \(U_i\) to the corresponding power of \(f_i\). The main problem is to determine when this homomorphism is injective, and when its image is undistorted.

Before stating the technical embedding theorem in its full generality, we record two concrete consequences that are often easier to apply. In the first, weak commutativity is obtained from structural hypotheses on the ambient HHG. In the second, it follows from a stronger support condition on the chosen elements.

\subsection{Main results and applications}
We first treat the structural version. We impose three assumptions on \((G,\mathfrak S)\), stated precisely in Section~\ref{Subsec:AssumptionsOnHHGs}. Assumption~\ref{Assumption:Essential} requires every nesting-minimal domain to be unbounded, Assumption~\ref{Assumption:CoboundedAction} is a coboundedness condition on non-nesting-minimal domains, and Assumption~\ref{Assumption:CentralExtension} requires the algebraic structure of certain domain stabilizers to reflect the geometry of orthogonality.

\begin{maintheorem}[Theorem~\ref{Thm:MainTheoremAssumptiononG_Precise}]\label{Thm:MainTheoremAssumptiononG}
Let $(G,\mathfrak{S})$ be a hierarchically hyperbolic group satisfying Assumptions~\ref{Assumption:Essential},~\ref{Assumption:CoboundedAction} and~\ref{Assumption:CentralExtension}, and let $f_1,\dots,f_m\in G$ be axial elements, each fully supported on an unbounded domain $U_i\in\mathfrak{S}$.
If $\{f_1,\dots,f_m\}$ is geometrically irredundant, then there exists a constant $N = N(\{ f_i \})$ such that the $N$-th powers of the $f_i$ generate a subgroup $H$ isomorphic to $A(\mathcal{O}^{\Sigma})$, where \(\Sigma=\{U_i\}_{i=1}^m\).

Moreover, if each supporting domain $U_i$ of $f_i$ is not nesting-minimal, then after increasing $N$ if necessary, $H$ is undistorted.
\end{maintheorem}

The primary example for Theorem~\ref{Thm:MainTheoremAssumptiononG} is provided by mapping class groups equipped with their standard HHG structures; we recall these structures and fix the relevant conventions in Section~\ref{Subsec:MappingClassGroupCase}.
For a finite-type hyperbolic surface \(S\), every pure mapping class with connected support, including a Dehn twist about an essential curve, is fully supported on the corresponding unbounded domain in the standard HHG structure.
Thus Theorem~\ref{Thm:MainTheoremAssumptiononG} gives the following mapping class group consequence.

\begin{maincorollary}[Theorem~\ref{Thm:MCG}]\label{Cor:MCG}
Let $S$ be a finite-type orientable surface with $\chi(S)<0$, and let $\operatorname{MCG}(S)$ be its mapping class group. For $i=1,\dots,m$, let $f_i\in \operatorname{MCG}(S)$ be a pure mapping class supported on a connected essential subsurface $U_i\subseteq S$. If $\{ f_1, \dots,f_m\}$ is irredundant, then there exists a constant $D = D(\{ f_i \})>0$ such that if $d\ge D$, then 
\[H = \langle f_1^{d},\dots,f_m^{d}\rangle\cong A(\Gamma),\]
where $\Gamma$ is the graph with vertex set $\{ f_i \}$ in which two distinct vertices $f_i,f_j$ are adjacent if and only if $U_i$ and $U_j$ are disjoint.

Moreover, if each $f_i$ is supported on a non-annular subsurface, equivalently, if none of the $f_i$ is a nonzero power of a Dehn twist, then after increasing $D$ in a controlled way, the subgroup $H$ is undistorted.
\end{maincorollary}

\begin{remark}\label{Rmk:MCGcase}
The embedding statement in Corollary~\ref{Cor:MCG}, without the undistortion conclusion, recovers the embedding theorem stated by Runnels \cite{Run21}.
This statement refines earlier results of Koberda \cite{Kob12} and Seo \cite{Seo21}: compared with Koberda's theorem, it gives a more effective choice of powers, and compared with Seo's theorem, it applies beyond Dehn twists to pure mapping classes with connected support.

For the undistortion statement, Corollary~\ref{Cor:MCG} should be compared with the theorem of Clay--Leininger--Mangahas \cite{CLM12}: under the same non-annularity assumption on the support of each mapping class, Corollary~\ref{Cor:MCG} removes their additional pairwise non-nesting assumption. Runnels \cite{Run21} also stated a more general undistortion theorem allowing pure mapping classes with arbitrary connected support, including annuli.

However, the proof of the embedding statement in \cite{Run21} does not explicitly justify one technical point; the required estimate is supplied by the embedding part of our HHG-level argument, and the statement itself does not need to be modified. 

The proof of the undistortion statement in \cite{Run21} contains two further issues: one leaves a gap in the argument, and the other relies on a bounded-orbit assertion that fails in the required generality. A detailed discussion of these issues is given in Remark~\ref{Rmk:RunnelsIssues}. Accordingly, Corollary~\ref{Cor:MCG} records the undistortion conclusion currently justified by the available techniques, including the HHG-level arguments developed here.
\end{remark}

The second version imposes no structural assumption on the ambient HHG, but instead requires the chosen elements to be rigidly fully supported.
An axial element \(f\in G\) is \emph{rigidly} fully supported on an unbounded domain \(U\in\mathfrak{S}\) if it is fully supported on \(U\) and acts trivially on \(\mathcal{C}V\) for every unbounded \(V\bot U\) (see Definition~\ref{Def:FullySupportedElement}).

\begin{maintheorem}[Theorem~\ref{Thm:MainTheoremStrong_Precise}]\label{Thm:MainTheoremStrong}
Let $(G,\mathfrak{S})$ be a hierarchically hyperbolic group.
Let $f_1,\dots,f_m\in G$ be axial elements, each rigidly fully supported on an unbounded domain $U_i\in\mathfrak{S}$, and suppose that $\{f_1,\dots,f_m\}$ is geometrically irredundant. 
Then there exists a constant $N = N(\{ f_i \})$ such that the $N$-th powers of the $f_i$ generate an undistorted subgroup isomorphic to $A(\mathcal{O}^{\Sigma})$, where \(\Sigma=\{U_i\}_{i=1}^m\).
\end{maintheorem}

The primary example for Theorem~\ref{Thm:MainTheoremStrong} is given by RAAGs equipped with their standard HHG structures associated to standard sub-RAAGs; see Section~\ref{Subsec:RAAGCase}.
In this setting, every unbounded domain admits an axial element rigidly fully supported on it, and the orthogonality graph contains the Kim--Koberda extension graph as an induced subgraph. 
Thus Theorem~\ref{Thm:MainTheoremStrong} recovers the Kim--Koberda extension graph theorem.

\begin{maincorollary}[Kim--Koberda \cite{KK13}; Theorem~\ref{Thm:KKthm}]\label{Cor:RAAG}
Let $\Gamma$ and $\Lambda$ be finite simplicial graphs. 
If $\Lambda$ is an induced subgraph of the extension graph $\Gamma^e$, then $A(\Gamma)$ contains an undistorted subgroup isomorphic to $A(\Lambda)$.
\end{maincorollary}

Theorem~\ref{Thm:MainTheoremStrong} also applies to the family of HHGs satisfying the algebraic conditions in \cite[Definitions~3.1 and~3.7]{OPRAAG}. 
This family includes, in particular, fundamental groups of compact special cube complexes equipped with their rich-family HHG structures \cite[Proposition~3.10(2)]{OPRAAG}. We return to these examples in Section~\ref{Subsect:Remarks and questions}.

\smallskip

The recovery of the Kim--Koberda theorem suggests that the orthogonality graph should play, for HHGs, a role analogous to the extension graph for RAAGs.
Indeed, under suitable hypotheses, finite induced subgraphs of the orthogonality graph of \(\mathfrak{S}\) produce RAAG subgroups of \(G\) through our construction.
In the stronger setting where every unbounded domain admits an axial element rigidly fully supported on it, weak commutativity holds automatically; see Proposition~\ref{Prop:weakCommutativityForStrongFSE}.
This gives the following Kim--Koberda-style criterion.

\begin{convention}
For an HHG \((G,\mathfrak{S})\), we denote by ${\mathfrak{S}_{\infty}}$ the set of all unbounded domains in ${\mathfrak{S}}$.     
\end{convention}

\begin{maintheorem}\label{Thm:KKStyle}
Let \((G,\mathfrak{S})\) be an HHG, and let \(\Lambda\) be a finite induced subgraph of \(\mathcal{O}^{\mathfrak{S}_{\infty}}\).
Then:
\begin{enumerate}
\item\label{Item:WeakerOne} If \((G,\mathfrak{S})\) satisfies the weak commutativity property and every unbounded domain admits an axial element fully supported on it, then there exists a subgroup of \(G\) isomorphic to \(A(\Lambda)\).
\item\label{Item:StrongerOne} If every unbounded domain admits an axial element \emph{rigidly} fully supported on it, then there exists an undistorted subgroup of \(G\) isomorphic to \(A(\Lambda)\).
\end{enumerate}
\end{maintheorem}
\begin{proof}
\noindent\emph{Item~\eqref{Item:WeakerOne}.}
For each vertex \(v_i\in\mathcal V(\Lambda)\), let \(U_i\in\mathfrak S_\infty\) be the corresponding domain, and choose an axial element \(f_i\in G\) fully supported on \(U_i\). Since \(\Lambda\) is an induced subgraph of \(\mathcal O^{\mathfrak S_\infty}\), its distinct vertices correspond to distinct domains. Hence the bigsets of the \(f_i\) are distinct, so the collection \(\{f_i\}\) is geometrically irredundant.

If \(v_i\) and \(v_j\) are adjacent in \(\Lambda\), then \(U_i\bot U_j\). Since \((G,\mathfrak S)\) satisfies the weak commutativity property, the embedding part of Theorem~\ref{Thm:MainEmbeddingThm}, stated in the next subsection, implies that sufficiently large powers of the \(f_i\) generate a subgroup of \(G\) isomorphic to \(A(\Lambda)\).

\smallskip
\noindent\emph{Item~\eqref{Item:StrongerOne}} 
Now choose each \(f_i\) to be rigidly fully supported on \(U_i\). As above, the collection \(\{f_i\}\) is geometrically irredundant. Theorem~\ref{Thm:MainTheoremStrong} then implies that sufficiently large powers of the \(f_i\) generate an undistorted subgroup of \(G\) isomorphic to \(A(\Lambda)\).
\end{proof}

Theorem~\ref{Thm:KKStyle} can be further exploited to detect specific subgroups of interest, such as surface subgroups. 
By a \emph{surface group}, we mean the fundamental group of a closed hyperbolic surface. Combining Theorem~\ref{Thm:KKStyle} with well-known sufficient conditions for the existence of surface subgroups inside RAAGs, we obtain the following immediate corollary.

\begin{maincorollary}\label{Cor:Surface subgroup}
Let \((G,\mathfrak{S})\) be an HHG, and suppose that $\mathcal{O}^{\mathfrak{S}_{\infty}}$ contains an induced cycle of length at least $5$. Then:
\begin{enumerate}
\item If $(G,\mathfrak{S})$ satisfies the weak commutativity property and every unbounded domain admits an axial element fully supported on it, then $G$ contains a surface subgroup.
\item If every unbounded domain admits an axial element rigidly fully supported on it, then $G$ contains an undistorted surface subgroup.
\end{enumerate}
\end{maincorollary}
\begin{proof}
Let \(C_n\) be an induced cycle in \(\mathcal O^{\mathfrak S_\infty}\), where \(n\ge5\).
By Theorem~\ref{Thm:KKStyle}, under the hypothesis of Item~\eqref{Item:WeakerOne}, the group \(G\) contains a subgroup isomorphic to \(A(C_n)\); under the hypothesis of Item~\eqref{Item:StrongerOne}, it contains such a subgroup quasi-isometrically embedded.

By \cite{DSS89}, the RAAG \(A(C_n)\) contains the fundamental group of a closed orientable surface of genus $1+(n-4)2^{n-3}$; since \(n\ge 5\), this surface is hyperbolic.
Moreover, the surface in their construction is realized by a compact square complex whose inclusion into the standard cubical complex for \(A(C_n)\) is a local isometry. 
Hence, by \cite[Proposition~II.4.14]{BH99}, the resulting surface subgroup is quasi-isometrically embedded in \(A(C_n)\).

Composing these embeddings gives a surface subgroup of \(G\) in the first case, and an undistorted surface subgroup of \(G\) in the second case.
\end{proof}

\subsection{The technical embedding theorem}\label{Subsec:TechnicalEmbeddingIntro}

We now state the technical theorem underlying the results above. 
It isolates two hypotheses on a geometrically irredundant collection of axial elements, each fully supported on an unbounded domain.
The first is a virtual commutation condition for elements whose supporting domains are orthogonal. The second is a bounded-orbit condition needed to promote the resulting embedding to an undistorted embedding.
A more precise formulation appears as Theorem~\ref{Thm:MainEmbeddingThm_Precise}, after the relevant terminology has been introduced.

\begin{maintheorem}[Main embedding theorem]\label{Thm:MainEmbeddingThm}
Let $(G,\mathfrak{S})$ be an HHG. Fix elements $f_1,\dots,f_m\in G$ such that
\begin{itemize}
\item each $f_i$ is an axial element \emph{fully supported} on an unbounded domain $U_i \in \mathfrak{S}$, and
\item $\{f_1,\dots,f_m\}$ is \emph{geometrically irredundant}.
\end{itemize}    
Set \(\Sigma=\{U_i\}_{i=1}^m\), allowing repetitions.

Suppose that there exists \(k_1>0\) such that \(f_i^{k_1}\) and \(f_j^{k_1}\) commute whenever \(U_i\bot U_j\). This condition is automatic if \((G,\mathfrak S)\) satisfies the weak commutativity property. Then there exists a constant \(N=N(\{f_i\})\) such that the \(N\)-th powers of the \(f_i\) generate a subgroup \(H\) isomorphic to \(A(\mathcal O^\Sigma)\).

If, in addition, there exists \(k_2>0\) such that, for each \(i\), every element of the subgroup \(\langle f_j^{k_2}\mid U_j\bot U_i\rangle\) acts on \(\mathcal CU_i\) with uniformly bounded orbits, then, after increasing \(N\) if necessary, \(H\) is undistorted in \(G\).
\end{maintheorem}

The proof of Theorem~\ref{Thm:MainEmbeddingThm}, given in Section~\ref{Sec:Main Proof}, has two parts. The embedding part is a ping-pong argument on projections to quasi-axes in the hyperbolic spaces associated to the supporting domains. This part is inspired by strategies used in the mapping class group setting by Runnels, but is formulated entirely in the language of HHGs.

The undistortion part is the main technical point of the proof. After writing an element of the source RAAG in syllable form, each syllable determines a translated supporting domain and a translated quasi-axis in the associated hyperbolic space. The difficulty is that distinct syllables may determine quasi-axes in the same associated hyperbolic space, so their contributions cannot be summed independently in the distance formula. We overcome this difficulty by using the syllable order in the RAAG normal form to organize these repeated contributions, thereby obtaining the lower bound required for undistortion. A detailed roadmap for this argument is given at the beginning of Section~\ref{Subsec:UndistortionPart}.


Theorems~\ref{Thm:MainTheoremAssumptiononG} and~\ref{Thm:MainTheoremStrong} are then obtained by verifying, under their respective hypotheses, the two technical conditions in Theorem~\ref{Thm:MainEmbeddingThm}.

\begin{figure}[ht]
\centering
\begin{tikzpicture}[scale=1.1]
\draw[rounded corners=20pt, thick] (-5,-2.4) rectangle (2.7,2);
\node[above] at (-4,1.95) {$\mathcal{P}_1$: Weak commutativity};

\draw[rounded corners=20pt, thick] (-3,-2.7) rectangle (4,1);
\node[above] at (2,-3.25) {$\mathcal{P}_2$: $\forall$ $U\in\mathfrak{S}$ unbounded, $\exists$ an element fully supported on $U$};

\draw[rounded corners=10pt, red, thick] (-4,-1.5) rectangle (1.3,1.7);
\node[above, red] at (-1.8,1.2) {$\mathcal{P}_3$: Assumptions~\ref{Assumption:Essential},~\ref{Assumption:CoboundedAction}~\&~\ref{Assumption:CentralExtension}};

\draw[rounded corners=10pt, green!70!black, thick] (-2.3,-2.2) rectangle (2.5,0.5);
\node[above, green!70!black] at (0.2,-2.2) {$\mathcal{P}'_2$: rigidly fully supported};
\end{tikzpicture}
\caption{Families of HHGs considered in Theorems~\ref{Thm:MainTheoremAssumptiononG},~\ref{Thm:MainTheoremStrong}, and~\ref{Thm:KKStyle}}
\label{Fig:Venn_diagram}
\end{figure}

\subsection{Remarks and questions}\label{Subsect:Remarks and questions}

The results above depend on two kinds of input.
First, the ambient HHG must have enough algebraic structure for orthogonality of domains to be reflected in commutation, at least after passing to powers.
Second, the chosen elements must have sufficiently localized hierarchical dynamics, as expressed by requiring them to be fully supported or rigidly fully supported on unbounded domains.
We organize the following remarks and questions around these two aspects, and then mention several further directions.

\begin{enumerate}[wide]
\item \textbf{Families of HHGs.}
We have introduced several families of HHGs, each reflecting a different aspect of our construction.
Their relationships are summarized in Figure~\ref{Fig:Venn_diagram}; we use the following notation:
\begin{itemize}[leftmargin=2em]
\item \(\mathcal{P}_1\): HHGs satisfying weak commutativity.
\item \(\mathcal{P}_3\): HHG structures satisfying Assumptions~\ref{Assumption:Essential},~\ref{Assumption:CoboundedAction}, and~\ref{Assumption:CentralExtension}; by Proposition~\ref{Prop:Weak commutativity of FSE},
\(\mathcal P_3\subseteq\mathcal P_1\).
\item \(\mathcal{P}_2\) (\(\mathcal{P}'_2\), resp.): HHGs in which every unbounded domain admits an axial element fully supported (rigidly fully supported, resp.) on it; by Proposition~\ref{Prop:weakCommutativityForStrongFSE}, \(\mathcal{P}'_2\subseteq \mathcal{P}_1\cap\mathcal{P}_2\).
\end{itemize}
Since an HHG is a group equipped with a chosen HHG structure, membership in these families is a property of the pair \((G,\mathfrak{S})\), rather than of the abstract group \(G\) alone.
Accordingly, whenever we say that a group belongs to one of these families, the HHG structure under consideration is understood.

We first discuss known and expected examples.

The family \(\mathcal{P}'_2\) is already quite large. In \cite[Definitions~3.1 and~3.7]{OPRAAG}, the authors introduce a class \(\Xi\) of HHG structures satisfying the \(\mathbf F_U\) stabilizers property and the orthogonal decomposition property, motivated by \cite{AB23}.
The \(\mathbf F_U\) stabilizers property guarantees that every unbounded domain \(U\in\mathfrak{S}\) admits an axial element rigidly fully supported on \(U\); see Remark~\ref{Rmk:F_Ustabilizer}.
Hence any HHG satisfying the \(\mathbf F_U\) stabilizers property satisfies the hypotheses of Theorem~\ref{Thm:KKStyle}\eqref{Item:StrongerOne}.
With the HHG structures specified in \cite[Proposition~3.10]{OPRAAG}, the class \(\Xi\) includes hyperbolic groups, compact special groups, finite direct products of groups in \(\Xi\), and groups hyperbolic relative to finite collections of groups in \(\Xi\). In particular, \(\Xi\subseteq\mathcal{P}'_2\).

Mapping class groups of finite-type hyperbolic surfaces with their standard HHG structures provide important examples in \(\mathcal P_3\); more precisely, non-hyperbolic mapping class groups lie in \((\mathcal P_2\cap\mathcal P_3)\setminus\mathcal P'_2\) (see Section~\ref{Subsec:MappingClassGroupCase}).
Several standard examples of short HHGs discussed in \cite[Sections~2.3.2--2.3.4]{Man24} also provide promising candidates for \(\mathcal P_3\), including Artin groups of large and hyperbolic type, many non-geometric graph manifold groups, and extensions of Veech groups. Their standard structures satisfy Assumptions~\ref{Assumption:CoboundedAction} and~\ref{Assumption:CentralExtension}; the remaining issue concerns compatibility with Assumption~\ref{Assumption:Essential}, as explained in Remark~\ref{Rmk:StandardShortHHGs}.
We do not recall the definition of a short HHG here, but note that short HHGs form a special family of the combinatorial HHGs introduced in \cite{BHMS24}; moreover, many HHGs admit combinatorial HHG structures in the sense of \cite{HMS26}.


These examples suggest comparing the element-rich family \(\mathcal{P}'_2\) or \(\mathcal{P}_2\) with the structurally defined family \(\mathcal{P}_3\), and more broadly understanding how these families depend on the chosen HHG structure.


\begin{question}
Which short HHGs, or more generally which combinatorial HHGs, lie in \(\mathcal{P}_3\)? Does every short HHG lie in \(\mathcal P_1\)?
\end{question}

\begin{question}
If \((G,\mathfrak{S})\) and \((G,\mathfrak{S}')\) are two HHG structures on the same group and \((G,\mathfrak{S})\) lies in one of the families above, must \((G,\mathfrak{S}')\) lie in the same family?
\end{question}

\item \textbf{Choice of elements.}
Alongside identifying suitable families of ambient HHGs, it is important to determine which axial elements are fully supported or rigidly fully supported on unbounded domains.
Even when the chosen HHG structure belongs to one of the families above, applying our results requires identifying a suitable collection of such elements.

For the standard short HHG examples discussed above, generators of the infinite cyclic directions already provide fully supported elements. In particular, this applies to the standard generators of Artin groups of large and hyperbolic type; see \cite[Section~2.3.2]{Man24}.
For these examples, the outstanding structural issue is therefore whether one can realize membership in \(\mathcal P_3\) while retaining these fully supported elements, as discussed in Remark~\ref{Rmk:StandardShortHHGs}.
Resolving this issue for these Artin groups would allow our structural theorem to provide another construction of RAAG subgroups, independent of \cite{JS22}.

The search for suitable elements also suggests studying decompositions of powers of axial elements into factors fully supported on single unbounded domains.
This perspective is closely related to the \emph{orthogonal decomposition property} introduced by Abbott--Behrstock \cite[Definition~3.7]{AB23}.
In the formulation used to define \(\Xi\) in \cite[Definition~3.7]{OPRAAG}, the factors are explicitly required to be rigidly fully supported.
These formulations motivate the following question.

\begin{question}
For which HHGs can every infinite-order element, after passing to a power, be expressed as a product of elements fully supported on single unbounded domains?
Which families in Figure~\ref{Fig:Venn_diagram} satisfy this property?
\end{question}

At the same time, one could instead allow an axial element to be fully supported on its \emph{bigset}~\eqref{Eq:bigset}; this could be formalized by modifying Definition~\ref{Def:FullySupportedElement}.
The single-domain assumption is useful here because it provides a quasi-axis in one associated hyperbolic space.
If these quasi-axes could be replaced by hierarchy paths in the proof, then the argument might extend to elements fully supported on their bigsets.
Such an extension could recover a broader range of known RAAG embedding results into RAAGs, and in favorable cases upgrade them to quasi-isometric embeddings.

\item \textbf{Further directions.}
Beyond finding suitable HHGs and suitable elements, the present framework suggests several further questions.

First, Theorem~\ref{Thm:MainTheoremAssumptiononG} requires the supporting domains to be non-nesting-minimal in order to obtain undistortion, whereas Theorem~\ref{Thm:MainTheoremStrong} does not.
This distinction reflects the bounded-orbit condition needed in the undistortion part.
We expect that the non-nesting-minimal condition in Theorem~\ref{Thm:MainTheoremAssumptiononG} should be removable, leading to the following conjecture.

\begin{conjecture}
For an HHG $(G,\mathfrak{S})$ satisfying Assumptions~\ref{Assumption:Essential},~\ref{Assumption:CoboundedAction} and~\ref{Assumption:CentralExtension}, if $f_1,\dots,f_m\in G$ are axial elements, each fully supported on a single unbounded domain, and if the collection $\{f_1,\dots,f_m\}$ is geometrically irredundant, then sufficiently large powers of the $f_i$ generate an \emph{undistorted} subgroup isomorphic to a RAAG.
\end{conjecture}

The loxodromic annular behavior exhibited in \cite[proof of Lemma~5.2]{Mou18}, and discussed in Remark~\ref{Rmk:RunnelsIssues}, illustrates the main obstruction in the nesting-minimal case.
If the conjecture holds, it would provide an HHG-level route toward the annular case in the undistortion statement of Runnels \cite{Run21}.
Moreover, the RAAG embeddings into Artin groups of large and hyperbolic type mentioned above would be quasi-isometric embeddings.

Second, Theorem~\ref{Thm:KKStyle} provides a combinatorial criterion for detecting RAAG subgroups via the orthogonality graph. This suggests studying explicit HHG structures whose orthogonality graphs admit natural combinatorial descriptions. In ongoing work, the authors apply this framework to families of HHGs in \(\mathcal P_3\), including \emph{graph braid groups} and \emph{right-angled Coxeter groups}, for which both the hierarchical structures and the associated orthogonality graphs admit natural combinatorial descriptions.

Finally, besides the method in Corollary~\ref{Cor:Surface subgroup}, other approaches for constructing surface subgroups in RAAGs are known; see \cite{CW, CSS08}.
However, a complete characterization of RAAGs containing surface subgroups remains open.

\begin{question}
Is there a sufficient condition ensuring that an HHG contains a surface subgroup, independent of embeddings of RAAGs known to contain such subgroups?
\end{question}
\end{enumerate}

\subsection{Outline of the paper}
Section~\ref{Sec:Preliminary} introduces various constants and quasi-axes associated to loxodromic isometries of Gromov hyperbolic spaces, provides a focused overview of the HHS/HHG theory with emphasis on aspects most relevant to this work, and recalls normal forms in RAAGs.

Section~\ref{Sec:SufficientConditions} develops the notion of \emph{geometric irredundancy} and the support conditions imposed on axial elements, namely being \emph{fully supported} or \emph{rigidly fully supported} on an unbounded domain, and states a more precise version of Theorem~\ref{Thm:MainEmbeddingThm}.

In Section~\ref{Sec:Projections and MSBI}, we analyze projections onto uniform quasi-geodesics in the hyperbolic spaces associated to unbounded domains and establish a generalized consistency inequality (Lemma~\ref{Lem:MSBI}).

Section~\ref{Sec:Main Proof} contains the proof of Theorem~\ref{Thm:MainEmbeddingThm}, obtained by controlling projections between quasi-axes associated to the specific elements, using the conditions and inequalities developed in the preceding sections.

Section~\ref{Sec:WithAssumptionsonG} proves the two main variants (Theorems~\ref{Thm:MainTheoremAssumptiononG} and~\ref{Thm:MainTheoremStrong}) of the main embedding theorem: one under structural assumptions on the ambient HHG, and the other under a stronger support condition on the chosen axial elements.

Finally, Section~\ref{Sec:Comparison} compares our results with existing work in the settings of mapping class groups and RAAGs (with respect to their standard HHG structures).

\subsection*{Acknowledgements}
We thank Montserrat Casals-Ruiz, Mark Hagen, and Ilya Kazachkov for valuable discussions and for sharing insights from their work in progress. We are also grateful to Giorgio Mangioni and Harry Petyt for their helpful feedback and corrections on an early draft. We thank Jason Behrstock and Inhyeok Choi for sharing their expertise on mapping class groups. Finally, we thank the anonymous referees for their careful reading and detailed comments, which helped us improve the exposition and clarify the scope of the results.

The first draft of this paper was prepared while the first author was a postdoctoral researcher at the University of the Basque Country, supported by the Basque Government grant IT1483-22. The second author was supported by the Samsung Science and Technology Foundation under Project Number SSTF-BA2202-03.

\section{Preliminaries}\label{Sec:Preliminary}
\subsection{Hyperbolic geometry}\label{Sec:ActionsonHyperbolicSpaces}
We recall some standard properties of Gromov-hyperbolic spaces that will be used repeatedly in this paper; for further details, see \cite{BH99} and \cite{KS22}.
Throughout this subsection, we fix a $\delta$-hyperbolic geodesic space $X$ and denote its Gromov boundary by $\upartial X$. Since we often consider quasi-convex subsets of $X$ and nearest-point projections onto them, we adopt the following convention.

\begin{convention}\label{Conv:ProjectionToQC}
For a quasi-convex subspace \(Q\) of a hyperbolic space \(X\), we choose a
map $\mathfrak p_Q:X\to Q$ such that $\dist_X\bigl(x,\mathfrak p_Q(x)\bigr)\le \dist_X(x,Q)+1$ for every \(x\in X\).
We refer to \(\mathfrak p_Q\) as a \emph{nearest-point projection} to \(Q\), with the understanding that it is defined only up to uniformly bounded error by Lemma~\ref{Lem:CoarselyLipschitz}.
For a subset \(A\subseteq X\), we write $\mathfrak p_Q(A)=\{\mathfrak p_Q(a)\mid a\in A\}$.
For \(E\ge 0\) and a subset \(A\subseteq X\), we denote by \(\mathcal N_E(A)\) the \(E\)-neighborhood of \(A\) in \(X\).
When we write \(\mathfrak p_Q(-)\) or \(\mathcal N_E(A)\), the ambient hyperbolic space is understood from context.
\end{convention}

\begin{lemma}\cite[Proposition~3.11, Ch.III.\(\Gamma\)]{BH99}
\label{Lem:CoarselyLipschitz}
Let \(Q\subseteq X\) be \(\lambda\)-quasi-convex.
Then there exists a constant \(R=R(\delta,\lambda)>0\) such that any nearest-point projection \(\mathfrak p_Q:X\to Q\) is \((1,R)\)-coarsely Lipschitz.
Moreover, any two such projections are pointwise \(R\)-close.
\end{lemma}

The next two lemmas concern the composition of nearest-point projections and the quasi-convexity of their images.

\begin{lemma}\cite[Corollary~1.107]{KS22}\label{Lem:Composition of Projection}
Let $A,B\subseteq X$ be $\lambda$-quasi-convex subsets such that $A\subseteq B$. Then there exists a constant $E_{1}=E_{1}(\delta,\lambda)$ such that, for every $x\in X$, we have $$\dist_{X}(\mathfrak{p}_{A}(x),\,\mathfrak{p}_{A}\circ\mathfrak{p}_{B}(x))\le E_{1}.$$
\end{lemma}

\begin{lemma}\cite[Lemma~1.117]{KS22}\label{Lem:Projection of quasiconvex subsets}
Let $Y,Z\subset X$ be $\lambda$-quasi-convex subsets. Then there exists a constant $E'_1=E'_1(\delta,\lambda)$ such that the projection $\mathfrak{p}_{Y}(Z)$ is $E'_1$-quasiconvex in $X$.
\end{lemma}

Finally, the following lemma provides a dichotomy for projections between two quasi-convex subsets, depending on the distance between them. Here, $\dist_{Haus}$ denotes the Hausdorff distance.

\begin{lemma}\cite[Corollaries~1.140 and~1.143]{KS22}\label{Lem:Projections of two quasiconvex}
Let $U,V\subset X$ be $\lambda$-quasi-convex subspaces. Then there exist constants $E_{2}=E_{2}(\delta,\lambda)>0$ and $E_{3}=E_{3}(\delta,\lambda)>0$ such that either
\begin{itemize}
\item 
both projection images $\mathfrak{p}_{U}(V)$ and $\mathfrak{p}_{V}(U)$ have diameter at most $E_{2}$, or 
\item 
$\dist_{Haus}(\mathfrak{p}_{U}(V),\mathfrak{p}_{V}(U))\leq E_{3}$ and $\mathfrak{p}_{U}(V)\subset \mathcal{N}_{E_{3}}(V)\cap U$.
\end{itemize}
\end{lemma}

Next, we examine the isometries of $X$. An isometry $g$ of $X$ is said to be 
\begin{itemize}
\item \emph{elliptic} if it has bounded orbits;
\item \emph{loxodromic} if the map $\Z\to X$ defined by $n\mapsto g^nx$ is a quasi-isometric embedding for some (equivalently, any) $x\in X$;
\item \emph{parabolic} otherwise.
\end{itemize}
Equivalently, $g$ can be characterized by the cardinality of its limit set in $\upartial X$: it is elliptic, parabolic, or loxodromic if the limit set of $g$ has cardinality $0$, $1$, or $2$, respectively. If the limit set of $g$ has cardinality $2$, then we call these limit points $g^{\pm\infty}$.

\begin{definition}\label{Def:StableTranslationLength}
Let $(Y,\dist_Y)$ be a metric space. For an isometry $h$ of $Y$, its \emph{stable translation length} $\tau_Y(h)$ is defined as follows: 
\[\tau_Y(h):=\lim_{n\to\infty}\frac{\dist_Y(y,h^ny)}{n},\]
for some (equivalently, any) $y\in Y$.
\end{definition}

It follows directly from the definition that stable translation length is homogeneous: \(\tau_Y(h^n)=|n|\tau_Y(h)\) for every \(n\in\mathbb{Z}\).

If $g$ is a loxodromic isometry of $X$, then $\tau_X(g)$ is strictly positive. 
Moreover, $g$ acts on $X$ as a translation along a quasi-geodesic axis connecting its two limit points $g^{\pm\infty}$ in $\upartial X$. 
Up to taking a power of $g$, one can choose such an axis to be a quasi-geodesic with constants depending only on the hyperbolicity constant $\delta$, and independent of the particular loxodromic isometry. 
An explicit construction of such an axis appears in \cite[Section~3]{Cou16}, with a brief version given in \cite[Section~1.1]{AB23}.
Accordingly, we state the existence of such an axis without proof.

\begin{lemma}\label{Lem:UniformQuasiAxis}
Let $X$ be a $\delta$-hyperbolic space.
Then there exists a constant $L=L(\delta)>0$ such that for any loxodromic isometry $g$ of $X$ with $\tau_X(g)\ge L$, there exists a $(2,\delta)$-quasi-geodesic $\gamma$ in $X$ satisfying the following:
\begin{enumerate}
\item $\gamma$ connects the limit points of $g$ in $\upartial X$;
\item $\gamma$ is $(L+8\delta)$-quasi-convex and preserved by $g$;
\item $\dist_X(x,g^n(x))\ge |n| \tau_X(g)$ for any $x\in X$ and $n\in\Z$.
\end{enumerate}
\end{lemma}

\begin{definition}[Quasi-axis]\label{Def:Quasi-axis}
Let $g$ be a loxodromic isometry of a $\delta$-hyperbolic space $X$ with $\tau_X(g)\ge L$ for the constant $L$ in Lemma~\ref{Lem:UniformQuasiAxis}. 
Any quasi-geodesic obtained in Lemma~\ref{Lem:UniformQuasiAxis} will be called a \emph{quasi-axis} of \(g\).    
\end{definition}




\subsection{Brief introduction to HHS/HHG}\label{Sec:IntroofHHG}
Hierarchically hyperbolic spaces (HHSs) are quasi-geodesic spaces whose geometry is organized via projections onto a family of hyperbolic spaces.
The definition of an HHS was introduced in \cite{BHS17I} and axiomatized in a more efficient fashion in \cite{BHS19}, as follows. 
Recall that a \emph{quasi-geodesic space} is a metric space where any two points can be connected by a uniform-quality quasi-geodesic.

\begin{definition}\label{Def:HHS}
A quasi-geodesic space $(\mathcal{X},\dist_\mathcal{X})$ is a \emph{hierarchically hyperbolic space (HHS)} if there are a constant $E>0$, an index set $\mathfrak{S}$, whose elements are called \emph{domains}, and a set $\{ \mathcal{C}U \mid U\in \mathfrak{S} \}$ of $E$-hyperbolic spaces $(\mathcal{C}U, \dist_{U})$ such that the following conditions are satisfied:
\begingroup
\renewcommand\labelenumi{(\theenumi)}
\begin{enumerate}
\item \textbf{Projections.} For each $U\in\mathfrak{S}$, there is a \emph{projection} $\pi_{U} : \mathcal{X} \rightarrow 2^{\mathcal{C}U}$ such that
\begin{itemize}
\item each point $x\in\mathcal{X}$ is mapped by $\pi_{U}$ to a set of diameter at most $E$;
\item $\pi_{U}$ is $(E,E)$-coarsely Lipschitz and the image of $\pi_U$ is $E$-quasi-convex in $\mathcal{C}U$.
\end{itemize}
A domain $U$ is said to be \emph{unbounded} (\emph{bounded}, resp.) if $\mathcal{C}U$ has infinite (finite, resp.) diameter.
\item \textbf{Nesting.} $\mathfrak{S}$ is equipped with a partial order $\sqsubseteq$, and either $\mathfrak{S}=\emptyset$ or $\mathfrak{S}$ contains a unique $\sqsubseteq$-maximal domain (usually, denoted by $S$). When $U\sqsubseteq V$, we say that $U$ is \emph{nested} in $V$. For each $U\in\mathfrak{S}$, we denote by $\mathfrak{S}_{U}$ the set of all  $V\in\mathfrak{S}$ such that $V\sqsubseteq U$. Moreover, for all $U,V\in\mathfrak{S}$ with $U\sqsubsetneq V$, there is a specified subset $\rho^{U}_{V}\subset \mathcal{C}V$ with $\diam_{\mathcal{C}V}(\rho^{U}_{V})\leq  E$. There is also a projection map $\rho^{V}_{U}$ from $\mathcal{C}V$ into $2^{\mathcal{C}U}$.

\item \textbf{Orthogonality.} $\mathfrak{S}$ has a relation $\bot$, called \emph{orthogonality}, satisfying the following:
\begin{itemize}
\item $\bot$ is symmetric and anti-reflexive;
\item if $U\sqsubseteq V$ and $V\bot W$, then $U\bot W$; 
\item if $U\bot V$, then $U$ and $V$ are not $\sqsubseteq$-comparable.
\end{itemize}

We also require that for each $T\in\mathfrak{S}$ and each $U\in\mathfrak{S}_{T}$ such that $\{ V\in\mathfrak{S}_{T} \mid V\bot U \}\neq \emptyset$, there exists $W\in\mathfrak{S}_{T}-\{T\}$ such that whenever $V\bot U$ and $V\sqsubseteq T$, we have $V\sqsubseteq W$. This $W$ is called \emph{a container for the orthogonal complement of $U$} in $T$ (it need not be unique, or orthogonal to $U$).

\item \textbf{Transversality and consistency.} If $U,V\in\mathfrak{S}$ are not orthogonal and neither is nested in the other, then we say that $U,V$ are \emph{transverse}, denoted by $U \pitchfork V$. Moreover, if $U\pitchfork V$, then there are non-empty sets $\rho^{U}_{V}\subseteq \mathcal{C}V$ and $\rho^{V}_{U}\subseteq \mathcal{C}U$, each of diameter at most $E$ and satisfying
\begin{center}
$\min \{\dist_{V}(\pi_{V}(x),\rho^{U}_{V}), \dist_{U}(\pi_{U}(x),\rho^{V}_{U})  \}\leq E\quad$ \text{for all }$x\in \mathcal{X}$. 
\end{center}

For two domains $U,V\in\mathfrak{S}$ satisfying $U\sqsubsetneq V$, we have 
\begin{center}
$\min\{\dist_{V}(\pi_{V}(x),\rho^{U}_{V}), \diam_{\mathcal{C}U}(\pi_{U}(x)\cup \rho^{V}_{U}(\pi_{V}(x))  \}\leq E\quad$ \text{for all }$x\in \mathcal{X}$.
\end{center}

The preceding two inequalities are the \emph{consistency inequalities} for points in $\mathcal{X}$.

Finally, if $U\sqsubseteq V$, then $\dist_{W}(\rho^{U}_{W},\rho^{V}_{W})\leq E$ whenever $W\in\mathfrak{S}$ satisfies either $V\sqsubsetneq W$ or $V\pitchfork W$ and $W$ is not orthogonal to $U$.
    
\item \textbf{Finite complexity.} Any set of pairwise $\sqsubseteq$-comparable domains has cardinality $\le E$.
	
\item \textbf{Large links.} Let $U\in \mathfrak{S}$ and $x,y\in \mathcal{X}$. Let $N=E \dist_{U}(\pi_{U}(x),\pi_{U}(y))+E$. Then there exists $\{ T_{i} \}_{i=1,\dots,\lfloor N\rfloor}\subseteq \mathfrak{S}_{U}-\{U\}$ such that for all $T\in\mathfrak{S}_{U}-\{U\}$, either $T\sqsubseteq T_{i}$ for some $i$, or $\dist_{T}(\pi_{T}(x),\pi_{T}(y))\leq E$. 

\item \textbf{Partial Realization.} Let $\{ U_{i}\}$ be a finite collection of pairwise orthogonal domains of $\mathfrak{S}$ and let $p_{i}\in  \mathcal{C}U_{i}$ for each $i$. Then there exists $x\in \mathcal{X}$ such that: 
\begin{itemize}
\item $\dist_{U_{i}}(\pi_{U_{i}}(x), p_{i})\leq E$ for all $i$, and
\item for each $i$ and each $V\in\mathfrak{S}$ with $U_{i}\sqsubsetneq V$ or $U_{i}\pitchfork V$, we have $\dist_{V}(\pi_{V}(x),\rho^{U_{i}}_{V})\leq E$.
\end{itemize}

\item \textbf{Bounded geodesic image.}
If \(U\sqsubsetneq V\), then every geodesic \(\gamma\subseteq\mathcal C V\) either intersects \(\mathcal N_E(\rho_V^U)\), or $\diam_{\mathcal C U}\bigl(\rho_U^V(\gamma)\bigr)\le E$.
 
\item \textbf{Uniqueness.} For each $\kappa\geq 0$, there exists $\theta_{u}=\theta_{u}(\kappa)$ such that if $x,y\in \mathcal{X}$ and $\dist_{\mathcal{X}}(x,y)\geq \theta_{u}$, then there exists $U\in\mathfrak{S}$ such that $\dist_{U}(\pi_{U}(x),\pi_{U}(y))\geq \kappa$.
\end{enumerate}
\endgroup

We often refer to $\mathfrak{S}$, together with the nesting and orthogonality relations, and the projections, as a \emph{hierarchically hyperbolic structure} for the space $\mathcal{X}$.
To emphasize its hierarchically hyperbolic structure $\mathfrak{S}$, the hierarchically hyperbolic space $\mathcal{X}$ is usually written as $(\mathcal{X},\mathfrak{S})$.
\end{definition}

\begin{notation}\label{Notation:Omitting Projection maps}
When we measure distance between a pair of sets, we take the minimum distance between the two sets.
Where it causes no confusion, we often omit projection maps when referring to distance in $\mathcal{C}U$ for $U\in\mathfrak{S}$. For instance, given $x,y\in \mathcal{X}$ and $p\in\mathcal{C}U$, we write $\dist_{U}(x,y)$ in place of $\dist_U(\pi_U(x),\pi_U(y))$ and $\dist_U(x,p)$ in place of $\dist_U(\pi_U(x), p)$.
When necessary for clarity, we may also write $\mathcal{C}(U)$ instead of $\mathcal{C}U$.
\end{notation}

\begin{remark}\label{Rmk:RestrictedDomain}
By the consistency and bounded geodesic image axioms, when $U\sqsubsetneq V$, the map $\rho_{U}^{V}$ need only be defined on $\mathcal{C}V- \mathcal{N}_E(\rho_{V}^U)$.
Accordingly, when referring to $\rho_{U}^{V}$ for $U\sqsubsetneq V$, we take its domain to be $\mathcal{C}V-\mathcal{N}_E(\rho_{V}^U)$, rather than all of $\mathcal{C}V$.
\end{remark}


\begin{remark}
The definition of an HHS in \cite[Definition~1.1]{BHS19} involves several constants in addition to $E$. However, by \cite[Remark~1.6]{BHS19}, one may assume without loss of generality that a single sufficiently large constant dominates all of them. Accordingly, we adopt the above definition using such a universal constant $E$.
\end{remark}

An important consequence of a space being hierarchically hyperbolic is the existence of the distance formula \cite[Theorem~4.5]{BHS19}, which enables one to measure distances in $\mathcal{X}$ via distances in the hyperbolic spaces $\mathcal{C}U$. Given $a,b \ge 0$, the notation $\{\!\!\{a\}\!\!\}_b$ denotes the quantity which is $a$ if $a \ge b$, and is $0$ otherwise. Given $K>0$ and $C\ge 0$, we say $a\asymp_{K,C}b$ if $K^{-1}a-C\le  b \le Ka+C$. 

\begin{theorem}[Distance formula]\label{Thm:DistanceFormula}
Let $(\mathcal{X},\mathfrak{S})$ be an HHS. Then there exists a constant $s_{0}$ such that for all $s\geq s_{0}$, there exists a constant $K_s\ge 1$ such that, for all $x,y\in \mathcal{X}$, 
$$\dist_{\mathcal{X}}(x,y)\asymp_{K_s,K_s} \sum_{U\in\mathfrak{S}} \bigl\{\!\!\bigl\{ \dist_{U}(x,y)\bigr\}\!\!\bigr\}_{s}.$$
\end{theorem}

\begin{definition}[Relevant domain]\label{Def:RelevantDomain}
Let $(\mathcal{X},\mathfrak{S})$ be an HHS.
Given two points $x,y\in \mathcal{X}$ and $\theta>0$, we say that $U\in\mathfrak{S}$ is \emph{$\theta$-relevant} (with respect to $x,y$) if $\dist_{U}(x,y)>\theta$. We denote by $\operatorname{Rel}(x,y,\theta)$ the set of $\theta$-relevant domains. 
\end{definition}

In other words, the distance formula expresses $\dist_{\mathcal{X}}(x,y)$ as being coarsely equivalent to the sum of the distances in $\mathcal{C}U$ between $\pi_U(x)$ and $\pi_U(y)$, taken over the relevant domains $U \in \operatorname{Rel}(x,y,s)$.


For a group $G$ to be called a hierarchically hyperbolic \emph{group}, it must not only admit a hierarchically hyperbolic structure, but also act on this structure by automorphisms. In other words, the defining axioms are required to be $G$-equivariant, with each element of $G$ inducing an automorphism of the underlying hierarchically hyperbolic space.

\begin{definition}[Automorphism]\label{Def:Automorphism}
An \emph{automorphism} of an HHS $(\mathcal{X},\mathfrak{S})$ consists of a quasi-isometry $f:\mathcal{X}\to\mathcal{X}$, a bijective map $f^{\sharp}:\mathfrak{S}\to \mathfrak{S}$ preserving the relations $\sqsubseteq$, $\bot$ and $\pitchfork$, and a set $\{ f^{\diamond}(U):\mathcal{C}U\to \mathcal{C}(f^{\sharp}(U)) \mid U\in\mathfrak{S}\}$ of isometries such that the following two diagrams coarsely commute (with uniform constants) for all $U,V\in \mathfrak{S}$ with $U\not\bot V$:
\[
\begin{tikzcd}[column sep=3pc, row sep=2pc]
\mathcal{X} \arrow[r, "f"]  \arrow[d, "\pi_U"']  &  \mathcal{X} \arrow[d, "\pi_{f^{\sharp}(U)}"]\\
\mathcal{C}U \arrow[r, "{f^{\diamond}(U)}"] & \mathcal{C}(f^{\sharp}(U))
\end{tikzcd}
\quad\quad\quad
\begin{tikzcd}[column sep=3pc, row sep=2pc]
\mathcal{C}U \arrow[r, "{f^{\diamond}(U)}"]  \arrow[d, "\rho^{U}_{V}"'] & \mathcal{C}(f^{\sharp}(U)) \arrow[d, "{\rho^{f^{\sharp}(U)}_{f^{\sharp}(V)}}"]\\
\mathcal{C}V \arrow[r, "{f^{\diamond}(V)}"] & \mathcal{C}(f^{\sharp}(V))
\end{tikzcd}
\]
where $\rho^{\star}_{\ast}:\mathcal{C}(\star)\rightarrow \mathcal{C}(\ast)$ is the (coarse) map given in Definition~\ref{Def:HHS}.
\end{definition}

Whenever no ambiguity arises, we will abuse notation by omitting superscripts and simply writing all maps as $f$. In particular, when we say that $f$ acts on $\mathcal{C}U$, we mean that $f^{\sharp}$ stabilizes a domain $U\in\mathfrak{S}$ and that we are considering the isometry $f^{\diamond}(U):\mathcal{C}U\to\mathcal{C}(f^{\sharp}(U))=\mathcal{C}U$.

\begin{definition}[Hierarchically hyperbolic group]\label{Def:HHG}
A finitely generated group $G$ is a \emph{hierarchically hyperbolic group (HHG)} if it admits an $E$-HHS structure $(G, \mathfrak{S})$ for some constant $E$ and a word metric $\dist_G$, such that the following conditions are satisfied.
\begin{itemize}
\item 
Each element $g\in G$ is an automorphism of $(G,\mathfrak{S})$ via left multiplication;
\item 
The action of $G$ on $\mathfrak{S}$ is cofinite;
\item For all $g,h\in G$ and $U\in \mathfrak{S}$, we have $g\pi_{U}(h)=\pi_{gU}(gh)$. Moreover, for all $U,V\in\mathfrak{S}$ with $U\pitchfork V$ or $V\sqsubsetneq U$, we have $g\rho^{V}_{U}=\rho^{gV}_{gU}$.
\end{itemize}
\end{definition}

\begin{remark}
Definition~\ref{Def:HHG} follows the formulation given in \cite{PS23}. The original definition was introduced in \cite{BHS17I, BHS19}, and it is shown in \cite{DHS20} that the two definitions are equivalent.
\end{remark}


\textbf{For the remainder of this paper, we fix an HHG, together with a universal HHS constant \(E\) and a word metric \(\dist_G\), chosen once and for all.}

The following fact motivates why it is natural to study embeddings of RAAGs into HHGs: RAAGs interpolate between free groups and free abelian groups.

\begin{theorem}[Tits alternative for HHGs, {\cite[Theorem~4.1]{DHS20}}]\label{Thm:TitsAlternative}
Let \((G,\mathfrak S)\) be an HHG, and let \(H\le G\) be a finitely generated subgroup. Then either \(H\) contains a nonabelian free group, or \(H\) is virtually abelian.
\end{theorem}


For an element \(g\in G\), its \emph{bigset} is the collection of domains
\begin{equation}\label{Eq:bigset}
\B(g):=\{\,U\in\mathfrak{S}\mid \diam_{\mathcal{C}U}(\langle g\rangle.x)=\infty
\text{ for some }x\in G\,\},
\end{equation}
where \(\langle g\rangle\) is the cyclic subgroup generated by \(g\).
Since boundedness of the projected orbit is independent of the basepoint, one may equivalently use any \(x\in G\) in this definition. 
We say that \(g\) is \emph{elliptic} if its orbits in \(G\) are bounded, and \emph{axial} if the orbit map \(n\mapsto g^n x\) is a quasi-isometric embedding for some (equivalently, any) \(x\in G\). The following proposition records the resulting axial--elliptic dichotomy and the corresponding behavior on the bigset.

\begin{proposition}[Axial--elliptic dichotomy]\label{Prop:AxialEllipticDichotomy}
An element \(g\in G\) is elliptic if and only if \(\B(g)=\emptyset\), if and only if \(g\) has finite order.

If \(g\) has infinite order, then \(g\) is axial. Moreover,
\begin{itemize}
\item \(g\) permutes \(\B(g)\), which is a nonempty finite pairwise orthogonal subset of \(\mathfrak S\);
\item there exists \(M=M(\mathfrak S)>0\) such that, for every \(U\in\B(g)\), the element \(g^M\) fixes \(U\) and acts loxodromically on \(\mathcal CU\).
\end{itemize}
\end{proposition}
\begin{proof}
The first assertion follows from \cite[Proposition~6.4]{DHS17} and \cite[Lemma~1.21]{AB23}. Now suppose that \(g\) has infinite order. By \cite[Theorem~3.1]{DHS20}, the subgroup \(\langle g\rangle\) is undistorted, so \(g\) is axial, while the first assertion gives \(\B(g)\neq\emptyset\).

The set \(\B(g)\) is \(g\)-invariant by equivariance and is pairwise orthogonal by \cite[Lemma~6.7]{DHS17}. In particular, it is finite, and its cardinality is bounded by a constant \(r=r(\mathfrak S)\); see \cite[Lemma~2.1]{BHS19}. Since \(g\) permutes \(\B(g)\), the power \(g^{|\B(g)|!}\) fixes every element of \(\B(g)\).
Moreover, \cite[Theorem~3.1]{DHS20} shows that \(g^t\) acts loxodromically on \(\mathcal CU\), for every \(U\in\B(g)\), whenever \(t\) is a nonzero multiple of \(|\B(g)|!\). Hence one may take \(M=r!\).
\end{proof}

\subsection{Normal forms in RAAGs}\label{Subsec:Normal forms}
Let $A(\Gamma)$ be a RAAG with \emph{standard generating set} $\mathcal{V}=\{v_1,\dots,v_n\}$, i.e., the vertex set of $\Gamma$. In \cite[Section~4]{CLM12}, Clay--Leininger--Mangahas introduced a \emph{strict partial order} on the set of syllables for any nontrivial element $g\in A(\Gamma)$. In this subsection, we recall this order and introduce the notion of a central form of $g$.

Consider a word $w=x_1^{e_1}\cdots x_k^{e_k}$ with $x_i \in \mathcal{V}$ and $e_i \in \Z$. 
Each $x_i^{e_i}$ is called a \emph{syllable} of $w$. The following \emph{elementary moves} transform a word into another representing the same element of $A(\Gamma)$:
\begin{enumerate}
\item Delete a syllable $x_i^{e_i}$ if $e_i=0$.
\item If $x_i=x_{i+1}$, then replace $x_i^{e_i}x_{i+1}^{e_{i+1}}$ by $x_i^{e_i+e_{i+1}}$.
\item (\emph{Shuffling}) If $[x_i,x_{i+1}]=1$, then replace $x_i^{e_i}x_{i+1}^{e_{i+1}}$ by $x_{i+1}^{e_{i+1}}x_i^{e_i}$.
\end{enumerate}

Let $\Min(g)$ denote the set of words in $\mathcal{V}$ representing $g$ with the fewest syllables. The following theorem of Hermiller--Meier is by now standard.

\begin{theorem}[Hermiller--Meier {\cite{HM95}}]\label{Thm:HM}
Any word representing $g\in A(\Gamma)$ can be transformed to any word in $\Min(g)$ via a sequence of elementary moves, none of which increases the number of syllables.

In particular, if $w,w'\in \Min(g)$, then $w$ can be transformed to $w'$ by a sequence of shufflings.
\end{theorem}

As a consequence, the \emph{syllable length} of $g\in A(\Gamma)$ (with respect to $\mathcal{V}$) is well-defined: it equals the number of syllables in any $w\in\Min(g)$.

For $w = x_1^{e_1}\cdots x_k^{e_k}\in \Min(g)$, denote the set of syllables of $w$ by
\[\syl(w)=\{x_1^{e_1},\dots,x_k^{e_k}\}.\]
Formally, one should record the position of each syllable, i.e.\ $\syl(w) = \{(x_i^{e_i},i)\}_{i=1}^k$, but we shall abuse notation slightly and write simply $x_i^{e_i}\in \syl(w)$. By Theorem~\ref{Thm:HM}, for any $w,w'\in \Min(g)$ there is a natural bijection $\syl(w)\leftrightarrow \syl(w')$. Thus the syllable set of $g$, denoted by $\syl(g)$, is well-defined up to these canonical identifications. This allows us to define a strict partial order on $\syl(g)$.

\begin{definition}[Strict partial order on syllables {\cite{CLM12}}]\label{Def:strictpartialorder}
For $g\in A(\Gamma)$, define a strict partial order $\prec_s$ on $\syl(g)$ by declaring
\[x_i^{e_i} \prec_s x_j^{e_j} \quad \iff \quad x_i^{e_i} \text{ precedes } x_j^{e_j} \text{ in every } w\in \Min(g).\]
\end{definition}

\begin{remark}\label{Rmk:IncomparableSyllablesCommute}
If two syllables of an element of \(A(\Gamma)\) are incomparable with respect to \(\prec_s\), then they are powers of commuting standard generators. 
Moreover, after shuffling, they can be made adjacent in some word in \(\Min(g)\).
\end{remark}

Any $w\in \Min(g)$ can be shuffled to a particularly useful representative, called its \emph{central form}. Recall that for a word $w$ in $\mathcal{V}$, the \emph{support} of $w$, denoted by $\operatorname{supp}(w)$, is the set of standard generators appearing in $w$.

\begin{definition}[Central form]\label{Def:Central form}
Let $g\in A(\Gamma)$. A word $w\in \Min(g)$ is a \emph{central form} of $g$ if
\[w = u_k v_{i_k}^{t_k} u_{k-1} v_{i_{k-1}}^{t_{k-1}}\cdots u_1 v_{i_1}^{t_1},\]
where $v_{i_j}\in\mathcal{V}$ and $u_j$ is a subword such that
\begin{itemize}
\item $\operatorname{supp}(u_j v_{i_j})$ spans a \emph{clique} (or a complete subgraph) in $\Gamma$ for $1\le j\le k$, and
\item $v_{i_j}$ does not commute with $v_{i_{j-1}}$ for $2\le j\le k$.
\end{itemize}
\end{definition}

\begin{lemma}\label{Lem:Central form}
Every element $g\in A(\Gamma)$ admits a central form.
\end{lemma}
\begin{proof}
This follows directly from Theorem~\ref{Thm:HM}.
\end{proof}

Lemma~\ref{Lem:Central form} is used in Section~\ref{Subsec:EmbeddingPart} to establish an embedding of a RAAG into an HHG, while Definition~\ref{Def:strictpartialorder} is crucial in Section~\ref{Subsec:UndistortionPart} for promoting this embedding to an undistorted one.

\section{Main embedding theorem}\label{Sec:SufficientConditions}

In this section, we introduce terminology for axial elements that facilitates their expected generation of a right-angled Artin subgroup in an HHG, and then present a refined version of Theorem~\ref{Thm:MainEmbeddingThm}. 

\subsection{Geometric irredundancy and support conditions for axial elements}\label{Subsec:Free subgroups}

Let $f$ and $g$ be axial elements in $G$.
If $\B(f)$ and $\B(g)$ are distinct, then $f$ and $g$ clearly cannot share a common power since $\B(f)=\B(f^n)$ for every nonzero integer $n$.
When $\B(f)=\B(g)$, however, an additional condition is required to ensure that $f$ and $g$ have no common power. The following fact shows that if at least one unbounded domain exhibits a form of `freeness', then sufficiently large powers of $f$ and $g$ generate an undistorted subgroup isomorphic to the free group of rank $2$; in particular $f$ and $g$ share no common power. 
This criterion is motivated by the work of Kent–Leininger \cite{KL08}, which establishes the prevalence of undistorted free subgroups in mapping class groups, highlighting the richness of their subgroup structure.

\begin{proposition}\label{Lem:KL08}
Let $(G,\mathfrak{S})$ be an HHG and $\{f_1,\dots,f_m\}$ a collection of axial elements in $G$.
Suppose that there exists a domain $U\in\bigcap_i\B(f_i)$ such that each $f_i$ acts loxodromically on $\mathcal{C}U$.
If the sets of limit points of the $f_i$ in $\upartial\mathcal{C}U$ are all disjoint, then there exists $N>0$ such that for any $n\ge N$, the $n$-th powers of the $f_i$ generate a quasi-isometrically embedded free subgroup of rank $m$.
\end{proposition}
\begin{proof}
Following the proof of \cite[Theorem 1.4]{KL08}, there exists a constant $N$ such that if $n\ge N$, then the abstract free group $\langle f^n_1,\dots,f^n_m\rangle$ is quasi-isometrically embedded in $\mathcal{C}U$ via the orbit map based at some (indeed any) point.
In particular, there exists $C=C(n)>0$ with $\lim_{n\to\infty} C=\infty$ such that for any $g\in G$ and any cyclically reduced word $w$ in $\{f^n_1,\dots,f^n_m\}$, we have $\dist_{U}(g,w.g)\ge C|w|$, where $|w|$ is the word length of $w$. 
It follows by Theorem~\ref{Thm:DistanceFormula} that (after increasing $N$ if necessary) there exists $D$, depending on $C$ and the HHS constant $E$, such that $\dist_{G}(g,w.g)\ge D|w|$. Therefore, the subgroup generated by $\{f^n_1,\dots,f^n_m\}$ is the desired free subgroup.
\end{proof}

The setting of \cite[Theorem 1.4]{KL08} is the case where $G$ is the mapping class group of a closed hyperbolic surface and each $f_i$ is a pseudo-Anosov mapping class in the above proposition.

Based on Proposition~\ref{Lem:KL08}, we introduce the following terminology, which provides a sufficient condition for finitely many axial elements of $G$ to have no pair of distinct elements sharing a common power.

\begin{definition}[Geometric irredundancy]\label{Def:irredundant}
Let $(G,\mathfrak{S})$ be an HHG. A finite collection of axial elements $f_i$ in $G$ is said to be \emph{geometrically irredundant} if, for any two distinct $f_i$ and $f_j$, at least one of the following holds:
\begin{itemize}
\item $\B(f_i)$ and $\B(f_j)$ are not identical, or
\item there is a domain $U\in \B(f_i)\cap \B(f_j)$ such that the sets of limit points of $f_i$ and $f_j$ in $\upartial\mathcal{C}U$ are disjoint. 
\end{itemize}
\end{definition}

Suppose that $f$ and $g$ are axial elements, and let $U\in\mathfrak{S}$ be a domain belonging to neither $\B(f)$ nor $\B(g)$. Assume further that both $f$ and $g$ stabilize $U$; in particular, any word $w$ in $\{f,g\}$ also stabilizes $U$.
By definition, the actions of $f$ and $g$ on $\mathcal{C}U$ are elliptic but the action of $w$ on $\mathcal{C}U$ need not be elliptic.
To begin addressing such situations, we introduce the following notion, analogous to a partial pseudo-Anosov mapping class supported on an essential subsurface or to a Dehn twist in a mapping class group.

\begin{definition}[Being fully supported and rigidly fully supported]\label{Def:FullySupportedElement}
Let $(G,\mathfrak{S})$ be an HHG. An axial element $g\in G$ is said to be \emph{fully supported} on an unbounded domain $U\in\mathfrak{S}$ if 
\begin{enumerate}
\item\label{Item:stabilize} $g$ stabilizes both $U$ and every unbounded domain $V$ orthogonal to $U$,
\item\label{Item:Loxodromic} the action of $g$ on $\mathcal{C}U$ is loxodromic, and
\item\label{Item:UniformlyBoundedOrbit} for each unbounded $V\bot U$, there exists a constant $Q_{g,V}\ge 0$ such that every $g$-orbit in $\mathcal{C}V$ has diameter at most $Q_{g,V}$.
\end{enumerate}
If, moreover, the action of $g$ on $\mathcal{C}V$ is trivial for every unbounded $V\bot U$ (equivalently, one may take $Q_{g,V}=0$ for every such $V$), then $g$ is said to be \emph{rigidly fully supported} on $U$.

When \(g\) is fully supported on \(U\), we refer to a quasi-axis for the action of \(g\) on \(\mathcal{C}U\), in the sense of Definition~\ref{Def:Quasi-axis}, simply as a \emph{quasi-axis of \(g\)}.
\end{definition}

\begin{remark}\label{Rmk:BigsetOfFullySupported}
If an axial element \(g\) is fully supported on an unbounded domain \(U\), then
\(\B(g)=\{U\}\).
Indeed, \(U\in\B(g)\) by the loxodromic action of \(g\) on \(\mathcal{C}U\).
Since the domains in \(\B(g)\) are pairwise orthogonal \cite[Lemma~6.7]{DHS17}, any other domain in \(\B(g)\) would be an unbounded domain orthogonal to \(U\), contradicting the bounded-orbit condition in Definition~\ref{Def:FullySupportedElement}.
\end{remark}

\begin{remark}\label{Rmk:UnboundedDomainsInFullSupport}
The conditions in Definition~\ref{Def:FullySupportedElement} are imposed only on unbounded domains orthogonal to \(U\). Indeed, only unbounded domains can belong to the bigset of an element, and all orthogonal-domain control used in the proof of the main embedding theorem concerns the unbounded supporting domains of the chosen elements.
Thus no condition on bounded domains orthogonal to \(U\) is needed for the arguments of this paper.
\end{remark}

\begin{definition}[Weak commutativity]\label{Def:WeakCommutativity}
We say that an HHG $(G,\mathfrak{S})$ satisfies the \emph{weak commutativity property} if, whenever $g,h\in G$ are axial elements fully supported on orthogonal unbounded domains, there exists some $N>0$ such that $g^{N}$ and $h^{N}$ commute. 
\end{definition}

The next lemma shows that, once axial elements fully supported on pairwise orthogonal domains commute, sufficiently large common powers generate an undistorted free abelian subgroup.
Moreover, the action of this subgroup on the hyperbolic space associated to any unbounded domain orthogonal to all of the supporting domains has uniformly bounded orbits.

\begin{lemma}\label{Lem:UndistortedFreeAbelianSubgroup}
Let $f_1,\dots,f_m \in G$ be commuting axial elements, each fully supported on an unbounded domain $U_i \in \mathfrak{S}$ with $U_i \bot U_j$ whenever $i \neq j$. Then there exists a constant $N$ such that if $n\ge N$, then:
\begin{enumerate}
\item The subgroup $H=\langle f_1^n,\dots,f_m^n\rangle$ is an undistorted free abelian subgroup of rank $m$.
\item For any unbounded domain $V$ with $V \bot U_i$ for all $i=1,\dots,m$, the action of every element of $H$ on $\mathcal{C}V$ has orbits of diameter at most $Q_V$, where $Q_V$ depends only on $\{Q_{f_i,V}\}_{i=1}^m$ and the HHS constant $E$.
\end{enumerate}
\end{lemma}

Note that $m$ is bounded by the HHS constant $E$ (see \cite[Lemma~2.1]{BHS19}).

\begin{proof}
Let $s_0$ and $K_{s_0}$ be the constants in Theorem~\ref{Thm:DistanceFormula}, and define
\[Q=\max\bigl(\{\,Q_{f_j,U_i}\mid 1\le i,j\le m,\ i\neq j\,\}\cup\{0\}\bigr).\]
Since the stable translation length $\tau_{\mathcal{C}U_i}(f_i)$ is positive, we can choose a sufficiently large constant $N>0$ such that, for any $n\ge N$, $\tau_{\mathcal{C}U_i}(f_i^n)\ge s_0+2E+mQ$ for each $i$. 

Let $w=f_1^{na_1}\cdots f_m^{na_m}$ be a word in $\{f_1^n,\dots,f_m^n\}$ with $a_i\in\Z$, and let $x\in G$.
By Definitions~\ref{Def:FullySupportedElement} and~\ref{Def:HHG}, we have 
\begin{equation}\label{Eq:CoarselyEqual}
\begin{aligned}
\pi_{U_i}(wx)&=f_1^{na_{1}}f_2^{na_{2}}\cdots f_m^{na_m}\pi_{U_i}(x)\\
&\stackrel{Q}{=}  f_1^{na_1}f_2^{na_2}\cdots f_{m-1}^{na_{m-1}}\pi_{U_i}(x)\stackrel{Q}{=}\cdots \\ 
&\stackrel{Q}{=} f_1^{na_1}f_2^{na_2}\cdots f_{i}^{na_{i}}\pi_{U_i}(x) = f_1^{na_1}f_2^{na_2}\cdots f_{i-1}^{na_{i-1}}\pi_{U_i}(f_{i}^{na_{i}}x)\\
&\stackrel{Q}{=} f_1^{na_1}f_2^{na_2}\cdots f_{i-2}^{na_{i-2}}\pi_{U_i}(f_{i}^{na_{i}}x) \stackrel{Q}{=}\cdots\stackrel{Q}{=} \pi_{U_i}(f^{na_i}_ix),    
\end{aligned}
\end{equation}
where $A \stackrel{Q}{=} B$ indicates that $\dist_{\mathrm{Haus}}(A,B) \le Q$.
By Lemma~\ref{Lem:UniformQuasiAxis}, it follows that  
\[\dist_{U_i}(x,wx)\asymp_{1,mQ}\dist_{U_i}(x,f_i^{na_i}x)\ge |a_i|\tau_{\mathcal{C}U_i}(f_i^n)-2E.\]
Applying the distance formula, we have 
\[\dist_G(x, wx)\ge \frac{1}{K_{s_0}}\!\left( \sum_{i=1}^{m}\Bigl\{\!\!\Bigl\{|a_i|\tau_{\mathcal{C}U_i}(f_i^n)-2E-mQ\Bigr\}\!\!\Bigr\}_{s_0} \!\right)-K_{s_0}\ge \frac{s_0}{K_{s_0}}\!\left(\sum_{i=1}^m|a_i|\!\right)-K_{s_0},\]
where $\sum_{i=1}^m|a_i|$ is equal to the word length of $w$.
Thus, the subgroup $H$ generated by $\{f_1^n,\dots,f_m^n\}$ is quasi-isometrically embedded.

Suppose that $V\in\mathfrak{S}$ is an unbounded domain orthogonal to each $U_i$.
Set \[q_{V}=\max\!\left\{\,Q_{f_i,V}\mid 1\le i\le m\,\right\}.\]
By the same reasoning as in Equation~\eqref{Eq:CoarselyEqual}, for any $h\in H$ and $x\in G$, we have \[\pi_{V}(hx)=h\pi_{V}(x)\stackrel{Q_V}{=}\pi_{V}(x),\]
where $Q_V=m q_V$. Therefore, the action of each $h\in H$ on $\mathcal{C}V$ has uniformly bounded orbits, with diameter independent of the choice of $h$.
\end{proof}

\subsection{Precise statement of the main embedding theorem}\label{Subsec:Homomorphism}

We now give the precise version of Theorem~\ref{Thm:MainEmbeddingThm}.
The first condition below ensures that suitable powers of the chosen elements define a homomorphism from the RAAG determined by their supporting domains, while geometric irredundancy allows us to prove that this homomorphism is injective.
The second condition is the bounded-orbit hypothesis needed for undistortion; as suggested by Lemma~\ref{Lem:UndistortedFreeAbelianSubgroup}, this condition is automatic in certain pairwise orthogonal situations, but it need not hold in general.

\begin{theorem}\label{Thm:MainEmbeddingThm_Precise}
Let $(G,\mathfrak{S})$ be a hierarchically hyperbolic group.
For each $i=1,\dots,m$, let $f_i\in G$ be an axial element fully supported on an unbounded domain $U_i\in\mathfrak{S}$, and suppose that the collection $\{f_1,\dots,f_m\}$ is geometrically irredundant.

Consider the following two conditions: 
\begin{enumerate}[label=(A.\arabic*), ref=A.\arabic*]
\item\label{Item:Sec3PowersCommute} 
There exists $N>0$ such that, for every positive integer \(d\), there is a homomorphism
\begin{equation}\label{Eq:Homomorphism}
\phi_{d N}:A(\mathcal{O}^{\Sigma})\to G \quad\text{defined by}\quad v_i\mapsto f^{d N}_i,
\end{equation}
where $\mathcal{O}^{\Sigma}$ is the orthogonality graph of $\Sigma=\{U_1,\dots,U_m\}$ (Definition~\ref{Def:OrthogonalityGraph}) and $v_i$ is the standard generator of $A(\mathcal{O}^{\Sigma})$ corresponding to $U_i$.
\item\label{Item:Sec3SubgroupHasUniformlyBoundedOrbits} 
There exist $N_1>0$ and $Q\ge 0$ such that, for each $i$, the action of any element in the subgroup $\langle \{f^{N_1}_j\mid U_j\bot U_i\}\rangle$ on $\mathcal{C}U_i$ has orbits of diameter at most $Q$.
\end{enumerate}

If Condition~\eqref{Item:Sec3PowersCommute} holds, then there exists a constant $D = D(\{ f_i \})>0$ such that the homomorphism~\eqref{Eq:Homomorphism} is injective for any $d\ge D$. 

If, in addition, Condition~\eqref{Item:Sec3SubgroupHasUniformlyBoundedOrbits} holds, then after replacing $N$ by $\lcm(N,N_1)$ and possibly increasing $D$, the homomorphism~\eqref{Eq:Homomorphism} is moreover a quasi-isometric embedding.
\end{theorem}

\begin{remark}
If each $f_i$ is \emph{rigidly} fully supported on $U_i$, the statement simplifies, as the two technical conditions in Theorem~\ref{Thm:MainEmbeddingThm_Precise} become automatic. This is the content of Theorem~\ref{Thm:MainTheoremStrong}, and a precise formulation with proof will be given in Section~\ref{Subsec:AssumptionsOnElements}.
\end{remark}

\begin{remark}[Variants of conditions]
Condition~\eqref{Item:Sec3PowersCommute} in Theorem~\ref{Thm:MainEmbeddingThm_Precise} can be omitted if the HHG $(G,\mathfrak{S})$ satisfies the weak commutativity property.
One might hope to replace \emph{geometric irredundancy} by the weaker notion of \emph{irredundancy}, meaning that no pair of distinct elements shares a common nonzero power. For a general action on a hyperbolic space, however, irredundancy does not imply the quasi-axis separation used in our proof. 
Indeed, two translations of the same line whose translation lengths have irrational ratio commute and share no common nonzero power, while their quasi-axes coincide. Thus irredundancy alone is insufficient to guarantee the geometric separation required by our argument.

If, however, the HHG \((G,\mathfrak S)\) is \emph{hierarchically acylindrical} in the sense of \cite[Definition~9.18]{DHS17}, then geometric irredundancy can be replaced by irredundancy. Indeed, for each domain \(U\in\mathfrak S\), the image of \(\operatorname{Stab}_G(U)\) in \(\operatorname{Isom}(\mathcal CU)\) acts acylindrically. Hence two loxodromic elements acting on \(\mathcal CU\) whose limit sets are not disjoint are commensurable. Consequently, for axial elements fully supported on the same domain, irredundancy implies that their limit sets are disjoint, and therefore implies geometric irredundancy.
\end{remark}

Theorem~\ref{Thm:MainEmbeddingThm_Precise} naturally divides into two parts: the \emph{embedding} part (with the first condition) and the \emph{undistortion} part (with both conditions). The proofs of these parts are given in Section~\ref{Sec:Main Proof}, in that order.

\section{Projections onto quasi-geodesics}\label{Sec:Projections and MSBI}

In an HHS \((\mathcal{X},\mathfrak{S})\), the projection maps satisfy
consistency inequalities and coarsely commute with automorphisms.
In this section, we show that analogous properties hold after postcomposing
these projections with nearest-point projections onto uniform
quasi-geodesics in the associated hyperbolic spaces.

Let \(g\) be an automorphism of \((\mathcal{X},\mathfrak{S})\), and let
\(\alpha\subseteq\mathcal{C}U\) be a quasi-geodesic for
\(U\in\mathfrak{S}\).
Since \(g\) induces an isometry
\(\mathcal{C}U\to\mathcal{C}(gU)\), the maps
$g\circ\mathfrak p_\alpha\circ\pi_U$ and $\mathfrak p_{g\alpha}\circ\pi_{gU}\circ g$ coarsely coincide.
For HHGs, we will use the following special case.

\begin{lemma}\label{lemma:commuting projection maps}
For an HHG \((G,\mathfrak{S})\) with constant $E$, let \(g\in G\) stabilize a domain
\(U\in\mathfrak{S}\), and let
\(\alpha\subseteq\mathcal{C}U\) be a \(g\)-invariant $\lambda$-quasi-axis.
Then, for every \(n\in\mathbb{Z}\) and every \(x\in G\),
\[
\dist_{\operatorname{Haus}}\!\left(
g^n\circ\mathfrak p_\alpha\circ\pi_U(x),\,
\mathfrak p_\alpha\circ\pi_U(g^n x)
\right)
\le R,
\]
where \(R=R(E,\lambda)\) is the constant from Lemma~\ref{Lem:CoarselyLipschitz},
applied to \(\alpha\).
\end{lemma}

\begin{proof}
Since \(g\alpha=\alpha\), we have \(g^n\alpha=\alpha\) for every \(n\in\mathbb{Z}\). Moreover, \(g^n\) acts isometrically on \(\mathcal{C}U\), and the equivariance of the HHG projections gives $g^n\pi_U(x)=\pi_U(g^nx)$.
Thus, for each \(p\in\pi_U(x)\), both $g^n\mathfrak p_\alpha(p)$ and $\mathfrak p_\alpha(g^np)$ are nearest-point projections of \(g^np\) to \(\alpha\), which are \(R\)-close by Lemma~\ref{Lem:CoarselyLipschitz}.
Therefore, the Hausdorff-distance conclusion follows.
\end{proof}

The consistency inequalities also remain valid after postcomposing with nearest-point projections onto quasi-geodesics; in the mapping class group case, this was called the \emph{multi-scale Behrstock inequality} \cite[Lemma~5]{Run21}.
However, postcomposition requires introducing a sufficiently large constant, determined by those in Section~\ref{Sec:ActionsonHyperbolicSpaces}.

Using the HHS constant \(E\), fix a constant $E'\ge L+8E$, where $L= L(E)$ is the constant given in Lemma~\ref{Lem:UniformQuasiAxis}.
Then choose 
\[ E_0\ge 2\max\{R(E,E'),\,E_1(E,E'),\,E_2(E,E'),\,E_3(E,E')\}\] 
for all constants and additive errors arising from Lemmas~\ref{Lem:CoarselyLipschitz}, \ref{Lem:Composition of Projection}, and \ref{Lem:Projections of two quasiconvex}.

\begin{notation}\label{Notation:OntoQuasigeodesics}
When measuring distances between images under nearest-point projections onto quasi-geodesics, we often omit the projection notation in line with Notation~\ref{Notation:Omitting Projection maps}.
For instance, let $\alpha\subseteq \mathcal{C}U$ be a quasi-geodesic.
Given $x,y\in \mathcal{X}$ and $\gamma\subseteq\mathcal{C}V$, we write 
$$\dist_{\alpha}(x,y)\ \text{ for }\ \dist_{U}(\mathfrak{p}_{\alpha}\circ\pi_{U}(x),\,\mathfrak{p}_{\alpha}\circ\pi_{U}(y))\quad\text{and}\quad\dist_{\alpha}(x,\gamma)\ \text{ for }\ \dist_U(\mathfrak{p}_{\alpha}\circ\pi_U(x),\,\mathfrak{p}_{\alpha}\circ\rho_U^V(\gamma)).$$
If $U=V$, then $\rho_U^V$ is the identity on $\mathcal{C}V$.
\end{notation}

\begin{lemma}[Multi-scale consistency inequality]\label{Lem:MSBI}
Let \((\mathcal{X},\mathfrak{S})\) be an HHS with HHS constant \(E\), and let \(E', E_0\) be the constants fixed above.
For non-orthogonal domains \(U,V\in\mathfrak{S}\), let \(\alpha\subseteq\mathcal{C}U\) and \(\beta\subseteq\mathcal{C}V\) be \((2,E')\)-quasi-geodesics.
Then, for every \(x\in\mathcal{X}\),
\begin{equation}\label{Inequality}
\min \{\dist_{\alpha}(x,\beta),\, \dist_{\beta}(x,\alpha)\}\le2E_0. 
\end{equation}
\end{lemma}
\begin{proof}
Since the conclusion is symmetric in \((U,\alpha)\) and \((V,\beta)\), it suffices to consider the following three cases: \(U=V\), \(U\pitchfork V\), and \(U\sqsubsetneq V\).

\smallskip
\noindent{\textbf{Case 1 (\(U=V\)).}}
Suppose that \(\dist_{\alpha}(x,\beta)>E_0\).
Choose \(z\in\pi_U(x)\), and let \(\gamma\) be a geodesic joining \(z\) to \(\mathfrak p_\beta(z)\).
Since $\diam_{\mathcal{C}U}(\mathfrak{p}_{\alpha}(\gamma))$ must be greater than $E_0$, Lemma~\ref{Lem:Projections of two quasiconvex} implies that $\mathcal{N}_{E_{0}}(\gamma)\cap \alpha\neq \emptyset$.
Choose \(a\in\alpha\) and \(r\in\gamma\) such that \(\dist_U(a,r)\le E_0\).
Since \(\gamma\) joins \(z\) to \(\mathfrak p_\beta(z)\), the point \(\mathfrak p_\beta(z)\) is $R$-close to \(\mathfrak{p}_\beta(r)\).
Hence, by Lemma~\ref{Lem:CoarselyLipschitz}, it follows that $\dist_\beta(x,\alpha)\le 2E_0$.

\smallskip
\noindent{\textbf{Case 2 ($U\pitchfork V$).}}
In this case, by the transversality and consistency axiom, we have 
$$\min \bigl\{\, \dist_{U}(x,\rho_U^V),\, \dist_{V}(x,\rho_V^U)\,\bigr\}<E.$$
Since $\mathfrak{p}_{\alpha}$ and $\mathfrak{p}_{\beta}$ are $(1,E_0)$-coarsely Lipschitz, Inequality~\eqref{Inequality} holds.

\smallskip
\noindent{\textbf{Case 3 ($U\sqsubsetneq V$).}}
Suppose that $\dist_{\alpha}(x,\beta)>2E_0$. 
If $\pi_{V}(x)\cap \mathcal{N}_{E}(\rho^{U}_{V})\neq\emptyset$, then 
$\dist_\beta(x,\alpha)\le2E_0$.
On the other hand, if $\dist_{V}(x,\rho^{U}_{V})>E$, by the consistency axiom, we have $$\diam_{\mathcal{C}U}\bigl(\rho_U^V(\pi_V(x))\cup\pi_U(x)\bigr)\le E.$$
Choose \(z\in\pi_V(x)\), and let \(\gamma\) be a geodesic joining
\(z\) to \(\mathfrak p_\beta(z)\). 
If $\gamma$ intersects $\mathcal{N}_{E}(\rho^{U}_{V})$, then $\dist_{\beta}(x,\alpha)\le E_0$.
Otherwise, the bounded geodesic image axiom gives $\diam_{\mathcal C U}\bigl(\rho_U^V(\gamma)\bigr)\le E$.
Since \(z\in\pi_V(x)\), \(q=\mathfrak p_\beta(z)\in\beta\), and $\diam_{\mathcal C U}\bigl(\rho_U^V(\pi_V(x))\cup\pi_U(x)\bigr)\le E$, it follows that
$\dist_\alpha(x,\beta)\le 2E_0$, contrary to our assumption.
\end{proof}

\begin{remark}\label{Rmk:Bounded}
Let $U,V\in\mathfrak{S}$ be non-orthogonal domains, and let $\alpha\subseteq\mathcal{C}U$ and $\beta\subseteq\mathcal{C}V$ be quasi-geodesics. 
Then the projection $\mathfrak{p}_{\beta}\circ\rho^{U}_{V}(\alpha)$ from $\alpha$ to $\beta$ has finite diameter unless $U=V$ and $\alpha,\beta$ share a limit point in $\upartial\mathcal{C}U$.

Indeed, if $U\pitchfork V$ or $U\sqsubsetneq V$, then $\diam_{\mathcal{C}V}(\rho^{U}_{V})<E$, hence $\diam_{\mathcal{C}V}(\mathfrak{p}_{\beta}\circ\rho^{U}_{V}(\alpha))< 2E_0$.
If $V\sqsubsetneq U$, then the bounded geodesic image axiom implies that each connected component of $\alpha-\mathcal{N}_{E}(\rho^{V}_{U})$ (recall Remark~\ref{Rmk:RestrictedDomain}) has image of diameter at most $E$ under $\rho^{U}_{V}$, so $\mathfrak{p}_{\beta}\circ\rho^{U}_{V}(\alpha)$ has finite diameter.  
If $U=V$, then $\mathfrak{p}_{\beta}\circ\rho^{U}_{V}(\alpha)=\mathfrak{p}_{\beta}(\alpha)$, which has finite diameter exactly when $\alpha$ and $\beta$ have disjoint limit sets in $\upartial\mathcal{C}U$.
\end{remark}

Motivated by the notion of relevant domains (Definition~\ref{Def:RelevantDomain}), we define \emph{relevant quasi-geodesics}, which will play a key role in the next section.

\begin{definition}[Relevant quasi-geodesic]\label{Def:RelevantQuasigeodesic}
Let $x,y\in\mathcal{X}$ and $\theta>0$.  
A quasi-geodesic $\gamma\subseteq\mathcal{C}U$ for $U\in\mathfrak{S}$ is said to be \emph{$\theta$-relevant (with respect to $x,y$)} if $\dist_{\gamma}(x,y)>\theta$. 
\end{definition}

\section{Proof of Theorem~\ref{Thm:MainEmbeddingThm_Precise}}\label{Sec:Main Proof}

Throughout this section, we make the following assumptions:
\begin{itemize}[leftmargin=2em]
\item \((G,\mathfrak{S})\) is an HHG with HHS constant \(E\), and, for each \(i\in\mathcal{I}=\{1,\dots,m\}\), the element \(f_i\in G\) is axial and fully supported on an unbounded domain \(Z_i\in\mathfrak{S}\).
\item We set $\Sigma=\{Z_i\}_{i\in\mathcal{I}}$, allowing repetitions, and let \(\Gamma=\mathcal{O}^{\Sigma}\) be its orthogonality graph.
\item The collection \(\{f_1,\dots,f_m\}\) is geometrically irredundant, and Condition~\eqref{Item:Sec3PowersCommute} holds.
\item Unless otherwise stated, we fix a point \(x\in G\).
\end{itemize}

Let \(N>0\) be the constant given by Condition~\eqref{Item:Sec3PowersCommute}.
Thus, for every positive integer \(d\), there is a homomorphism
\begin{equation}\label{Eq:HomomorphismAgain}
\phi=\phi_{dN}:A(\Gamma)\longrightarrow H\leq G
\qquad\text{defined by}\qquad v_i\longmapsto f_i^{dN},
\end{equation}
where each vertex \(v_i\in\Gamma\) corresponds to \(Z_i\in\Sigma\).
For each \(i\in\mathcal{I}\), define
\[J'(i)=\{\,j\in\mathcal{I}\mid Z_j\not\bot Z_i \text{ and } j\neq i\,\}.\]
Equivalently, \(J'(i)\) is the set of indices \(j\neq i\) such that the standard generators \(v_j\) and \(v_i\) do not commute in \(A(\Gamma)\).

We first record the bounded-orbit constants needed below. Let \(\triangle\) be a clique in \(\Gamma\), and let \(v_i\in\triangle\).
Since the elements $\{f_j^N\mid v_j\in\triangle\setminus\{v_i\}\}$ pairwise commute, Lemma~\ref{Lem:UndistortedFreeAbelianSubgroup} gives constants \(D_{\triangle,i}>0\) and \(Q_{\triangle,i}\geq 0\) such that, for every \(d\geq D_{\triangle,i}\), every element of $\langle\, f_j^{dN}\mid v_j\in\triangle\setminus\{v_i\}\,\rangle$ acts on \(\mathcal{C}Z_i\) with orbits of diameter at most \(Q_{\triangle,i}\).
Set
\[
D_0=\max\{\,D_{\triangle,i}\mid v_i\in\triangle,\, 
\triangle\text{ a clique in }\Gamma\,\}
\quad\text{and}\quad
Q=\max\{\,Q_{\triangle,i}\mid v_i\in\triangle,\, 
\triangle\text{ a clique in }\Gamma\,\}.
\]

We next record a uniform bound associated to quasi-axes.
For each \(i\in\mathcal{I}\), every sufficiently large positive power of \(f_i\) has the same pair of limit points in \(\upartial\mathcal{C}Z_i\), and any quasi-axis obtained from Lemma~\ref{Lem:UniformQuasiAxis} is a uniform quasi-geodesic joining these limit points.
Hence, by Remark~\ref{Rmk:Bounded}, geometric irredundancy, and the Morse lemma, there exists a constant \(B\geq 0\), independent of \(d\), such that whenever \(d\) is large enough, say $d\ge D_1$, for every \(f_i^{dN}\) to admit a quasi-axis \(\alpha_i'\subseteq\mathcal{C}Z_i\), we have
\begin{equation}\label{Inequality:B}
\diam_{\mathcal{C}Z_i}\left(\mathfrak{p}_{\alpha_i'}(\pi_{Z_i}(x)) \cup \bigcup_{j\in J'(i)} \mathfrak{p}_{\alpha_i'} \bigl(\rho^{Z_j}_{Z_i}(\alpha_j')\bigr)\right)\leq B
\end{equation}
for every \(i\in\mathcal{I}\).

Let \(E'_1=E'_1(E,E')\) be the constant from Lemma~\ref{Lem:Projection of quasiconvex subsets} with $E'\ge L+8E$ ($L$ is given in Lemma~\ref{Lem:UniformQuasiAxis}), enlarged if necessary so that \(E'_1\ge E'\); in particular, \(E'_1\) bounds the quasi-convexity constants of all quasi-axes and all projection images considered in this section.
We now fix a sufficiently large constant \(\Theta\) dominating the following quantities:
\begin{itemize}[leftmargin=2em]
\item $\max\{s_0,\,2E_0,\,Q,\,B,\,1\}$, where \(s_0\) is the constant from Theorem~\ref{Thm:DistanceFormula} and \(E_0\) is the constant from Lemma~\ref{Lem:MSBI};
\item $2\max\{R(E,E'_1),\,E_1(E,E'_1),\,E_2(E,E'_1),\,E_3(E,E'_1)\}$, where \(R,E_1,E_2,E_3\) are the constants from Lemmas~\ref{Lem:CoarselyLipschitz}, \ref{Lem:Composition of Projection}, and \ref{Lem:Projections of two quasiconvex}.
\end{itemize}

Finally, choose \(D\geq \max\{D_0,D_1\}\) sufficiently large so that, for every \(d\geq D\) and every \(i\in\mathcal{I}\), $\tau_{\mathcal{C}Z_i}(f_i^{dN})>120\Theta$.

For the proof, fix \(d\geq D\) and choose a quasi-axis $\alpha_i'\subseteq\mathcal{C}Z_i$ of \(f_i^{dN}\) for each \(i\in\mathcal{I}\). Since \(\Theta\geq 2E_0\), Lemma~\ref{Lem:MSBI} applies throughout this section with \(\Theta\) in place of \(2E_0\). Moreover, Inequality~\eqref{Inequality:B} holds after replacing $B$ by $\Theta$.

\smallskip

As recalled in Section~\ref{Subsec:Normal forms}, every nontrivial element of $A(\Gamma)$ admits a representative with the fewest syllables, unique up to shuffling. 
By appropriately choosing shufflings, such a representative can be expressed in a central form. In the first subsection, devoted to the proof of the embedding part of Theorem~\ref{Thm:MainEmbeddingThm_Precise}, the central form falls within the scope of Lemma~\ref{Lem:UndistortedFreeAbelianSubgroup}, yielding that its image under $\phi$ is nontrivial. 
In the second subsection, concerning the proof of the undistortion part of Theorem~\ref{Thm:MainEmbeddingThm_Precise}, it is necessary to track every syllable of the element, and we also establish the necessity of this bookkeeping.

\subsection{Embedding part}\label{Subsec:EmbeddingPart}
As indicated above, by applying a central form of an element of $A(\Gamma)$ (see Definition~\ref{Def:Central form}) together with Lemma~\ref{Lem:UndistortedFreeAbelianSubgroup}, we prove the first part of Theorem~\ref{Thm:MainEmbeddingThm_Precise}.

\begin{proof}[Proof of the embedding part of Theorem~\ref{Thm:MainEmbeddingThm_Precise}]
By Lemma~\ref{Lem:Central form}, any nontrivial element in $A(\Gamma)$ admits a central form $w=u_{n}h_{n}\cdots u_{1}h_{1}$ such that
\begin{itemize}
\item For $1\le j\le n$, $h_j$ is a power of a standard generator $v_{i_j}$ and $\operatorname{supp}(u_j h_{j})$ spans a clique in $\Gamma$, and
\item $h_{j}$ does not commute with $h_{{j-1}}$ or, equivalently, $i_{j}\in J'(i_{j-1})$ for $2\le j\le n$. 
\end{itemize}
To complete the proof, it suffices to show that $\phi(w)\in G$ is nontrivial. 
If $n=1$, i.e., $w=u_1h_1$, then Lemma~\ref{Lem:UndistortedFreeAbelianSubgroup} ensures that $\phi(w)$ is nontrivial, and we are done. Thus we may assume $n\ge 2$. Since conjugation does not affect the triviality of $\phi(w)$, we may, after conjugating $w$ by $u_1h_1$ if necessary, further assume that $h_{1}$ does not commute with $h_{n}$.

For convenience, for each $1\le j\le n$, set 
\[U_j := Z_{i_j},\quad J(j):=\{ k \mid 1\le k\le n\text{ and $h_j$ and $h_k$ do not commute} \}\quad\text{and}\quad \alpha_j := \alpha'_{i_j}.\] 
(Note that $J(j)$ is non-empty by the definition of a central form and the assumption that $n\ge 2$.)
With this notation, \(\phi(h_j)\) is fully supported on \(U_j\), preserves \(\alpha_j\), and acts loxodromically on it with sufficiently large stable translation length.

We now prepare a ping-pong argument. For each $1\leq i \leq n$, define
\begin{align*}
X_{i}&=\left\{  y\in G \mid \min\left\{ \dist_{\alpha_{i}}(y,\alpha_{i'}) \mid i'\in J(i)\right\}>2\Theta  \right\}.      
\end{align*}

\begin{step}\label{Step:1}
Each $X_i$ is non-empty.    
\end{step} 
For each $1\le i\le n$ and ${i'}\in J(i)$, either $U_i$ and $U_{i'}$ are distinct and $U_i\notbot U_{i'}$, or $U_i=U_{i'}$ but the sets of endpoints of $\alpha_i$ and $\alpha_{i'}$ in $\upartial \mathcal{C}U_i$ are disjoint. 
Then by Remark~\ref{Rmk:Bounded} and our choice of $\Theta$, we have $$\diam_{\mathcal{C}U_i}\left(\bigcup_{i'\in J(i)} \mathfrak{p}_{\alpha_{i}}\left(\rho^{U_{i'}}_{U_{i}}(\alpha_{i'})\right)\right)<\Theta.$$
Since $\alpha_i$ has infinite diameter in $\mathcal{C}U_i$, we can choose $p_i\in \alpha_i$ so that
\[\dist_{\alpha_{i}}(p_{i},\alpha_{i'})>4\Theta\quad \text{for all }i'\in J(i).\] 
Let $y\in G$ be a partial realization of the singleton coordinate $\{p_i\}$. Then, by our choice of $\Theta$,
$$\dist_{U_{i}}(y,\alpha_{i'})\geq \dist_{U_{i}}(p_{i},\alpha_{i'})-\dist_{U_{i}}(p_{i},y)-\diam_{\mathcal{C}U_i}(\pi_{U_i}(y))>2\Theta,$$
for each $i'\in J(i)$ and hence $y\in X_{i}$. 

\begin{step}\label{Step:2}
If $i'\in J(i)$, then $X_i \cap X_{i'}=\emptyset$.    
\end{step} 

For any $y\in X_{i}$, since $\dist_{\alpha_{i}}(y,\alpha_{i'})>2\Theta$, by Lemma~\ref{Lem:MSBI}, $\dist_{\alpha_{i'}}(y,\alpha_{i})<\Theta$ and thus $y\notin X_{i'}$.

\begin{step}\label{Step:3}
For each $2\le i\le n$, $\phi(u_i h_i)$ maps $X_{i-1}$ into $X_i$.
\end{step}

First, by Lemma~\ref{Lem:UndistortedFreeAbelianSubgroup} and our choice of $\Theta$, we have \(\dist_{U_{i}}(x,\phi(u_{i}) x)<\Theta\) for all \(x\in G\), since \(\operatorname{supp}(u_i h_i)\) spans a clique and the standard generator occurring in \(h_i\) does not occur in \(u_i\).
Now let $y\in X_{i-1}$. Since $i-1\in J(i)$ (and thus $i\in J(i-1)$), we obtain
\[\dist_{\alpha_{i-1}}(y,\alpha_{i})>2\Theta \quad \implies \quad \dist_{\alpha_{i}}(y,\alpha_{i-1})<\Theta\] 
by Lemma~\ref{Lem:MSBI}. 
Since the stable translation length of $\phi(h_i)$ on $\alpha_i$ exceeds $120\Theta$, for each $k\in J(i)$, by Lemma~\ref{Lem:UniformQuasiAxis}, we have
\begin{align*}
\dist_{\alpha_i}\bigl(\phi(u_i h_i)y,\,\alpha_k\bigr)
&\ge \dist_{\alpha_i}\bigl(\phi(h_i)y,\,y\bigr) -\dist_{\alpha_i}\bigl(\phi(h_i)y,\,\phi(u_ih_i)y\bigr) - \dist_{\alpha_i}(y,\,\alpha_{i-1}) \\
&\quad\quad - \diam_{\mathcal{C}U_i}\Bigl(\mathfrak{p}_{\alpha_i}(\rho^{U_{i-1}}_{U_i}(\alpha_{i-1})) \cup \mathfrak{p}_{\alpha_i}(\rho^{U_k}_{U_i}(\alpha_k))\Bigr)-2\Theta \\
&> \dist_{\alpha_i}(\phi(h_i)y,\,y) -5\Theta >2\Theta.
\end{align*}
Hence, $\phi(u_i h_i)(y)\in X_i$, as required.

\begin{step} 
Run the ping-pong argument to complete the proof.
\end{step}
Choose $p\in \alpha_1$ such that $5\Theta<\dist_{\alpha_{1}}(p,\alpha_{k})<10\Theta$ for each $k\in J(1)$, which is possible by our choice of $\Theta$. Then a partial realization $y\in G$ of $\{p\}$ lies in $X_1$ by the argument in Step~\ref{Step:1}. For each $k\in J(1)$, we have
\begin{align*}
\dist_{\alpha_1}\bigl(\phi(u_1 h_1)y,\,\alpha_k\bigr)
&\ge \dist_{\alpha_1}\bigl(\phi(h_1)y,\,y\bigr) -\dist_{\alpha_1}\bigl(\phi(h_1)y,\,\phi(u_1h_1)y\bigr) -\dist_{\alpha_1}(y,\,\alpha_{k})-2\Theta \\
&> \dist_{\alpha_1}(\phi(h_1)y,\,y) -\Theta-(10\Theta+\Theta) -2\Theta >2\Theta,
\end{align*}
and thus $\phi(u_1 h_1)(y)\in X_1$.
By Step~\ref{Step:3}, successive application of $\phi(u_i h_i)$ sends 
\[\phi(u_{i-1}h_{i-1}\cdots u_1 h_1)y \in X_{i-1}\quad\text{to}\quad\phi(u_i h_i\cdots u_1 h_1)y \in X_i\] 
for each $2\le i\le n$ and thus we have $\phi(w)y \in X_n$.
Finally, by Step~\ref{Step:2} and the assumption that $h_{1}$ and $h_{n}$ do not commute, it follows that $X_1\cap X_n=\emptyset$ and hence $\phi(w)y\neq y$. Therefore $\phi(w)$ is a nontrivial element in $G$, completing the proof.
\end{proof}


\subsection{Undistortion part}\label{Subsec:UndistortionPart}
The previous subsection shows that the homomorphism~\eqref{Eq:HomomorphismAgain} is injective.
In this subsection, we further assume that Condition~\eqref{Item:Sec3SubgroupHasUniformlyBoundedOrbits} holds.
As a consequence, the constant $N$ in the homomorphism~\eqref{Eq:HomomorphismAgain} is replaced by $\lcm(N,N_1)$, and we assume that the constant $Q$ appearing in the construction of $\Theta$ is replaced by the constant $Q$ appearing in Condition~\eqref{Item:Sec3SubgroupHasUniformlyBoundedOrbits}.
The parameter $d$, introduced in the first paragraph of this section, remains sufficient once it is redefined using this updated value of $\Theta$.

Throughout this subsection, we fix a nontrivial element $y\in A(\Gamma)$.
For a word $w=x^{e_1}_1\cdots x^{e_k}_k$ in $\Min(y)$, where $x_j=v_{i_j}$ is a standard generator of $A(\Gamma)$, consider $$\mathbf{w}=\phi(x^{e_1}_1)\cdots\phi(x^{e_k}_k)$$ as a word in $\{f^{dN}_1,\dots,f^{dN}_m\}$. For convenience, we write $\mathbf{w}=g_1\cdots g_k$, where $$g_j=\phi(x^{e_j}_j)=(f^{dN}_{i_j})^{e_{j}},\qquad j=1,\dots,k.$$ 




\begin{convention}
With a slight abuse of terminology, we adopt the following conventions, noting that $\phi$ is now injective:
\begin{enumerate}
\item\label{Enum:GammaCommute} Two elements $h_1,h_2\in H$ are said to \emph{$\Gamma$-commute} if there exist $a_1,a_2\in A(\Gamma)$ with $\phi(a_1) = h_1$ and $\phi(a_2) = h_2$ such that $a_1$ and $a_2$ commute in $A(\Gamma)$. 
In particular, \((f_i^{dN})^{e_i}\) and \((f_{i'}^{dN})^{e_{i'}}\) \(\Gamma\)-commute if and only if either \(Z_i\bot Z_{i'}\) or \(i=i'\).

\item\label{Enum:Syllablelength}
The \emph{$\Gamma$-syllable length} of $\mathbf{w}$ is defined as the syllable length of $w$ with respect to the standard generating set of $A(\Gamma)$.
Moreover, we write $g_i \prec_s g_j$ if $x_i^{e_i} \prec_s x_j^{e_j}$, where $\prec_s$ is a strict partial order on $\syl(y)= \{x_1^{e_1}, \dots, x_k^{e_k}\}$ as introduced in Definition~\ref{Def:strictpartialorder}.
\end{enumerate}
\end{convention}

\begin{notation}\label{Notation:Section5}
We change notation by setting $$U_j=Z_{i_j}\qquad\text{and}\qquad\alpha_j=\alpha'_{i_j},$$ that is, $g_j$ is fully supported on $U_j$ and $\alpha_j\subseteq \mathcal{C}U_j$ is the quasi-axis of $g_j$.
We denote by $\mathbf{w}_{\ell}=g_{1}\cdots g_{\ell}$ the subword of $\mathbf{w}$ consisting of its first $\ell$ syllables, and let $\mathbf{w}_0$ denote the empty word.
We then define the \emph{translated} supporting domain $V_j=\mathbf{w}_{j-1}U_{{j}}$ and the \emph{translated} quasi-axis $\beta_{j}=\mathbf{w}_{j-1}\alpha_{j}$ for $j=1,\dots,k$, and set
$$\mathcal{D}_{\mathbf{w}}=\{V_{j}\mid j=1,\dots,k\}\subseteq\mathfrak{S}\quad\text{and}\quad\mathcal{Q}_{\mathbf{w}}=\{\beta_{j}\mid j=1,\dots,k\}.$$
We emphasize that $\mathcal{D}_{\mathbf{w}}$ is considered as a \emph{subset}, so repetitions are not allowed.
\end{notation}

The following lemma is \cite[Lemma~5.1]{CLM12}, reformulated in our notation. We include its short proof for completeness. 

\begin{lemma}\label{Lem:TranslatedDomainsInvariantUnderShuffling} Let \(w,w'\in\Min(y)\) differ by a sequence of shufflings. Under the canonical identification \(\syl(w)\leftrightarrow\syl(w')\) described in Section~\ref{Subsec:Normal forms}, corresponding syllables determine the same translated supporting domain. In particular, \(\mathcal{D}_{\mathbf{w}}\) is invariant under shufflings. 
\end{lemma}
\begin{proof}
It suffices to consider a single shuffling. Write
\[w=a\,x_r^{e_r}x_{r+1}^{e_{r+1}}\,b\quad\text{and}\quad w'=a\,x_{r+1}^{e_{r+1}}x_r^{e_r}\,b,\]
where \(x_r\) and \(x_{r+1}\) commute in \(A(\Gamma)\). Let
\[\mathbf a=\phi(a),\qquad g_r=\phi(x_r^{e_r}),\qquad g_{r+1}=\phi(x_{r+1}^{e_{r+1}}). \]
Since \(w\in\Min(y)\), \(x_r\neq x_{r+1}\).
Hence the vertices \(x_r\) and \(x_{r+1}\) are adjacent in \(\Gamma\), so their supporting domains \(U_r\) and \(U_{r+1}\) are orthogonal.
Since \(g_r\) and \(g_{r+1}\) are fully supported on \(U_r\) and \(U_{r+1}\), respectively, we have
\[g_rU_{r+1}=U_{r+1}\qquad\text{and}\qquad g_{r+1}U_r=U_r. \]
Therefore, the translated supporting domains associated to the two shuffled syllables agree before and after the shuffling:
\[\mathbf a U_r=\mathbf a g_{r+1}U_r, \qquad \mathbf a g_rU_{r+1}=\mathbf a U_{r+1}.\]
The translated supporting domains associated to syllables preceding the shuffled pair are unaffected.
For every syllable following the shuffled pair, the two corresponding
prefixes have the same image under \(\phi\), since $g_rg_{r+1}=g_{r+1}g_r$.
Hence its translated supporting domain is also unchanged.
\end{proof}

We briefly explain how the translated quasi-axes introduced above are used to prove undistortion.
Lemma~\ref{Lemma:Relative quasigeodesic} first shows that each syllable \(g_i\) makes progress, linear in \(\lvert e_i\rvert\), along its translated quasi-axis \(\beta_i\).
If the translated supporting domains \(V_i\) were all distinct, the distance formula would immediately give the desired lower bound.

The issue is that distinct syllables may determine quasi-axes in the same hyperbolic space \(\mathcal{C}V\).
Lemmas~\ref{Lemma:bound} and~\ref{Lem:diameter of projection is small} show that such repetitions are compatible with the syllable order and that the corresponding quasi-axes have uniformly bounded mutual projections.
Lemma~\ref{Lem:order1} and Proposition~\ref{Lem:order2} then produce a strict total order on the quasi-axes contained in each fixed \(\mathcal{C}V\).
Finally, Proposition~\ref{Lem:overlap} assigns to these quasi-axes pairwise disjoint subsegments of a geodesic from \(\pi_V(x)\) to \(\pi_V(\mathbf{w}x)\), giving the required lower bound after summing over \(V\in\mathcal{D}_{\mathbf{w}}\) and applying the distance formula.

\begin{lemma}[$\beta_i$ is relevant]\label{Lemma:Relative quasigeodesic}
For each $\beta_{i}\in\mathcal{Q}_{\mathbf{w}}$, we have $\dist_{\beta_{i}}(x,\mathbf{w}x)> 100\vert e_{i} \vert \Theta$. In particular, each $V_{i}\in\mathcal{D}_{\mathbf{w}}$ is a $90\vert e_{i} \vert \Theta$-relevant domain with respect to $x$ and $\mathbf{w}x$.
\end{lemma}
\begin{proof}
We use induction on the syllable length $k$ of $w$ (or the $\Gamma$-syllable length of $\mathbf{w}$) for the proof of the first inequality. 

Suppose that \(k=1\).
By Lemmas~\ref{lemma:commuting projection maps} and~\ref{Lem:UniformQuasiAxis}, we have
\begin{equation}\label{Inequality:100}
\dist_{\alpha_1}(x,\mathbf{w}x) \ge 120\lvert e_1\rvert\Theta-\Theta-2\Theta >110\lvert e_1\rvert\Theta,
\end{equation}
as desired.

Suppose that \(k>1\). Fix \(1\le i\le k\), and define
\begin{align*}
r &= \max\Bigl(\{\,t<i\mid g_t \text{ does not \(\Gamma\)-commute with } g_i\,\}
\cup\{0\}\Bigr)\quad\mathrm{and}\\
s &=\min\Bigl(\{\,t>i\mid g_t \text{ does not \(\Gamma\)-commute with } g_i\,\}
\cup\{k+1\}\Bigr).
\end{align*}
Write $\mathbf w=ab_1g_ib_2c$, where
\[a=g_1\cdots g_r,\qquad b_1=g_{r+1}\cdots g_{i-1},\qquad b_2=g_{i+1}\cdots g_{s-1},\qquad c=g_s\cdots g_k,\]
with the convention that a product over an empty range is the identity.
Thus every syllable of \(b_1\) and \(b_2\) \(\Gamma\)-commutes with \(g_i\).
Moreover, if \(a\) is nontrivial, then its last syllable \(g_r\) does not \(\Gamma\)-commute with \(g_i\), and if \(c\) is nontrivial, then its first syllable \(g_s\) does not \(\Gamma\)-commute with \(g_i\).

By the triangle inequality, we have
\[\dist_{\beta_i}(x,\mathbf wx)=\dist_{\alpha_i}(\mathbf w_{i-1}^{-1}x,g_ib_2cx)\ge \dist_{\alpha_i}(b_2cx,g_ib_2cx)-\dist_{\alpha_i}(b_1^{-1}a^{-1}x,x) -\dist_{\alpha_i}(x,b_2cx)-2\Theta.\]
By the action of \(g_i\), the term $\dist_{\alpha_i}(b_2cx,g_ib_2cx)$ is bounded below by \(110|e_i|\Theta\), as in Inequality~\eqref{Inequality:100}.
Since every syllable of \(b_1\) and \(b_2\) \(\Gamma\)-commutes with \(g_i\), Condition~\eqref{Item:Sec3SubgroupHasUniformlyBoundedOrbits}, Lemma~\ref{Lem:CoarselyLipschitz}, and our choice of \(\Theta\) give
\[\dist_{\alpha_i}(b_1^{-1}a^{-1}x,x)\asymp_{1,\Theta} \dist_{\alpha_i}(a^{-1}x,x)\quad\text{and}\quad \dist_{\alpha_i}(x,b_2cx)\asymp_{1,\Theta} \dist_{\alpha_i}(x,cx).\]
Therefore, it remains to prove
\begin{equation}\label{Eq:bounded translation}
\max\bigl\{
\dist_{\alpha_i}(a^{-1}x,x),\,
\dist_{\alpha_i}(x,cx)
\bigr\}<2\Theta.
\end{equation}

If \(c\) is trivial, then \(\dist_{\alpha_i}(x,cx)=0\). Suppose that \(c\) is nontrivial.
Its syllable length is strictly smaller than \(k\), and its first syllable \(g_s\) does not \(\Gamma\)-commute with \(g_i\).
By the induction hypothesis, we have $\dist_{\alpha_s}(x,cx)>100|e_s|\Theta$.
Since \(g_s\) does not \(\Gamma\)-commute with \(g_i\), our choice of \(\Theta\) gives \(\dist_{\alpha_s}(x,\alpha_i)<\Theta\).
Hence
\[\dist_{\alpha_s}(cx,\alpha_i)\ge \dist_{\alpha_s}(x,cx)-\dist_{\alpha_s}(x,\alpha_i)-\Theta >98\Theta. \]
By Lemma~\ref{Lem:MSBI}, we obtain $\dist_{\alpha_i}(cx,\alpha_s)<\Theta$.
Since
\[\diam_{\mathcal{C}U_i}\Bigl(\mathfrak p_{\alpha_i}(\pi_{U_i}(x))\cup\mathfrak p_{\alpha_i}\bigl(\rho_{U_i}^{U_s}(\alpha_s)\bigr) \Bigr)<\Theta,\]
it follows that $\dist_{\alpha_i}(x,cx)<2\Theta$.

The estimate $\dist_{\alpha_i}(a^{-1}x,x)<2\Theta$ is proved similarly: it is immediate if \(a\) is trivial, and otherwise one applies the induction hypothesis to \(a^{-1}\), whose first syllable is \(g_r^{-1}\) and does not \(\Gamma\)-commute with \(g_i\).
This proves~\eqref{Eq:bounded translation}, and hence completes the inductive step.
\end{proof}

We begin with the observation that the projections of the $\mathbf{w}_i x$ onto $\beta_j$ exhibit a dichotomy: each is close to the projection onto $\beta_j$ of either $x$ or $\mathbf{w}x$.
This will be established using the same underlying philosophy as in the proof of Inequality~\eqref{Eq:bounded translation}.

\begin{lemma}\label{Lemma:bound}
Suppose that $1\le i<j\le k$. 
For indices $i\le j'<j''\le k$ and $1\le i'<i''< j$, the following inequality holds: $$\max\bigl\{\dist_{\beta_{i}}(\mathbf{w}_{j'}x,\,\mathbf{w}_{j''}x),\ \dist_{\beta_{j}}(\mathbf{w}_{i'}x,\,\mathbf{w}_{i''}x)\bigr\}<10\Theta.$$
\end{lemma}
\begin{proof}
We first consider the term involving $\beta_i$. Observe that
\[\dist_{\beta_{i}}(\mathbf{w}_{j'}x,\,\mathbf{w}_{j''}x)=\dist_{\alpha_{i}}(\mathbf{w}_{i-1}^{-1}\mathbf{w}_{j'}x,\,\mathbf{w}_{i-1}^{-1}\mathbf{w}_{j''}x)=\dist_{\alpha_{i}}(g_{i}\cdots g_{j'}x,\,g_{i}\cdots g_{j''}x).\]
If $g_i$ $\Gamma$-commutes with $g_{i+1},\dots,g_k$, then the last term is bounded above by $2\Theta$ since $$\pi_{U_i}(g_{i}\cdots g_{j'}x)=g_{i}\cdots g_{j'}\pi_{U_i}(x)\stackrel{\Theta}{=}g_{i}\cdots g_{j''}\pi_{U_i}(x)=\pi_{U_i}(g_{i}\cdots g_{j''}x).$$
Otherwise, let $m\ge i+1$ be the smallest index such that $g_{m}$ does not $\Gamma$-commute with $g_i$.
Define
\begin{align*}
A= \begin{cases}
\dist_{\alpha_{i}}(g_{m}\cdots g_{j'}x,g_{m}\cdots g_{j''}x) & \text{if $m\le j'$};\\
\dist_{\alpha_{i}}(x,g_{m}\cdots g_{j''}x) & \text{if $j'<m\le j''$};\\
0 & \text{if $j''<m$}.
\end{cases}
\end{align*}
Note that in all cases, we have $\dist_{\alpha_{i}}(g_{i}\cdots g_{j'}x, g_{i}\cdots g_{j''}x)\le A +4\Theta$. 

Since $g_m$ does not commute with $g_i$, by the argument used in the proof of Inequality~\eqref{Eq:bounded translation}, the distance in $\mathcal{C}U_i$ between the projections of $x$ and any term appearing in $A$ onto $\alpha_i$ is at most $2\Theta$.
Since $\diam_{\mathcal{C}U_i}(\mathfrak{p}_{\alpha_i}\circ\pi_{U_i}(x))<2\Theta$, it follows that $A<2\Theta+2\Theta+2\Theta$. Hence, \[\dist_{\beta_{i}}(\mathbf{w}_{j'}x,\,\mathbf{w}_{j''}x)<10\Theta.\]

The second inequality follows by an entirely analogous argument.
\end{proof}

We now show that when $V_i=V_j$, the projection of $\beta_j$ onto $\beta_i$ is uniformly bounded and has diameter significantly smaller than the lower bound $100|e_i|\Theta$ appearing in Lemma~\ref{Lemma:Relative quasigeodesic}. Using this fact, we will show that there is a strict total order on a subset of $\mathcal{Q}_{\mathbf{w}}$, each of whose elements lies in the same $\mathcal{C}V_i$.

\begin{lemma}\label{Lem:diameter of projection is small}
If $V=V_{i}=V_{j}$ for some $1\le i< j\le k$, then $g_i\prec_s g_j$ and $\diam_{\mathcal{C}V_{i}}(\mathfrak{p}_{\beta_{i}}(\beta_{j}))<40\Theta$. In particular, $\beta_i\neq \beta_j$.
\end{lemma}
\begin{proof}
We first show that \(g_i\prec_s g_j\).
Suppose, for contradiction, that \(g_i\) and \(g_j\) are not \(\prec_s\)-comparable.
Equivalently, the corresponding syllables \(x_i^{e_i}\) and \(x_j^{e_j}\) of \(w\) are not \(\prec_s\)-comparable.
By Remark~\ref{Rmk:IncomparableSyllablesCommute}, there exists \(w'\in\Min(y)\), obtained from \(w\) by shufflings, in which the syllables corresponding to \(x_i^{e_i}\) and \(x_j^{e_j}\) are adjacent and commute in \(A(\Gamma)\).
By Lemma~\ref{Lem:TranslatedDomainsInvariantUnderShuffling}, these two syllables still determine the translated supporting domains \(V_i\) and \(V_j\), respectively.

Let \(\mathbf{u}\) be the image under \(\phi\) of the prefix of \(w'\) preceding these two adjacent syllables. Since \(w'\in\Min(y)\), these syllables are powers of distinct standard generators.
Since they commute in \(A(\Gamma)\), their supporting domains \(U_i\) and \(U_j\) are orthogonal. Up to interchanging \(i\) and \(j\), the corresponding translated supporting domains are therefore $\mathbf{u}U_i$ and $\mathbf{u}g_iU_j=\mathbf{u}U_j$, where 
the equality follows because \(U_j\) is unbounded and orthogonal to \(U_i\), so \(g_iU_j=U_j\).
Hence these translated supporting domains are orthogonal. Thus $V_i\bot V_j$, contradicting the assumption that \(V_i=V_j\). Therefore \(g_i\) and \(g_j\) are \(\prec_s\)-comparable. Since \(i<j\) in the chosen representative \(w\), it follows that $g_i\prec_s g_j$.

\smallskip

Next, set $\gamma_i:=\mathfrak p_{\beta_i}(\beta_j)\subseteq\beta_i$ and $\gamma_j:=\mathfrak p_{\beta_j}(\beta_i)\subseteq\beta_j$, and assume for contradiction that $\diam_{\mathcal{C}V_i}(\gamma_i)\geq 40\Theta$. 
By Lemma~\ref{Lem:Projections of two quasiconvex}, we have $$\dist_{Haus}(\mathfrak{p}_{\beta_i}(\beta_j),\, \mathfrak{p}_{\beta_j}(\beta_i))<\Theta \implies \diam_{\mathcal{C}V}(\gamma_j)\ge 38\Theta\text{ and } \gamma_j\subseteq \mathcal{N}_{\Theta}(\beta_i).$$
We will then derive a contradiction using $\gamma_j$.

Note that $\gamma_j$ is $\Theta$-quasi-convex by Lemma~\ref{Lem:Projection of quasiconvex subsets}, and that $\mathfrak{p}_{\gamma_j}$ is $(1,\Theta)$-coarsely Lipschitz by Lemma~\ref{Lem:CoarselyLipschitz}, both with our choice of $\Theta$.
By Lemma~\ref{Lem:Projections of two quasiconvex} and the fact that $\gamma_j$ is the projection of a quasi-geodesic $\beta_i$, one can choose two points $\gamma^{\pm}_j\in\gamma_j$ such that
\begin{equation}\label{Ineq:gamma_j}
|\dist_{V}(\gamma^{+}_j,\,\gamma^{-}_j)-\diam_{\mathcal{C}V}(\gamma_j)|< 2\Theta.
\end{equation}

Now we need the following observation.

\begin{observation}
For each nonzero integer \(m\), let \(w^{j,m}\) be the word obtained from \(w\) by replacing its \(j\)-th syllable \(x_j^{e_j}\) with \(x_j^{me_j}\), and let
\[\mathbf w^{j,m}=g_1\cdots g_{j-1}g_j^m g_{j+1}\cdots g_k\]
be its image under \(\phi\).
Changing the nonzero exponent of the \(j\)-th syllable does not alter which syllables can be shuffled past one another or combined.
Hence \(w^{j,m}\) has the same syllable length as \(w\), with its syllables indexed in the same way.
Moreover,
\[\mathbf w^{j,m}_{j-1}=\mathbf w_{j-1}\qquad\text{and}\qquad\mathbf w^{j,m}_{i-1}=\mathbf w_{i-1},\]
so the translated quasi-axes \(\beta_i\) and \(\beta_j\), and hence also \(\gamma_j\), are independent of \(m\).

Applying Lemma~\ref{Lemma:bound} to \(\mathbf w^{j,m}\), with \(j'=j-1\) and \(j''=j\), gives 
\[\dist_{\beta_i}\bigl(\mathbf w_{j-1}x,\,\mathbf w^{j,m}_j x\bigr)<10\Theta\]
for every nonzero integer \(m\).
Since \(\gamma_j\subseteq\mathcal N_\Theta(\beta_i)\), nearest-point projection to \(\gamma_j\) coarsely factors through nearest-point projection to \(\beta_i\).
Hence, by our choice of \(\Theta\),
\begin{equation}\label{Eq:gamma_j_bounded}
\dist_V\!\left(\mathfrak p_{\gamma_j}\circ\pi_V(\mathbf w_{j-1}x),\,\mathfrak p_{\gamma_j}\circ\pi_V(\mathbf w^{j,m}_j x)\right)<13\Theta
\end{equation}
for every nonzero integer \(m\).
\end{observation}

Note that $h_j=\mathbf w_{j-1}g_j\mathbf w_{j-1}^{-1}$ preserves \(\beta_j\) and acts loxodromically on it, while $\mathbf w^{j,m}_j x=h_j^m(\mathbf w_{j-1}x)$.
By Lemma~\ref{lemma:commuting projection maps}, the projections of \(\mathbf w^{j,m}_j x\) to \(\beta_j\) escape in the two directions of \(\beta_j\) as \(m\to+\infty\) and \(m\to-\infty\).
Since \(\gamma_j\subseteq\beta_j\) is quasi-convex, its projection coarsely clamps points lying sufficiently far beyond its two extremal regions to uniformly bounded neighborhoods of those regions.
Thus, for suitable distinct nonzero integers \(m_1,m_2\), setting $q_\ell=\pi_V(\mathbf w^{j,m_\ell}_j x)$ for $\ell=1,2$,
the two points
\[\mathfrak p_{\gamma_j}\circ\mathfrak p_{\beta_j}(q_1)
\quad\text{ and }\quad
\mathfrak p_{\gamma_j}\circ\mathfrak p_{\beta_j}(q_2)
\]
lie in the \(\Theta\)-neighborhoods of the two coarse endpoints \(\gamma_j^\pm\), respectively.

By Lemma~\ref{Lem:Composition of Projection} and our choice of \(\Theta\),
\[
\mathfrak p_{\gamma_j}(q_1)\stackrel{\Theta}{=}\mathfrak p_{\gamma_j}\circ\mathfrak p_{\beta_j}(q_1)\qquad\text{and}\qquad\mathfrak p_{\gamma_j}(q_2)\stackrel{\Theta}{=}\mathfrak p_{\gamma_j}\circ\mathfrak p_{\beta_j}(q_2).
\]
It follows from~\eqref{Ineq:gamma_j} that
\[\dist_V\bigl(\mathfrak p_{\gamma_j}(q_1),\,\mathfrak p_{\gamma_j}(q_2)\bigr)>32\Theta.\]
On the other hand, by applying Inequality~\eqref{Eq:gamma_j_bounded} twice and using the triangle inequality, we obtain
\begin{align*}
\dist_V\!\left(\mathfrak p_{\gamma_j}(q_1),\,\mathfrak p_{\gamma_j}(q_2)\right)
&\le
\dist_V\!\left(\mathfrak p_{\gamma_j}(q_1),\,\mathfrak p_{\gamma_j}\circ\pi_V(\mathbf w_{j-1}x)\right)
+\diam_{\mathcal C V}\!\left(\mathfrak p_{\gamma_j}\circ\pi_V(\mathbf w_{j-1}x)\right)\\
&\qquad+\dist_V\!\left(\mathfrak p_{\gamma_j}\circ\pi_V(\mathbf w_{j-1}x),\,\mathfrak p_{\gamma_j}(q_2)\right)\\
&<27\Theta,
\end{align*}
a contradiction.
Therefore, $\diam_{\mathcal C V_i}\bigl(\mathfrak p_{\beta_i}(\beta_j)\bigr)<40\Theta$.
\end{proof}

\begin{definition}\label{order}
For $V\in\mathcal{D}_{\mathbf{w}}$, let $\mathcal{Q}_{\mathbf{w}}(V)=\{\beta\in\mathcal{Q}_{\mathbf{w}}\mid \beta\subseteq \mathcal{C}V\}$.
We define a relation $\prec_V$ on $\mathcal{Q}_{\mathbf{w}}(V)$ by 
declaring
\[\beta_{i} \prec_V \beta_{j} \quad \iff \quad \dist_{\beta_{i}}(x,\beta_{j})\geq 10\Theta.\]
\end{definition}

\begin{lemma}\label{Lem:order1}
Suppose that $g_{i}\prec_s g_{j}$ and $V=V_i=V_j$ for some $i< j$. 
Then $\dist_{\beta_{i}}(x,\beta_{j})\geq 30\vert e_{i} \vert \Theta$. In particular, $\beta_{i}\prec_V\beta_{j}$.
\end{lemma}
\begin{proof}
We claim that $\dist_{\beta_{i}}(\mathbf{w}x,\beta_{j})<30\Theta$.
If the claim holds, then $\mathfrak{p}_{\beta_{i}}(\beta_{j})$ is contained in the $70\Theta$-neighborhood of $\mathfrak{p}_{\beta_{i}}(\pi_{V_{i}}(\mathbf{w}x))$ since $\diam_{\mathcal{C}V}(\mathfrak{p}_{\beta_{i}}(\beta_{j}))<40\Theta$ by Lemma~\ref{Lem:diameter of projection is small}.
By Lemma~\ref{Lemma:Relative quasigeodesic}, it follows that $\mathfrak{p}_{\beta_{i}}(\beta_{j})$ is at least $30\vert e_{i} \vert \Theta$-far from $\mathfrak{p}_{\beta_{i}}(\pi_{V_{i}}(x))$ as desired. 

Choose \(z_i\in\pi_V(\mathbf w_i x)\) and \(z_j\in\pi_V(\mathbf w_j x)\), and let \(\gamma\subseteq\mathcal C V\) be a geodesic joining \(z_i\) to \(z_j\).
By Lemmas~\ref{Lemma:bound} and~\ref{Lem:Projections of two quasiconvex}, 
\[\diam_{\mathcal{C}V}(\mathfrak{p}_{\beta_i}(\gamma))< 10\Theta+2\Theta+2\Theta=14\Theta.\]
By Lemma~\ref{Lemma:bound} again, $\mathfrak{p}_{\beta_i}(\gamma)$ is $10\Theta$-close to $\mathfrak{p}_{\beta_{i}}(\pi_{V}(\mathbf{w}x))$. 
On the other hand, since 
\begin{align*}
\diam_{\mathcal{C}V}(\mathfrak{p}_{\beta_{j}}(\gamma))&\ge\dist_{\beta_j}(\mathbf{w}_ix, \mathbf{w}_jx)\ge\dist_{\beta_j}(x, \mathbf{w}x)-\dist_{\beta_j}(x, \mathbf{w}_ix)-\dist_{\beta_j}(\mathbf{w}_jx, \mathbf{w}x)-2\Theta\\
&>100\Theta-10\Theta-10\Theta-2\Theta=88\Theta    
\end{align*}
by Lemmas~\ref{Lemma:Relative quasigeodesic} and~\ref{Lemma:bound}, $\mathfrak{p}_{\gamma}(\beta_{j})$ 
is contained in the $\Theta$-neighborhood of $\mathfrak{p}_{\beta_{j}}(\gamma)$ by Lemma~\ref{Lem:Projections of two quasiconvex}. 
In particular, the projections $\mathfrak{p}_{\beta_{i}}(\beta_{j})$ and $\mathfrak{p}_{\beta_{i}}(\gamma)$ are at most $2\Theta$-far. 
Therefore, we have 
$$\dist_{\beta_{i}}(\mathbf{w}x,\beta_{j})\le \dist_{V}(\mathfrak{p}_{\beta_{i}}(\beta_{j}),\,\mathfrak{p}_{\beta_{i}}(\gamma))+
\dist_{V}(\mathfrak{p}_{\beta_i}(\gamma),\,\mathfrak{p}_{\beta_i}\circ\pi_{V}(\mathbf{w}x))+
\diam_{\mathcal{C}V}(\mathfrak{p}_{\beta_i}(\gamma))<26\Theta,$$ which completes the proof of the claim.
\end{proof}

\begin{proposition}\label{Lem:order2}
The relation $\prec_V$ on $\mathcal{Q}_{\mathbf{w}}(V)$ is a strict total order. More precisely, for distinct $\beta_i,\beta_j\in\mathcal{Q}_{\mathbf{w}}(V)$, we have $\beta_i\prec_V\beta_j$ if and only if $i<j$.
\end{proposition}
\begin{proof}
Suppose first that \(\beta_i,\beta_j\in\mathcal{Q}_{\mathbf{w}}(V)\) with \(i<j\).
Then \(V_i=V_j=V\), and hence Lemma~\ref{Lem:diameter of projection is small} gives $g_i\prec_s g_j$.
Applying Lemma~\ref{Lem:order1}, we obtain $\beta_i\prec_V\beta_j$.

Conversely, suppose that \(\beta_i\prec_V\beta_j\).
If \(j<i\), then the preceding argument, with \(i\) and \(j\) interchanged, would give \(\beta_j\prec_V\beta_i\).
This is impossible by Lemma~\ref{Lem:MSBI}, since the two inequalities
\[\dist_{\beta_i}(x,\beta_j)\ge 10\Theta\qquad\text{and}\qquad
\dist_{\beta_j}(x,\beta_i)\ge 10\Theta\]
cannot hold simultaneously. Therefore \(i<j\).

Thus \(\prec_V\) agrees with the usual order of the indices of the elements of \(\mathcal{Q}_{\mathbf{w}}(V)\), and is therefore a strict total order.
\end{proof}

For each domain $V\in\mathcal{D}_{\mathbf{w}}$, define $\mathcal{I}^{\mathbf{w}}_{V}=\{ j \mid \beta_{j}\subseteq \mathcal{C}V\}$ as the set of indices of the $\beta_j$ contained in $\mathcal{C}V$; in particular, $|\mathcal{I}^{\mathbf{w}}_{V}|=|\mathcal{Q}_{\mathbf{w}}(V)|$.

\begin{figure}[ht]
\centering
\resizebox{\linewidth}{!}{
\begin{tikzpicture}[x=1cm,y=1cm,line cap=round,line join=round,font=\normalsize]
  \definecolor{axisblue}{RGB}{0,70,145}
  \definecolor{projectiongreen}{RGB}{0,120,92}
  \definecolor{ballgreen}{RGB}{105,184,36}
  \tikzset{axis/.style={draw=axisblue,line width=0.85pt},axis portion/.style={draw=axisblue,line width=1.8pt},connector/.style={draw=black,line width=0.5pt,dash pattern=on 4pt off 3pt},axis connector/.style={draw=axisblue,line width=0.45pt,dash pattern=on 3pt off 2pt},projection/.style={draw=projectiongreen,line width=0.6pt},label/.style={inner sep=1pt,text=black},small label/.style={label,font=\small}}
  \path[use as bounding box] (0,0) rectangle (16,6.3);

  \coordinate (x) at (0.9,2.9);
  \coordinate (wx) at (15.2,2.9);
  \coordinate (ix) at (2.65,1.95);
  \coordinate (iwx) at (7.8,1.35);
  \coordinate (jx) at (7.95,4.25);
  \coordinate (jwx) at (12.95,3.85);

  \draw[connector] (x) .. controls (1.45,2.72) and (2.04,2.25) .. (ix);
  \draw[connector] (x) .. controls (3.5,3.8) and (6.2,3.5) .. (jx);
  \draw[connector] (iwx) .. controls (10.3,2.3) and (12.6,2.85) .. (wx);
  \draw[connector] (jwx) .. controls (13.45,3.21) and (14.4,3.05) .. (wx);

  \draw[axis connector] (2.3,0.2) .. controls (2.53,1.38) and (2.47,2.18) .. (2.53,2.9);
  \draw[axis connector] (7.3,6.1) .. controls (7.43,5.3) and (7.3,3.75) .. (7.46,2.9);
  \draw[axis connector] (7.64,2.9) .. controls (7.69,2.0) and (7.96,0.95) .. (8.3,0.35);
  \draw[axis connector] (13.9,5.65) .. controls (13.71,4.73) and (13.59,3.68) .. (13.68,2.9);

  \draw[axis] (2.3,0.15) .. controls (2.4,0.8) and (2.45,1.56) .. (ix);
  \draw[axis portion] (ix) .. controls (2.98,2.44) and (3.29,2.58) .. (4.0,2.58) .. controls (4.95,2.55) and (5.8,2.55) .. (6.3,2.57) .. controls (6.86,2.59) and (7.0,2.57) .. (7.29,2.12) .. controls (7.48,1.85) and (7.65,1.62) .. (iwx);
  \draw[axis] (iwx) .. controls (8.0,0.92) and (8.14,0.59) .. (8.25,0.35);

  \draw[axis] (7.35,6.1) .. controls (7.53,5.48) and (7.81,4.88) .. (jx);
  \draw[axis portion] (jx) .. controls (8.11,3.51) and (8.63,3.54) .. (9.18,3.4) .. controls (9.83,3.23) and (11.4,3.27) .. (11.83,3.28) .. controls (12.37,3.3) and (12.65,3.55) .. (jwx);
  \draw[axis] (jwx) .. controls (13.42,4.41) and (13.7,4.98) .. (13.85,5.65);

  \draw[projection] (7.25,6.08) .. controls (7.23,5.35) and (6.74,3.25) .. (6.4,2.57);
  \draw[projection] (11.86,3.28) .. controls (10.26,3.17) and (8.49,2.64) .. (7.57,1.7);
  \draw[axisblue,line width=0.6pt,decorate,decoration={brace,mirror,amplitude=3pt}] (6.4,2.49) -- (7.48,1.65);

  \draw[ballgreen,line width=0.5pt] (ix) circle[radius=1.48];
  \draw[ballgreen,line width=0.5pt] (jx) circle[radius=1.56];

  \draw[black,line width=0.55pt] (x) -- (wx);
  \draw[black,line width=1.9pt] (2.53,2.86) -- (7.75,2.86);
  \draw[black,line width=1.9pt] (7.46,2.94) -- (13.68,2.94);
  \foreach \p in {x,wx} \fill (\p) circle[radius=1.35pt];
  \foreach \p in {(2.53,2.86),(7.75,2.86),(7.46,2.94),(13.68,2.94)} \fill \p circle[radius=1.6pt];
  \foreach \p in {ix,iwx,jx,jwx} \fill[axisblue] (\p) circle[radius=1.25pt];

  \pgfmathsetmacro{\sigmaIRight}{2.65+sqrt(1.48^2-(2.86-1.95)^2)}
  \pgfmathsetmacro{\sigmaJRight}{7.95+sqrt(1.56^2-(2.94-4.25)^2)}
  \coordinate (sigmaIL) at (2.53,2.86);
  \coordinate (sigmaIR) at (\sigmaIRight,2.86);
  \coordinate (sigmaJL) at (7.46,2.94);
  \coordinate (sigmaJR) at (\sigmaJRight,2.94);
  \draw[projection,line width=0.8pt] (sigmaIL) -- (sigmaIR);
  \draw[projection,decorate,decoration={brace,amplitude=4pt}] (2.53,2.92) -- (\sigmaIRight,2.92);
  \draw[projection,line width=0.8pt] (sigmaJL) -- (sigmaJR);
  \draw[projection,decorate,decoration={brace,amplitude=3pt}] (7.46,3) -- (\sigmaJRight,3);

  \coordinate (sigmaIProjL) at (2.99,2.332);
  \coordinate (sigmaIProjR) at (\sigmaIRight,2.5764);
  \draw[projection,densely dotted,line width=0.45pt] (sigmaIL) -- (sigmaIProjL);
  \draw[projection,densely dotted,line width=0.45pt] (sigmaIR) -- (sigmaIProjR);
  \fill[projectiongreen] (sigmaIProjL) circle[radius=1.15pt];
  \fill[projectiongreen] (sigmaIProjR) circle[radius=1.15pt];

  \node[label,anchor=south] at (0.90,3.07) {$\pi_V(x)$};
  \node[label,anchor=north east] at (15.95,2.76) {$\pi_V(\mathbf{w}x)$};
  \node[small label,anchor=north west,fill=white] at (2.70,1.9) {$\mathfrak{p}_{\beta_i}\circ\pi_V(x)$};
  \node[small label,anchor=north east] at (7.75,1.35) {$\mathfrak{p}_{\beta_i}\circ\pi_V(\mathbf{w}x)$};
  \node[label,anchor=south,fill=white] at (3.12,3.05) {$\sigma_i$};
  \node[label,anchor=south] at (8.09,3.1) {$\sigma_j$};
  \node[small label,anchor=south west,fill=white] at (8.04,4.26) {$\mathfrak{p}_{\beta_j}\circ\pi_V(x)$};
  \node[label,anchor=west,font=\large] at (13.30,5.55) {$\beta_j$};
  \node[label,anchor=west,font=\large] at (7.70,0.43) {$\beta_i$};
  \node[label,anchor=north east] at (6.91,2) {$\mathfrak{p}_{\beta_i}(\beta_j)$};
\end{tikzpicture}%
}
\caption{The thick black segments shown on the geodesic from $\pi_V(x)$ to $\pi_V(\mathbf{w}x)$ are $\mathfrak{p}_{\gamma}(\beta_i)$ and $\mathfrak{p}_{\gamma}(\beta_j)$. The left and right circles are centered at $\mathfrak{p}_{\beta_i}(\pi_V(x))$ and $\mathfrak{p}_{\beta_j}(\pi_V(x))$, and have radii $10|e_i|\Theta$ and $10|e_j|\Theta$, respectively.}
\label{figure1}
\end{figure}

\begin{proposition}\label{Lem:overlap}
For each $V\in\mathcal{D}_{\mathbf{w}}$, we have 
$$\sum_{i\in \mathcal{I}^{\mathbf{w}}_{V}}  \vert e_{i}\vert \leq \frac{2}{5\Theta}\dist_{V}(x,\mathbf{w}x).$$
\end{proposition}
\begin{proof}
Choose \(z_-\in\pi_V(x)\) and \(z_+\in\pi_V(\mathbf wx)\), and let
\(\gamma\subseteq\mathcal C V\) be a geodesic joining \(z_-\) to \(z_+\).
For each $\beta_{j}\in \mathcal{Q}_{\mathbf{w}}(V)$, by Lemmas~\ref{Lemma:Relative quasigeodesic} and~\ref{Lem:Projections of two quasiconvex}, we have 
$$\dist_{Haus}(\mathfrak{p}_{\beta_{j}}(\gamma),\mathfrak{p}_{\gamma}(\beta_{j}))\leq \Theta\quad\text{and}\quad\diam_{\mathcal{C}V}(\mathfrak{p}_{\gamma}(\beta_{j}))\geq 100\vert e_{j} \vert \Theta-2\Theta>98\vert e_{j} \vert \Theta.$$ 
Then there exists a segment $\sigma_j$ of $\gamma$ with diameter at least $5\vert e_{j} \vert \Theta$ such that 
$$\sigma_{j}\subseteq \mathcal{N}_{\Theta}\bigl(\mathfrak{p}_{\gamma}(\beta_{j})\bigr)\cap \mathcal{N}_{10\vert e_{j} \vert \Theta} \bigl(\mathfrak{p}_{\beta_{j}}(\pi_V(x))\bigr).$$

Suppose that there are two elements $i,j\in\mathcal{I}_V^{\mathbf{w}}$ with $i<j$. By Proposition~\ref{Lem:order2}, we have $\beta_i\prec_V\beta_j$.
Then $\mathfrak{p}_{\beta_i}(\sigma_i)$ and $\mathfrak{p}_{\beta_i}(\sigma_j)$ are at least $18|e_i|\Theta$-far for the following reasons (see Figure~\ref{figure1}):
\begin{itemize}
\item $\mathfrak{p}_{\beta_i}(\sigma_i)$ is contained in the $(10|e_i|\Theta)$-neighborhood of $\mathfrak{p}_{\beta_i}(\pi_V(x))$. 
\item 
Since $\mathfrak{p}_{\beta_i}(\beta_j)$ is contained in the $70\Theta$-neighborhood of $\mathfrak{p}_{\beta_i}(\pi_V(\mathbf{w}x))$ as in the proof of Lemma~\ref{Lem:order1} and $\sigma_j$ is contained in the $\Theta$-neighborhood of $\beta_j$, $\mathfrak{p}_{\beta_i}(\sigma_j)$ is contained in the $72\Theta$-neighborhood of $\mathfrak{p}_{\beta_i}(\pi_V(\mathbf{w}x))$.
\item
By Lemma~\ref{Lemma:Relative quasigeodesic}, $\mathfrak{p}_{\beta_i}(x)$ and $\mathfrak{p}_{\beta_i}(\mathbf{w}x)$ are at least $100|e_i|\Theta$-far. 
\end{itemize}
Thus, $\sigma_i$ and $\sigma_j$ must be disjoint.

By Proposition~\ref{Lem:order2}, the segments $\sigma_i$, for $i\in\mathcal{I}^{\mathbf{w}}_{V}$, lie disjointly and in linear order along $\gamma$.
Therefore, we have 
$$ \sum_{i\in \mathcal{I}^{\mathbf{w}}_{V}}  5\vert e_{i}\vert \Theta\leq \diam_{\mathcal{C}V}(\gamma)\leq 2\dist_{V}(x,\mathbf{w}x)$$ as desired.
\end{proof}

We are now ready to complete the proof of the undistortion part of Theorem~\ref{Thm:MainEmbeddingThm_Precise} by showing that $\dist_G(x,\phi(w)x)$ grows at least linearly with the word length $|w|$.

\begin{proof}[Proof of the undistortion part of Theorem~\ref{Thm:MainEmbeddingThm_Precise}]
In order to show that $H$ is an undistorted subgroup, it suffices to show that there exist constants $A\geq 1$ and $B\geq 0$, independent of $w$, such that 
\[\dist_{G}(x,\mathbf{w}x)\geq \frac{1}{A}|w|-B,\]
where $|w|$ is the word length of $w$ with respect to the standard generating set of $A(\Gamma)$.

By Lemma~\ref{Lemma:Relative quasigeodesic}, $\mathcal{D}_{\mathbf{w}}\subseteq\operatorname{Rel}(x,\mathbf{w}x,\Theta)$.
By Lemma~\ref{Lem:overlap} and the distance formula, we have
\begin{align*}
|w|&= \sum_{i=1}^{k} \vert e_{i}\vert =\sum_{V\in\mathcal{D}_{\mathbf{w}}}\sum_{i\in \mathcal{I}^{\mathbf{w}}_{V}}  \vert e_{i}\vert  \leq \frac{2}{5\Theta}\sum_{V\in \mathcal{D}_{\mathbf{w}}} \dist_{V}(x,\mathbf{w}x)\\
& \leq \frac{2}{5\Theta}\sum_{V\in \operatorname{Rel}(x,\mathbf{w}x,\Theta)} \dist_{V}(x,\mathbf{w}x)\leq 
\frac{2K_{\Theta}}{5\Theta}\dist_{G}(x,\mathbf{w}x)+\frac{2K_{\Theta}}{5\Theta},
\end{align*}
where $K_{\Theta}$ is the constant obtained from Theorem~\ref{Thm:DistanceFormula} based on the fact that $\Theta>s_0$ by our choice of $\Theta$.
This completes the proof.
\end{proof}

\begin{remark}
Indeed, using Lemma~\ref{Lem:MSBI}, we can generalize Lemmas~\ref{Lem:diameter of projection is small} and~\ref{Lem:order1} by allowing $V_i$ and $V_j$ to be distinct.
This leads to the construction of a strict partial order on $\mathcal{Q}_\mathbf{w}$ which respects the strict partial order on $\syl(y)$.
However, since this refinement is not relevant to the proof of undistortion, we omit the details.
\end{remark}

\section{Variants of the main embedding theorem}\label{Sec:WithAssumptionsonG}

Given a geometrically irredundant collection of finitely many axial elements in an HHG, each fully supported on an unbounded domain, verifying whether this collection satisfies the conditions in Theorem~\ref{Thm:MainEmbeddingThm_Precise} can be challenging.  
Here we describe two settings in which Conditions~\eqref{Item:Sec3PowersCommute} and~\eqref{Item:Sec3SubgroupHasUniformlyBoundedOrbits} of Theorem~\ref{Thm:MainEmbeddingThm_Precise} hold automatically, leading to variants of the main theorem.  
The first arises from imposing natural assumptions on HHGs, while the second follows from strengthening the conditions on each axial element.


\subsection{Assumptions on HHGs}\label{Subsec:AssumptionsOnHHGs}

Motivated by mapping class groups, we impose the following three structural assumptions on an HHG \((G,\mathfrak S)\). Before stating them, we recall that the first assumption can always be arranged by discarding inessential bounded domains.

\begin{lemma}[Passage to the essential domains, \cite{ABR25}]\label{Lem:EssentialReduction}
Let \((G,\mathfrak S)\) be an HHG, and let
\[\mathfrak S_{\mathrm{ess}}=\left\{U\in\mathfrak S \mid \text{there exists an unbounded domain }V\sqsubseteq U \right\}.\]
Then \((G,\mathfrak S_{\mathrm{ess}})\) is an HHG whose hyperbolic spaces, projections, and hierarchical relations are inherited from \((G,\mathfrak S)\). Moreover, every \(\sqsubseteq\)-minimal domain in \(\mathfrak S_{\mathrm{ess}}\) is unbounded.

If \(\mathfrak S_{\mathrm{ess}}=\emptyset\), then \(G\) is finite.
\end{lemma}

\begin{proof}
Since the action of \(G\) on \(\mathfrak S\) is cofinite, the diameters of the bounded hyperbolic spaces \(\mathcal C U\) are uniformly bounded. 
Thus \((G,\mathfrak S)\) satisfies the bounded domain dichotomy, and \cite[Lemma~4.2]{ABR25} implies that \((G,\mathfrak S_{\mathrm{ess}})\) is an HHS with the inherited hyperbolic spaces, projections, and hierarchical relations.
The set \(\mathfrak S_{\mathrm{ess}}\) is \(G\)-invariant, and the action of \(G\) on it is cofinite. Hence \((G,\mathfrak S_{\mathrm{ess}})\) is an HHG.

Every unbounded domain belongs to \(\mathfrak S_{\mathrm{ess}}\), by taking \(V=U\). If \(U\in\mathfrak S_{\mathrm{ess}}\) is \(\sqsubseteq\)-minimal, then there exists an unbounded \(V\sqsubseteq U\); minimality gives \(V=U\), so \(U\) is unbounded.

Finally, if \(\mathfrak S_{\mathrm{ess}}=\emptyset\), then the uniqueness axiom implies that \(G\) is bounded. Since \(G\) is finitely generated, it follows that \(G\) is finite.
\end{proof}


\begin{assumption}\label{Assumption:Essential}
Every \(\sqsubseteq\)-minimal domain is unbounded.
\end{assumption}

\begin{assumption}\label{Assumption:CoboundedAction}
There exists a constant $C$ such that for any non-$\sqsubseteq$-minimal domain $U\in\mathfrak{S}$, the set $\{\rho_U^W\mid W\sqsubsetneq U\}$ is $C$-coarsely dense in $\mathcal{C}U$.
\end{assumption}

\begin{assumption}\label{Assumption:CentralExtension}
For every unbounded \(\sqsubseteq\)-minimal domain \(U\in\mathfrak S\) and every domain \(T\in\mathfrak S\) with \(U\sqsubsetneq T\), if $\{\,V\in\mathfrak S\mid V\sqsubseteq T \text{ and } V\bot U\,\}\neq \emptyset$, then there exists a container \(U_T^\bot\) for the orthogonal complement of \(U\) in \(T\) such that \(U_T^\bot\bot U\); such a container is necessarily unique.

Furthermore, for each unbounded \(\sqsubseteq\)-minimal domain \(U\in\mathfrak S\), there exists an axial element \(f_U\in G\), fully supported on \(U\), such that
\[1\longrightarrow \langle f_U\rangle\longrightarrow \operatorname{Stab}_G(U) \xlongrightarrow{\varphi_U} H_U \longrightarrow 1\]
is a central extension, where \(H_U\) is an HHG with index set \(\mathfrak S_U^\bot:=\{W\in\mathfrak S\mid W\bot U\}\), whose HHG structure and stabilizer extensions are inherited from \((G,\mathfrak S)\).
\end{assumption}

\begin{remark}[Discussion of the assumptions]\label{Rmk:DiscussionOfAssumptions}
\begin{enumerate}[wide]
\item 
Assumption~\ref{Assumption:Essential}, considered by itself, does not restrict the underlying infinite group, by Lemma~\ref{Lem:EssentialReduction}. However, that lemma does not assert that Assumptions~\ref{Assumption:CoboundedAction} and~\ref{Assumption:CentralExtension} are preserved under passage to the essential domains. Accordingly, throughout the structural version of our theorem, all three assumptions are imposed on the same chosen HHG structure.
\item 
Assumption~\ref{Assumption:CoboundedAction} is implied by the more geometric requirement that, for every domain \(T\in\mathfrak S\), the stabilizer \(\operatorname{Stab}_G(T)\) act coboundedly on the standard product region associated to \(T\). We do not use this formulation here, since standard product regions have not otherwise been introduced; see \cite[\S15 and \S50]{CHK}.
\item 
Recall that an HHG \((G,\mathfrak S)\) has \emph{clean containers} if, whenever \(U\sqsubsetneq T\) and there exists a domain nested in \(T\) and orthogonal to \(U\), the orthogonal complement of \(U\) in \(T\) admits a container which is itself orthogonal to \(U\). Such an orthogonal container is necessarily unique. The first clause of Assumption~\ref{Assumption:CentralExtension} is this condition restricted to unbounded \(\sqsubseteq\)-minimal domains, and is therefore weaker than the full clean-container condition.
\item 
All three assumptions are inherited by the auxiliary HHGs \((H_U,\mathfrak S_U^\bot)\) appearing in Assumption~\ref{Assumption:CentralExtension}. First, suppose that \(V\) is \(\sqsubseteq\)-minimal in \(\mathfrak S_U^\bot\). If \(W\sqsubsetneq V\) in \(\mathfrak S\), then \(V\bot U\) implies \(W\bot U\), so \(W\in\mathfrak S_U^\bot\), contradicting the minimality of \(V\). Thus \(V\) is \(\sqsubseteq\)-minimal in \(\mathfrak S\), and hence is unbounded by Assumption~\ref{Assumption:Essential}. Therefore, Assumption~\ref{Assumption:Essential} is inherited.

Similarly, if \(T\in\mathfrak S_U^\bot\) and \(W\sqsubsetneq T\) in \(\mathfrak S\), then \(T\bot U\) implies \(W\bot U\). Hence the domains properly nested in \(T\) in \(\mathfrak S_U^\bot\) are exactly those in \(\mathfrak S\). It follows that Assumption~\ref{Assumption:CoboundedAction} is inherited with the same coarse-density constant.

Finally, consider the container clause of Assumption~\ref{Assumption:CentralExtension}. Let \(W,T\in\mathfrak S_U^\bot\), where \(W\) is an unbounded \(\sqsubseteq\)-minimal domain and \(W\sqsubsetneq T\), and suppose that there exists a domain nested in \(T\) and orthogonal to \(W\). The container \(W_T^\bot\) supplied by Assumption~\ref{Assumption:CentralExtension} in \((G,\mathfrak S)\) satisfies \(W_T^\bot\sqsubseteq T\) and \(W_T^\bot\bot W\).
Since \(T\bot U\), the relation \(W_T^\bot\sqsubseteq T\) also implies \(W_T^\bot\bot U\). Thus \(W_T^\bot\in\mathfrak S_U^\bot\), and it gives the required container in the auxiliary HHG. The stabilizer-extension clause is inherited by the recursive inheritance requirement in Assumption~\ref{Assumption:CentralExtension}.
Consequently, each auxiliary HHG \((H_U,\mathfrak S_U^\bot)\) satisfies all three assumptions.
\item The recursive inheritance established above allows the weak-commutativity argument to be applied inductively inside the auxiliary groups \(H_U\); see the proof of Proposition~\ref{Prop:Weak commutativity of FSE}.
\end{enumerate}
\end{remark}

\begin{remark}[Standard short HHGs]\label{Rmk:StandardShortHHGs}
The standard short HHG structures on the examples mentioned in the introduction satisfy Assumptions~\ref{Assumption:CoboundedAction} and~\ref{Assumption:CentralExtension}. Indeed, the cobounded-product-region property established in \cite[Lemma~2.13]{Man24} gives Assumption~\ref{Assumption:CoboundedAction}.
The unbounded \(\sqsubseteq\)-minimal domains are the quasiline domains associated to infinite cyclic directions, and the simplicial-container description in \cite[Lemma~2.12]{Man24} gives the container clause of Assumption~\ref{Assumption:CentralExtension} for these domains.
In the examples under consideration, each cyclic direction is central in the corresponding vertex stabilizer, and its generator is fully supported on the associated quasiline domain \cite[Section~2.2, Axioms~B--C, and Remark~2.7]{Man24}. The quotients carry the inherited HHG structures on the orthogonal factors, together with the required stabilizer extensions \cite[proof of Lemma~2.17, especially Claim~2.18]{Man24}.

These structures may nevertheless contain bounded \(\sqsubseteq\)-minimal domains, so they need not satisfy Assumption~\ref{Assumption:Essential}.
Although this assumption can always be arranged by passing to the essential domains, Lemma~\ref{Lem:EssentialReduction} does not by itself establish that Assumptions~\ref{Assumption:CoboundedAction} and~\ref{Assumption:CentralExtension} are preserved. Thus these standard short HHGs provide promising candidates for \(\mathcal P_3\), with the remaining task being to realize all three assumptions simultaneously in a single HHG structure.
\end{remark}

Under the assumptions above, the two technical conditions in Theorem~\ref{Thm:MainEmbeddingThm_Precise} follow from the stated hypotheses on the chosen elements. We therefore obtain the following more precise formulation of Theorem~\ref{Thm:MainTheoremAssumptiononG}.

\begin{theorem}\label{Thm:MainTheoremAssumptiononG_Precise}
Let $(G,\mathfrak{S})$ be a hierarchically hyperbolic group satisfying Assumptions~\ref{Assumption:Essential},~\ref{Assumption:CoboundedAction} and~\ref{Assumption:CentralExtension}.
Let $f_i\in G$ be an axial element fully supported on an unbounded domain $U_i\in\mathfrak{S}$ for $i=1,\dots,m$ such that $\{ f_1, \dots,f_m\}$ is geometrically irredundant.
Then there exist constants $N = N(\{ f_i \})>0$ and $D = D(\{ f_i \})>0$ such that for any $d\ge D$,  
$$H=\langle f_1^{d N},\dots,f_m^{d N}\rangle\cong A(\mathcal{O}^{\Sigma}),$$ 
where $\mathcal{O}^{\Sigma}$ is the orthogonality graph of $\Sigma=\{U_i\mid i=1,\dots,m\}$.

Moreover, if no $U_i$ is $\sqsubseteq$-minimal, after possibly increasing $D$ if necessary, the subgroup $H$ is quasi-isometrically embedded in $G$.
\end{theorem}

The proof of Theorem~\ref{Thm:MainTheoremAssumptiononG_Precise} relies on two facts, stated below.

\begin{lemma}\label{Lem:UniformlyBoundedOrbitsforNonNestingMinimal}
Let $(G,\mathfrak{S})$ be an HHG satisfying Assumption~\ref{Assumption:CoboundedAction}. 
If $U\in\mathfrak{S}$ is a non-$\sqsubseteq$-minimal domain, then an element $f\in\operatorname{Stab}_G(U)$ fixing every \(\sqsubseteq\)-minimal domain $V\sqsubsetneq U$ acts on $\mathcal{C}U$ with uniformly bounded orbits.
\end{lemma}
\begin{proof}
If $U$ is bounded, then the lemma is obvious.

Suppose that $U$ is unbounded and let $x\in\mathcal{C}U$. By Assumption~\ref{Assumption:CoboundedAction}, there exists a domain $V$ properly nested in $U$ such that $\rho_{U}^{V}$ is $C$-close to $x$ for the constant \(C\) in Assumption~\ref{Assumption:CoboundedAction}. 
If \(V\) is not \(\sqsubseteq\)-minimal, choose a \(\sqsubseteq\)-minimal domain \(V'\sqsubsetneq V\). Then \(\rho_U^{V'}\) and \(\rho_U^V\) are \(E\)-close, so after replacing \(V\) by \(V'\), we may assume that \(V\) is \(\sqsubseteq\)-minimal.

As $f$ fixes every \(\sqsubseteq\)-minimal domain properly nested in $U$, we have $f^nV=V$ for all $n\in\mathbb{Z}$, and hence $f^n(\rho_U^V)=\rho_U^{f^nV}=\rho_U^V$.
As $f$ acts on $\mathcal{C}U$ by isometries and $\diam_{\mathcal{C}U}(\rho_U^V)\le E$, it follows that $\dist_U(x,f^n x)\le 2C+6E$ for all $n\in\mathbb{Z}$. 
Thus the $f$-orbit of $x$ has uniformly bounded diameter in $\mathcal{C}U$.
\end{proof}




\begin{proposition}\label{Prop:Weak commutativity of FSE}
Let $(G,\mathfrak{S})$ be an HHG satisfying Assumptions~\ref{Assumption:Essential},~\ref{Assumption:CoboundedAction} and~\ref{Assumption:CentralExtension}.
Suppose that $f,g\in G$ are two axial elements fully supported on unbounded domains $U, V\in\mathfrak{S}$, respectively, with $U\bot V$.
Then there exists $N>0$ such that $f^{N}$ and $g^{N}$ commute.    
\end{proposition}
\begin{proof}
Suppose first that \(U\) is \(\sqsubseteq\)-minimal. Since \(f\) is fully supported on \(U\), we have \(\B(f)=\{U\}\), and \(f\) acts with uniformly bounded orbits on \(\mathcal{C}W\) for every unbounded \(W\bot U\).
Since the HHG structure of \(H_U\) is inherited from \((G,\mathfrak S)\), the image \(\varphi_U(f)\) has uniformly bounded orbits on \(\mathcal CW\) for every unbounded \(W\in\mathfrak S_U^\bot\). Hence \(\B_{H_U}(\varphi_U(f))=\emptyset\). 
It follows that \(\varphi_U(f)\) is elliptic and hence has finite order. Thus
\(f^N\in\ker\varphi_U=\langle f_U\rangle\) for some \(N>0\).
Since \(\langle f_U\rangle\) is central in \(\operatorname{Stab}_G(U)\) and \(g\in\operatorname{Stab}_G(U)\), it follows that \(f^N\) commutes with \(g\), and hence with \(g^N\).
The same argument applies if \(V\) is \(\sqsubseteq\)-minimal. Thus we may assume that neither \(U\) nor \(V\) is \(\sqsubseteq\)-minimal.


Define the following five subsets of $\mathfrak{S}$ whose union is all of $\mathfrak{S}$:
\begin{equation*}
\begin{gathered}
\mathfrak{S}_{1}=\mathfrak{S}_{U},\quad \mathfrak{S}_{2}=\mathfrak{S}_{V},\quad \mathfrak{S}_{3}=\{ W\in\mathfrak{S}  \mid W\bot \{U, V\}\},\\
\mathfrak{S}_{4}=\{ W\in\mathfrak{S}  \mid U\sqsubsetneq W \text{ or } U\pitchfork W\}\quad\text{and}\quad\mathfrak{S}_{5}=\{ W\in\mathfrak{S}  \mid V\sqsubsetneq W\text{ or }V\pitchfork W\}.
\end{gathered}    
\end{equation*}
Since $f$ and $g$ are fully supported, every element of $\langle f,g\rangle$ stabilizes both $U$ and $V$, and therefore acts by permutation on each $\mathfrak{S}_{i}$.

We begin by proving a special case of the proposition, which is the base step for induction.

\begin{claim}\label{claim}
Suppose $(G,\mathfrak{S})$, $f$, $g$, $U$, $V$, and the sets $\mathfrak{S}_1,\dots,\mathfrak{S}_5$ are as above.  
If $\mathfrak{S}_3$ contains no unbounded domain, then the proposition holds.
\end{claim}
\begin{proof}[Proof of Claim~\ref{claim}]
\renewcommand{\qedsymbol}{$\blacksquare$}
It suffices to show that $\langle f,g \rangle$ contains no free subgroup of rank $2$, since Theorem~\ref{Thm:TitsAlternative} would then imply that $\langle f,g \rangle$ is virtually abelian.

Suppose, for contradiction, that $a,b\in \langle f,g\rangle$ generate a free subgroup of rank $2$. Consider the commutator $h=[a,b]=aba^{-1}b^{-1}$.  
We will show that \(\B(h)=\emptyset\). By Proposition~\ref{Prop:AxialEllipticDichotomy}, this implies that \(h\) has finite order, contradicting the freeness of \(\langle a,b\rangle\).

\begin{stepA}
$h$ fixes every unbounded domain in $\mathfrak{S}_{1}\cup\mathfrak{S}_{2}$.
\end{stepA}
Writing
\begin{equation}\label{Eq:AandB}
a=f^{n_{1}}g^{n_{2}}\cdots f^{n_{2k-1}}g^{n_{2k}}\quad\text{and}\quad 
b=f^{m_{1}}g^{m_{2}}\cdots f^{m_{2l-1}}g^{m_{2l}},
\end{equation}
we note that $f^{n_i}$ permutes $\mathfrak{S}_{1}$ while $g^{m_j}$  fixes every unbounded domain in $\mathfrak{S}_{1}$.  
Hence, whenever \(W\in\mathfrak{S}_1\) is an unbounded domain, we have
\[a.W=f^{\sum_i n_{2i-1}}.W\quad\text{and}\quad b.W=f^{\sum_j m_{2j-1}}.W \]
and, in particular, \( h.W=W \).
The same holds for every unbounded domain \(W\in\mathfrak{S}_{2}\).

\begin{stepA}\label{S_1}
$h$ acts elliptically on $\mathcal{C}W$ for every unbounded domain $W\in\mathfrak{S}_1\cup\mathfrak{S}_2$.
\end{stepA}
If $W$ is an unbounded non–$\sqsubseteq$-minimal domain in \(\mathfrak{S}_1\), then $h$ fixes every \(\sqsubseteq\)-minimal domain properly nested in \(W\), and hence has uniformly bounded orbits on $\mathcal{C}W$ by Lemma~\ref{Lem:UniformlyBoundedOrbitsforNonNestingMinimal}. 
Indeed, every \(\sqsubseteq\)-minimal domain properly nested in \(W\) is unbounded by Assumption~\ref{Assumption:Essential}, and hence is fixed by \(h\) by the preceding step.

Suppose now that \(W\sqsubsetneq U\) is \(\sqsubseteq\)-minimal. By Assumption~\ref{Assumption:Essential}, \(W\) is unbounded. Since \(f\) acts loxodromically on \(\mathcal CU\), there exists \(k>0\) such that
\[\dist_U(\rho_U^W,\rho_U^{f^kW})>10E.\]
As \(f\) stabilizes \(U\) and preserves \(\sqsubseteq\)-minimality, the domains \(W\) and \(f^kW\) are distinct \(\sqsubseteq\)-minimal domains nested in \(U\), and hence neither is nested in the other. 
They are not orthogonal, since otherwise \cite[Lemma~1.5]{DHS17} would give \(\dist_U(\rho_U^W,\rho_U^{f^kW})\le 2E\). 
Thus \(W\pitchfork f^kW\). Moreover, both domains are unbounded members of \(\mathfrak S_1\), so the preceding step gives \(hW=W\) and \(h(f^kW)=f^kW\). By equivariance, 
\[h\rho_W^{f^kW}=\rho_{hW}^{h(f^kW)}=\rho_W^{f^kW}.\]
Since \(\rho_W^{f^kW}\) is a nonempty bounded subset of \(\mathcal CW\), the isometry induced by \(h\) has a bounded orbit, and hence \(h\) acts elliptically on \(\mathcal CW\).



\begin{stepA}\label{S_4}
$\B(h)\cap(\mathfrak{S}_3\cup\mathfrak{S}_4\cup\mathfrak{S}_5)=\emptyset$.
\end{stepA}
Let $W\in\mathfrak{S}_4$. Choose $p\in\mathcal{C}U$ and $q\in\mathcal{C}V$, and let $x\in G$ be a partial realization point of $\{p,q\}$.  
By the construction of $x$, we have $\dist_{W}(x,\rho^U_W)<E$ and thus $\dist_{hW}(hx,h\rho^U_W)<E$.
Moreover, since 
\[h\rho^U_W=\rho^{hU}_{hW}=\rho^U_{hW}\quad\text{and}\quad hW\in\mathfrak{S}_4,\]
we also have $\dist_{hW}(x,h\rho^U_{W})<E$, and thus $\dist_{hW}(hx,x)<3E$.
The same argument, with \(h^n\) in place of \(h\), gives $\dist_{h^nW}(h^nx,x)<3E$ for every \(W\in\mathfrak S_4\) and \(n\in\mathbb Z\).
Since \(h^n\) permutes \(\mathfrak S_4\), reindexing the domains yields \(\dist_W(x,h^nx)<3E\) for every \(W\in\mathfrak S_4\).
Therefore, we have $W\notin \B(h)$, and similarly, after replacing $U$ by $V$, any domain in $\mathfrak{S}_5$ is not in $\B(h)$.
Since $\mathfrak{S}_3$ is assumed to have no unbounded domains, clearly $\B(h)\cap\mathfrak{S}_3=\emptyset$.
Hence $\B(h)\cap(\mathfrak{S}_3\cup\mathfrak{S}_4\cup\mathfrak{S}_5)=\emptyset$.

\smallskip
By Steps~\ref{S_1} and~\ref{S_4}, $\B(h)$ is an empty set, which is a contradiction. 
Therefore, \(\langle f,g\rangle\) is virtually abelian. Hence, after passing to positive powers, \(f\) and \(g\) commute; that is, there exists \(N>0\) such that \(f^N\) and \(g^N\) commute.
\end{proof}

We argue by induction on the maximal cardinality of a collection of pairwise orthogonal unbounded \(\sqsubseteq\)-minimal domains in \(\mathfrak S_3\). 
By \cite[Lemma~2.1]{BHS19}, the cardinality of a pairwise orthogonal collection of domains is uniformly bounded, so the induction parameter is finite. By finite complexity, every domain contains a \(\sqsubseteq\)-minimal domain. Moreover, if an unbounded domain lies in \(\mathfrak S_3\), then every domain nested in it also lies in \(\mathfrak S_3\), and a minimal such domain is unbounded by Assumption~\ref{Assumption:Essential}. Hence the induction parameter is \(0\) if and only if \(\mathfrak S_3\) contains no unbounded domain, which is precisely the case covered by Claim~\ref{claim}.

Assume that the proposition holds whenever, in the corresponding setting, the subset of domains orthogonal to both supporting domains contains no collection of more than \(k\) pairwise orthogonal unbounded \(\sqsubseteq\)-minimal domains.
Suppose now that $\{W_1,\dots,W_{k+1}\}$ is a maximal collection of pairwise orthogonal unbounded $\sqsubseteq$-minimal domains in $\mathfrak{S}_3$.  
By Assumption~\ref{Assumption:CentralExtension}, there is a central extension
\[
1 \longrightarrow \langle w_1\rangle \longrightarrow \operatorname{Stab}_{G}(W_{1}) 
\xlongrightarrow{\varphi} H_{W_{1}} \longrightarrow 1.
\]
Since \(W_1\bot\{U,V\}\), both \(f\) and \(g\) lie in \(\operatorname{Stab}_G(W_1)\). 
By Remark~\ref{Rmk:DiscussionOfAssumptions}, the HHG \((H_{W_1},\mathfrak S_{W_1}^{\bot})\) satisfies Assumptions~\ref{Assumption:Essential} and~\ref{Assumption:CoboundedAction}; it satisfies Assumption~\ref{Assumption:CentralExtension} by the inheritance clause in that assumption.
Since the HHG structure of \((H_{W_1},\mathfrak S_{W_1}^{\bot})\) is inherited from \((G,\mathfrak S)\), the stabilization, loxodromic, and bounded-orbit conditions defining full support pass to the images of \(f\) and \(g\). Thus \(\varphi(f)\) and \(\varphi(g)\) are axial elements fully supported on \(U\) and \(V\), respectively.

Any collection of pairwise orthogonal unbounded \(\sqsubseteq\)-minimal domains in
\(\mathfrak S_{W_1}^{\bot}\cap\mathfrak S_3\) has cardinality at most \(k\); otherwise, adjoining \(W_1\) would contradict the maximality of \(\{W_1,\dots,W_{k+1}\}\).
By the induction hypothesis, there exists $N'>0$ with $\varphi(f^{N'})$ and $\varphi(g^{N'})$ commuting, and thus $[f^{N'},g^{N'}]\in\ker\varphi=\langle w_1\rangle$.
Since \(\langle w_1\rangle\) is central in \(\operatorname{Stab}_G(W_1)\), $\{ w_1,f^{N'},g^{N'}\}$ generates a nilpotent subgroup, which by Theorem~\ref{Thm:TitsAlternative} is virtually abelian.
Therefore, after passing to further positive powers if necessary, \(f^{N'}\) and \(g^{N'}\) commute. This completes the induction and hence the proof.
\end{proof}

Now, we can complete the proof of Theorem~\ref{Thm:MainTheoremAssumptiononG_Precise}.

\begin{proof}[Proof of Theorem~\ref{Thm:MainTheoremAssumptiononG_Precise}]
By Proposition~\ref{Prop:Weak commutativity of FSE}, Condition~\eqref{Item:Sec3PowersCommute} in Theorem~\ref{Thm:MainEmbeddingThm_Precise} holds and thus the first part of Theorem~\ref{Thm:MainTheoremAssumptiononG_Precise} holds by the embedding part of Theorem~\ref{Thm:MainEmbeddingThm_Precise}. 

For the moreover statement, suppose that every supporting domain \(U_i\) is non-\(\sqsubseteq\)-minimal. Fix \(i\), and let \(W\sqsubsetneq U_i\) be a \(\sqsubseteq\)-minimal domain. By Assumption~\ref{Assumption:Essential}, the domain \(W\) is unbounded.
If \(U_j\bot U_i\), then \(W\bot U_j\), and hence \(f_jW=W\) because \(f_j\) is fully supported on \(U_j\).
It follows that every element of \(\langle f_j\mid U_j\bot U_i\rangle\) fixes every \(\sqsubseteq\)-minimal domain properly nested in \(U_i\).
Hence Lemma~\ref{Lem:UniformlyBoundedOrbitsforNonNestingMinimal} shows that Condition~\eqref{Item:Sec3SubgroupHasUniformlyBoundedOrbits} holds with \(N_1=1\).
Therefore, the moreover conclusion follows from the undistortion part of Theorem~\ref{Thm:MainEmbeddingThm_Precise}.
\end{proof}

\begin{remark}\label{Rmk:WhyNoNestingMinimal}
Even under Assumptions~\ref{Assumption:Essential},~\ref{Assumption:CoboundedAction} and~\ref{Assumption:CentralExtension}, $\sqsubseteq$-minimal domains may not have a bounded-orbit condition as in Lemma~\ref{Lem:UniformlyBoundedOrbitsforNonNestingMinimal}.
More precisely, when $f_i$ is an axial element fully supported on $U_i$ for $i=1,2$ (possibly, $U_1=U_2$) such that both $U_1$ and $U_2$ are orthogonal to a $\sqsubseteq$-minimal domain $V$, some words in $f_1$ and $f_2$ may act loxodromically on $\mathcal{C}V$.
This phenomenon occurs even in mapping class groups; thus it is natural to exclude $\sqsubseteq$-minimal domains for the undistortion part of Theorem~\ref{Thm:MainTheoremAssumptiononG_Precise}.
\end{remark}

\subsection{Assumptions on elements}\label{Subsec:AssumptionsOnElements}
We now impose no structural assumptions on the ambient HHG, but instead strengthen the support condition on the chosen axial elements. The following is a more precise statement of Theorem~\ref{Thm:MainTheoremStrong}.


\begin{theorem}\label{Thm:MainTheoremStrong_Precise}
Let $(G,\mathfrak{S})$ be a hierarchically hyperbolic group, and let $f_i\in G$ be rigidly fully supported on an unbounded domain $U_i\in\mathfrak{S}$ for $i=1,\dots,m$ with the collection $\{ f_1, \dots,f_m\}$ geometrically irredundant. 
Then there exist constants $N = N(\{ f_i \})>0$ and $D = D(\{ f_i \})>0$ such that if $d\ge D$, then the map $v_i\mapsto f^{d N}_i$ for each standard generator $v_i$ of $A(\mathcal{O}^{\Sigma})$ induces an injective homomorphism $A(\mathcal{O}^{\Sigma})\to G$ which is a quasi-isometric embedding.
\end{theorem}


To prove the theorem, we rely on the following analogue of Proposition~\ref{Prop:Weak commutativity of FSE}, which asserts that, without any assumptions on the HHG, the set of axial elements rigidly fully supported on single domains satisfies the weak commutativity property.

\begin{proposition}\label{Prop:weakCommutativityForStrongFSE}
Let $(G,\mathfrak{S})$ be an HHG. If $f,g\in G$ are two axial elements rigidly fully supported on orthogonal unbounded domains $U, V\in\mathfrak{S}$, respectively, with $U\bot V$, then there exists $N>0$ such that $f^{N}$ and $g^{N}$ commute.
\end{proposition}
\begin{proof}
The proof follows the same strategy as the proof of Proposition~\ref{Prop:Weak commutativity of FSE}, but is simpler because the elements are rigidly fully supported.

Suppose, for contradiction, that the subgroup $\langle f, g\rangle$ contains a rank $2$ free subgroup $\langle a,b\rangle$, and let $h=[a,b]$.
Use the same decomposition \(\mathfrak{S}=\mathfrak{S}_1\cup\cdots\cup\mathfrak{S}_5\) as in the proof of Proposition~\ref{Prop:Weak commutativity of FSE}.
By Step~\ref{S_4} of that proof, which only uses that \(f\) and \(g\) are fully supported, it follows that $\B(h)\cap(\mathfrak{S}_4\cup\mathfrak{S}_5)=\emptyset$.
If $W\in\mathfrak{S}_3$ is an unbounded domain, then $h$ acts trivially on $\mathcal{C}W$ since both $f$ and $g$ act trivially on $\mathcal{C}W$ by definition.

It remains to consider an unbounded domain \(W\) contained in \(\mathfrak{S}_1\cup\mathfrak{S}_2\).
For \(W\in\mathfrak{S}_1\), the element \(g\) fixes every unbounded domain \(T\in\mathfrak{S}_1\) and acts trivially on \(\mathcal{C}T\), while the restriction of any word in \(f\) and \(g\) to \(\mathfrak{S}_1\) is determined, on unbounded domains, only by the total exponent of \(f\).
Thus the commutator \(h=[a,b]\) acts trivially on \(\mathcal{C}W\).
The same argument, with the roles of \(f\) and \(g\) interchanged, applies to every unbounded \(W\in\mathfrak{S}_2\).

Hence \(\B(h)=\emptyset\), so \(h\) has finite order. 
This contradicts the assumption that \(\langle a,b\rangle\) is a free subgroup of rank \(2\).
Therefore, \(\langle f,g\rangle\) contains no nonabelian free subgroup.
By Theorem~\ref{Thm:TitsAlternative}, \(\langle f,g\rangle\) is virtually abelian.
Consequently, after passing to positive powers, \(f\) and \(g\) commute.
\end{proof}

\begin{proof}[Proof of Theorem~\ref{Thm:MainTheoremStrong_Precise}]
Let $\{f_{1},\dots,f_{m}\}$ be a collection of axial elements that are rigidly fully supported on unbounded domains. By Theorem~\ref{Thm:MainEmbeddingThm_Precise}, it suffices to verify that this collection satisfies Conditions~\eqref{Item:Sec3PowersCommute} and~\eqref{Item:Sec3SubgroupHasUniformlyBoundedOrbits}.
Proposition~\ref{Prop:weakCommutativityForStrongFSE} ensures that Condition~\eqref{Item:Sec3PowersCommute} is satisfied. Meanwhile, Condition~\eqref{Item:Sec3SubgroupHasUniformlyBoundedOrbits} follows immediately from the definition of being rigidly fully supported on an unbounded domain.
\end{proof}

\begin{remark}[HHGs with axial elements rigidly fully supported on unbounded $U \in \mathfrak{S}$]\label{Rmk:F_Ustabilizer}
Although the condition of an axial element being rigidly fully supported is restrictive, such elements exist in many HHGs; for example, in HHGs with the \emph{\(\mathbf F_U\) stabilizers property} in the sense of \cite[Definition~3.1]{OPRAAG}.
For such an HHG, the orthogonal stabilizer \(G_U\) associated to \(U\) carries the induced HHG structure \((G_U,\mathfrak S_U)\).
By \cite[Lemma~3.4(1)]{OPRAAG}, if \(U\) is unbounded, then there exists an axial element \(f_U\in G_U\) which is rigidly fully supported on \(U\).
\end{remark}

\section{Comparison to known results}\label{Sec:Comparison}
We conclude by comparing our results with known RAAG embedding theorems for mapping class groups and RAAGs, the two motivating examples of HHGs.
This recovers the known theorems in these settings and, in several cases, gives slight generalizations.

\subsection{The mapping class group case}\label{Subsec:MappingClassGroupCase}
Building on \cite{MM98} and \cite{Behrstock06}, Behrstock--Hagen--Sisto constructed an HHS structure on the marking complex of a finite-type surface, and hence an HHG structure on its mapping class group.

\begin{theorem}\cite[Theorem~11.1]{BHS19}\label{Thm:StandardHHG of MCG}
Let \(S\) be a finite-type orientable surface with \(\chi(S)<0\), and let \(\mathcal{M}(S)\) be its marking complex. Then \((\mathcal{M}(S),\mathfrak{S})\) is an HHS with the following structure:
\begin{itemize}
\item \(\mathfrak{S}\) is the collection of isotopy classes of essential open subsurfaces of \(S\), possibly disconnected.
\item For each \(U\in\mathfrak{S}\), the associated space \(\mathcal{C}U\) is its curve graph.
\item The relations \(\sqsubseteq\), \(\bot\), and \(\pitchfork\) correspond respectively to nesting, disjointness, and overlapping.
\item For each \(U\in\mathfrak{S}\), the projection \(\pi_U:\mathcal{M}(S)\to \mathcal{C}U\) is the usual subsurface projection.
\item For \(U,V\in\mathfrak{S}\) with either \(U\sqsubsetneq V\) or \(U\pitchfork V\), the projection \(\rho^U_V\) is given by \(\pi_V(\partial U)\subseteq\mathcal{C}V\). If \(V\sqsubsetneq U\), then \(\rho^U_V:\mathcal{C}U\to 2^{\mathcal{C}V}\) is the subsurface projection.
\end{itemize}
Moreover, the action of \(\operatorname{MCG}(S)\) on \(\mathcal{M}(S)\) is proper and cocompact; in particular, \(\operatorname{MCG}(S)\) is an HHG.
\end{theorem}

Fix a marking \(\mu_0\in\mathcal{M}(S)\).
The orbit map $g\mapsto g\mu_0$ is an \(\operatorname{MCG}(S)\)-equivariant quasi-isometry.
We use this orbit map to transport the HHS structure above to an HHG structure on \(\operatorname{MCG}(S)\), which we refer to as the \emph{standard} HHG structure on \(\operatorname{MCG}(S)\).

\begin{remark}
As observed in \cite[Remark~9.3]{HMS26}, although the HHG structure above is typically described for surfaces without boundary, the same construction applies to surfaces with boundary, including braid groups.    
\end{remark}

Before applying our results to the standard HHG structure on \(\operatorname{MCG}(S)\), we compare the terminology used in \cite{Kob12,CLM12,Run21} with that in our paper, namely \emph{pure mapping classes with connected support} and \emph{axial elements in HHGs fully supported on unbounded domains}.

A mapping class \(f\in\operatorname{MCG}(S)\) is said to be \emph{pure} if there exists a (possibly empty) multicurve \(C\), fixed component-wise by \(f\), such that the restriction of \(f\) to each component of \(S\setminus C\) is either the identity or a pseudo-Anosov mapping class. In particular, every Dehn twist about an essential simple closed curve is pure.
If \(f\) is pseudo-Anosov on a component \(U\) of \(S\setminus C\), then \(f\) may also have a multitwist component supported on a submulticurve of \(C\). The \emph{support} of a pure mapping class \(f\) consists of the components of \(S\setminus C\) on which \(f\) acts as a pseudo-Anosov mapping class, together with the annular domains corresponding to those components of \(C\) about which \(f\) has nonzero twist exponent.
If $f$ is supported on an essential non-annular subsurface $U$ but acts with a nontrivial (multi)twist on the boundary of $U$, we regard its support as disconnected. 
Accordingly, by a \emph{pure mapping class with connected support}, we mean either a nonzero power of a Dehn twist about an essential simple closed curve, or the extension by the identity of a pseudo-Anosov mapping class on a connected essential non-annular subsurface, with no additional multitwist component along its boundary.

Every pure mapping class with connected support is an axial element fully supported on the corresponding unbounded domain in the standard HHG structure of \(\operatorname{MCG}(S)\).
Indeed, if \(f\in\operatorname{MCG}(S)\) is pure and supported on a connected essential subsurface \(U\subseteq S\), then the following hold:
\begin{enumerate}
\item
The induced action of \(f\) on \(\mathcal CU\) is loxodromic.
\item
The element \(f\) fixes every subsurface \(V\) disjoint from \(U\). Its action on \(\mathcal CV\) is trivial when \(V\) is non-annular, and has uniformly bounded orbits when \(V\) is an annulus whose core is not contained in \(\partial U\), by \cite[Lemma~3.1]{Mou18}.
\item
If \(V\) is a boundary annulus of \(U\), then the absence of an additional boundary multitwist component implies that the induced action of \(f\) on \(\mathcal CV\) has uniformly bounded orbits.
\end{enumerate}
However, we note that $f$ need not be rigidly fully supported, as explained below, and this is why non-hyperbolic $\operatorname{MCG}(S)$ with the standard HHG structure lies in $\mathcal{P}_2\setminus\mathcal{P}'_2$.

In much of the literature, two pure mapping classes are said to be \emph{irredundant} if they do not share a common nonzero power. We explain why, for pure mapping classes with connected support, irredundancy coincides with geometric irredundancy (Definition~\ref{Def:irredundant}) with respect to \((\operatorname{MCG}(S),\mathfrak S)\). Note that geometric irredundancy also depends on the HHG structure.

Let \(f\) be a pure mapping class with connected non-annular support \(U\), and let \(C\) be the multicurve consisting of the distinct essential curves represented by \(\partial U\).
Let \(Y_0,\ldots,Y_r\) be the components of the surface obtained by cutting \(S\) along \(C\) and replacing the newly created boundary circles by punctures, with \(Y_0\) corresponding to \(U\). The cutting homomorphism gives the central extension \cite[Lemma~2.1(1),(2)]{BLM83}
\[
1\longrightarrow \langle T_c\mid c\in C\rangle
\longrightarrow \operatorname{Stab}^{0}_{\operatorname{MCG}(S)}(U)
\xlongrightarrow{\varphi_U}
\operatorname{MCG}^{0}(Y_0)\times
\prod_{j=1}^{r}\operatorname{MCG}^{0}(Y_j)
\longrightarrow 1,
\]
where \(T_c\) is the Dehn twist about \(c\), the stabilizer preserves each cut component and each side of every curve in \(C\), and \(\operatorname{MCG}^{0}(Y_j)\) fixes each newly created puncture individually.
Then \(\varphi_U(f)=(\bar f,1,\ldots,1)\), where \(\bar f\in\operatorname{MCG}^{0}(Y_0)\) is pseudo-Anosov.

Under our connected-support convention, \(f\) has uniformly bounded orbits on every orthogonal domain, including the boundary annuli, and hence is fully supported on \(U\). Its action on a boundary annular curve graph is nevertheless nontrivial: the pseudo-Anosov restriction moves ideal endpoints in the compactified annular cover, and the defining homotopies fix these endpoints \cite[Section~2.4]{MM98}.
Thus full support here does not imply rigid full support.

Now consider pure mapping classes \(f_1,f_2\) with connected supports \(U_1,U_2\).
Distinct supports give distinct bigsets. If both supports are the same annulus, the elements share a common nonzero power.
If \(U_1=U_2=U\) is non-annular, then \(\operatorname{MCG}^{0}(Y_0)\curvearrowright\mathcal CU\) is acylindrical by \cite[Theorem~1.3]{Bow08}.
If the quasi-axes of the restrictions share an endpoint, then these restrictions generate a virtually cyclic subgroup by \cite[Proposition~6(3)]{BF02}, and hence share a common nonzero power. Such powers differ by a central boundary multitwist, which must be trivial by their bounded-orbit conditions on the boundary annuli. Therefore, for pure mapping classes with connected support, irredundancy coincides with geometric irredundancy.

We can now state an application of our results to mapping class groups using the standard terminology of mapping class group theory.

\begin{theorem}\label{Thm:MCG}
Let $S$ be a finite-type orientable surface with $\chi(S)<0$, and let $\operatorname{MCG}(S)$ denote its mapping class group.
Let $f_i\in \operatorname{MCG}(S)$ be a pure mapping class supported on a connected essential subsurface $U_i\subseteq S$ for $i=1,\dots,m$. 
If $\{ f_1, \dots,f_m\}$ is irredundant, then there exists a constant $D = D(\{ f_i \})>0$ such that for all $d\ge D$,
\[
H = \langle f_1^{d},\dots,f_m^{d}\rangle\cong A(\Gamma),
\]
where $\Gamma$ is the graph with vertex set $\{f_i\}$ in which two vertices $f_i,f_j$ are adjacent if and only if $U_i$ and $U_j$ are disjoint.

Moreover, if each $f_{i}$ is supported on an essential non-annular subsurface (equivalently, none of the $f_{i}$ is a nonzero power of a Dehn twist), then after enlarging $D$ in a controlled way, the subgroup $H$ is undistorted in $\operatorname{MCG}(S)$.
\end{theorem}
\begin{proof}
Consider the standard HHG structure \((\operatorname{MCG}(S),\mathfrak{S})\) fixed above.
The \(\sqsubseteq\)-minimal domains in this structure are annuli, whose associated curve graphs are unbounded. Hence this structure satisfies Assumption~\ref{Assumption:Essential}.

It also satisfies Assumptions~\ref{Assumption:CoboundedAction} and~\ref{Assumption:CentralExtension}. The former follows from the fact that the action of the stabilizer $\operatorname{Stab}_{\operatorname{MCG}(S)}(U)$ of each domain $U\in\mathfrak{S}$ on $\mathcal{C}U$ is cobounded (see also Remark~\ref{Rmk:DiscussionOfAssumptions}).
The clean-container condition (the first part of Assumption~\ref{Assumption:CentralExtension}) is established in \cite[Proposition~7.2]{ABD21}, while the classical cutting homomorphism gives the required \(\mathbb Z\)-central extension; see \cite[Lemma~2.1]{BLM83}. The quotient is naturally a finite-index subgroup of the mapping class group of the cut surface, with the inherited standard HHG structure, which verifies the latter part of Assumption~\ref{Assumption:CentralExtension}.

Note that each $U_i$ is an unbounded domain in $\mathfrak{S}$ such that $f_i$ is an axial element fully supported on $U_i$.
Since $f_i$ and $f_j$ commute in $\operatorname{MCG}(S)$ when $U_{i}$ and $U_{j}$ are orthogonal, the result follows directly from Theorem~\ref{Thm:MainTheoremAssumptiononG_Precise}, with the constant $N(\{f_i\})=1$.
\end{proof}

\begin{remark}\label{Rmk:RunnelsIssues}
We continue the discussion from Remark~\ref{Rmk:MCGcase} and explain the issues in the arguments of \cite{Run21} in relation to the present paper.

First, consider the central form \(w=u_nh_n\cdots u_1h_1\) used in the embedding proof (Section~\ref{Subsec:EmbeddingPart}). The ping-pong argument requires a uniform bound on \(\dist_{U_j}(\phi(h_j)x,\phi(u_jh_j)x)=\dist_{U_j}(x,\phi(u_j)x)\).
At the corresponding point in \cite{Run21}, the bounded-projection result \cite[Lemma~3.1]{Mou18} is invoked to control the commuting word \(u_j\). However, although the image of each syllable of \(u_j\) is supported on a connected subsurface, the total support of \(\phi(u_j)\) may be a disconnected union of pairwise disjoint subsurfaces, whereas Mousley's lemma is stated for a mapping class supported on a connected subsurface. Thus the cited result does not apply directly to the entire element \(\phi(u_j)\).

The missing justification can nevertheless be supplied without changing the statement. Mousley's lemma applies to the image of each individual syllable of \(u_j\), and Lemma~\ref{Lem:UndistortedFreeAbelianSubgroup} combines the resulting bounds to give a uniform estimate for the entire commuting word \(\phi(u_j)\). Since \(u_jh_j\) is a clique subword of a syllable-minimal representative, \(u_j\) has at most \(m\) syllables, so the resulting bound can, for example, be taken to be \(4m\). A disconnected-support version of Mousley's lemma with the same constant would give the sharper bound \(4\), but this sharpening is not needed. Lemma~\ref{Lem:UndistortedFreeAbelianSubgroup} formulates the required combination argument in the general setting of HHGs.


Second, there is a repeated-domain issue in the undistortion part. 
Distinct syllables in the source RAAG may determine translated quasi-axes lying in the same hyperbolic space of the target. 
Consequently, the corresponding projection contributions cannot simply be summed in the distance formula. 
This is the issue addressed in Section~\ref{Subsec:UndistortionPart}, where the strict partial order on syllables is related to a coarse ordering of translated quasi-axes.

Finally, the undistortion argument in \cite{Run21} also uses a bounded-orbit assertion for domains orthogonal to a supporting domain. 
For non-annular support, this is compatible with the argument used here, via Lemma~\ref{Lem:UniformlyBoundedOrbitsforNonNestingMinimal}. 
For annular supporting domains, however, the corresponding assertion fails in general. More precisely, in the proof of \cite[Lemma~5.2]{Mou18}, Mousley considers a Koberda embedding for which \(\phi(a)\) and \(\phi(b)\) are sufficiently high powers of Dehn twists about curves disjoint from \(\eta\), and constructs a word \(g=b^{c_1}a^{c_2}b^{c_3}\) such that \(\phi(g)\) acts loxodromically on the annular curve graph \(\mathcal C\eta\). The same construction yields an infinite family of such examples. Thus even a product of mapping classes supported on domains orthogonal to an annulus can have loxodromic dynamics on the curve graph of that annulus.
This example shows that Condition~\eqref{Item:Sec3SubgroupHasUniformlyBoundedOrbits} need not hold when the supporting domain is \(\sqsubseteq\)-minimal. This is why Theorem~\ref{Thm:MCG} retains the non-annularity assumption in its undistortion conclusion.
\end{remark}

\subsection{RAAG case}\label{Subsec:RAAGCase}
Although it was not the first work to address the embeddability problem between RAAGs \cite{DROMSRAAsubgroup, DYER90}, the paper \cite{KK13} by Kim--Koberda pioneered a combinatorial approach to this problem using the so-called extension graph, which has since influenced subsequent research in this direction \cite{CRDK13, KK14, CR15}. In this subsection, we recover their result through our framework while comparing it with the original proof.

\begin{theorem}\cite[Theorem~1.3]{KK13}\label{Thm:KKthm}
Let $\Gamma,\Lambda$ be finite simplicial graphs and $A(\Gamma),A(\Lambda)$ the associated RAAGs, respectively. 
If $\Lambda$ is an induced subgraph of the extension graph $\Gamma^{e}$ of $\Gamma$, then there exists an undistorted subgroup of $A(\Gamma)$ which is isomorphic to $A(\Lambda)$.
\end{theorem}

The \emph{extension graph} $\Gamma^e$ is defined as follows: its vertices are the conjugates $v^g=gvg^{-1}$ of the standard generators $v$ of $A(\Gamma)$ by elements $g\in A(\Gamma)$, and two vertices are joined by an edge whenever they commute.

Strictly speaking, the original statement of \cite[Theorem~1.3]{KK13} asserts only that there exists an injective homomorphism $A(\Lambda)\rightarrow A(\Gamma)$ whenever $\Lambda \leq \Gamma^e$. 
Here, we outline how the stronger quasi-isometric embedding version in Theorem~\ref{Thm:KKthm} follows from results in \cite{KK13}. By \cite[Lemma~3.1]{KK13}, $\Lambda$ can be realized as an induced subgraph of a graph $\Gamma'$ obtained from $\Gamma$ by finitely many doubling operations, where each step replaces the current graph $\Gamma_i$ with a new one $\Gamma_{i+1}$ such that $A(\Gamma_{i+1})$ is isomorphic to a finite-index subgroup of $A(\Gamma_i)$.
Consequently:
\begin{enumerate}
\item\label{Item:CanonicalEmbedding} there is a canonical embedding $A(\Lambda)\hookrightarrow A(\Gamma')$, which is a quasi-isometric embedding;
\item\label{Item:Finite-indexSubgroup} $A(\Gamma')$ is isomorphic to a finite-index subgroup of $A(\Gamma)$.
\end{enumerate}
Combining Items~\eqref{Item:CanonicalEmbedding} and~\eqref{Item:Finite-indexSubgroup}, we obtain Theorem~\ref{Thm:KKthm}.

Now we give another proof of Theorem~\ref{Thm:KKthm} using Theorem~\ref{Thm:KKStyle}, after recalling the \emph{standard} HHG structure on a RAAG $A(\Gamma)$, as developed in \cite{BHS17I,BR20}.

For an induced subgraph $\Lambda \leq \Gamma$, let $\langle \Lambda \rangle$ denote the subgroup of $A(\Gamma)$ generated by the vertices of $\Lambda$.  
Two left cosets $g\langle \Lambda \rangle$ and $h\langle \Lambda \rangle$ in $A(\Gamma)$ are called \emph{parallel} if $g^{-1}h\in\langle\Lambda\rangle\langle\operatorname{Link}(\Lambda)\rangle$, where $\operatorname{Link}(\Lambda)$ is the subgraph of $\Gamma$ spanned by the vertices adjacent to every vertex in $\Lambda$; the corresponding parallelism class is denoted $g[\Lambda]$.  
The collection
\[\mathfrak{S}_{\Gamma} = \{\, g[\Lambda] \mid g \in A(\Gamma),\ \emptyset\neq\Lambda \leq \Gamma \,\}\]
forms the index set for the \emph{standard} hierarchically hyperbolic structure on $A(\Gamma)$.  
The associated hyperbolic space $\mathcal{C}(g[\Lambda])$ is the factored contact graph \cite[Definition~8.14]{BHS17I}, which is isometric to $\mathcal{C}([\Lambda])$ for $A(\Lambda)$. It is bounded if and only if $\Lambda$ admits a nontrivial join decomposition.  

The relations $\sqsubseteq$ and $\bot$ are defined as follows:  
\begin{itemize}
\item 
$g[\Lambda_{1}] \sqsubseteq h[\Lambda_{2}]$ if $\Lambda_{1} \le \Lambda_{2}$ and there exists $k \in A(\Gamma)$ such that $k[\Lambda_{1}] = g[\Lambda_{1}]$ and $k[\Lambda_{2}] = h[\Lambda_{2}]$.
\item 
$g[\Lambda_{1}] \bot h[\Lambda_{2}]$ if $\Lambda_1$ and $\Lambda_2$ are disjoint, every vertex of $\Lambda_1$ is adjacent to every vertex of $\Lambda_2$, and there exists $k \in A(\Gamma)$ with $k[\Lambda_{1}] = g[\Lambda_{1}]$ and $k[\Lambda_{2}] = h[\Lambda_{2}]$. 
\end{itemize}

Projections of $A(\Gamma)$ to $\mathcal{C}(g[\Lambda])$ are given by nearest-point projection maps, viewing both $A(\Gamma)$ and $g\langle \Lambda\rangle$ as $\operatorname{CAT}(0)$ cube complexes; well-definedness of these projections is one of the main results of \cite{BHS17I}.  
Finally, the left action of $A(\Gamma)$ on itself makes $(A(\Gamma), \mathfrak{S}_\Gamma)$ an HHG.

It is known that $(A(\Gamma),\mathfrak{S}_{\Gamma})$ satisfies the $\mathbf{F}_U$ stabilizers property \cite[Proposition~3.10(2)]{OPRAAG}. Hence, by Remark~\ref{Rmk:F_Ustabilizer}, every unbounded domain admits an axial element that is rigidly fully supported on it.  
For clarity, we sketch a direct argument, without invoking the $\mathbf{F}_U$ stabilizers property, showing that each unbounded domain admits such an element.  
\begin{enumerate}
\item Let $g[\Lambda]\in \mathfrak{S}_{\Gamma}$ be an unbounded domain; in particular, $\Lambda$ does not admit a nontrivial join decomposition. By the Coarse Rank Rigidity Theorem \cite[Theorem~9.14]{DHS17}, there exists an axial element $t \in A(\Lambda)$ which acts loxodromically on $\mathcal{C}([\Lambda])$. Thus
\[gtg^{-1} \in g \langle \Lambda \rangle g^{-1} \le A(\Gamma)\]
acts loxodromically on $\mathcal{C}(g[\Lambda])$.
\item If $h[\Lambda'] \bot g[\Lambda]$, then $gtg^{-1}$ stabilizes $h[\Lambda']$. 
Indeed, by definition, there exists $k \in A(\Gamma)$ such that $h[\Lambda'] = k[\Lambda']$ and $gtg^{-1} = kt_0k^{-1}$ for some $t_0 \in \langle \Lambda \rangle$ that commutes with $\langle \Lambda' \rangle$.  
Hence $gtg^{-1}$ not only fixes $h[\Lambda']$ but also acts trivially on $\mathcal{C}(h[\Lambda'])$.
\end{enumerate}
Thus $gtg^{-1} \in g\langle \Lambda \rangle g^{-1}$ is rigidly fully supported on $g[\Lambda]$, as required.

\begin{proof}[Proof of Theorem~\ref{Thm:KKthm}]
Let $\mathcal{O}^{\mathfrak{S}_\Gamma}$ be the orthogonality graph of $\mathfrak{S}_\Gamma$.  
Define $$\phi:\Gamma^{e}\longrightarrow \mathcal{O}^{\mathfrak{S}_{\Gamma}}\qquad\text{by}\quad\phi(v^{g})=g[v],$$ which is well-defined since $v^{g}=v^{h}$ implies $g^{-1}h$ commutes with $v$, hence $g[v]=h[v]$.  
Moreover, any vertex of the image of $\phi$ corresponds to an unbounded domain in $\mathcal{O}^{\mathfrak{S}_\Gamma}$.
By the Centralizer Theorem \cite{Ser89}, two distinct conjugates \(v^g\) and \(w^h\) commute if and only if $v,w$ are adjacent in $\Gamma$ and there exists $k\in A(\Gamma)$ with $v^{k}=v^{g}$ and $w^{k}=w^{h}$, which is exactly the condition for $g[v]$ and $h[w]$ to be orthogonal. Thus $\phi$ embeds $\Gamma^{e}$ as an induced subgraph of $\mathcal{O}^{\mathfrak{S}_{\Gamma}}$.

If $\Lambda\leq \Gamma^{e}$, then $\phi(\Lambda)$ is an induced subgraph of $\mathcal{O}^{\mathfrak{S}_{\Gamma}}$ and any vertex in $\phi(\Lambda)$ corresponds to an unbounded domain in $\mathfrak{S}_{\Gamma}$. 
Since each unbounded domain of $(A(\Gamma),\mathfrak{S}_\Gamma)$ admits an axial element rigidly fully supported on it, Theorem~\ref{Thm:KKStyle} yields an undistorted subgroup of $A(\Gamma)$ isomorphic to $A(\Lambda)$.  
\end{proof}

\bibliographystyle{alpha}
\bibliography{bibliography.bib}

@article{AB23,
	title={{Conjugator lengths in hierarchically hyperbolic groups}},
	author={Abbott, Carolyn and Behrstock, Jason A.},
	journal={Groups, Geometry, and Dynamics},
	volume={17},
	number={3},
	pages={805--838},
	year={2023}
}

@article{ABD21,
  author       = {Abbott, Carolyn and Behrstock, Jason A. and Durham, Matthew G.},
  title        = {{Largest acylindrical actions and stability in hierarchically hyperbolic groups}},
  journal      = {Transactions of the American Mathematical Society, Series B},
  volume       = {8},
  pages        = {66--104},
  year         = {2021},
  note         = {With an appendix by Daniel Berlyne and Jacob Russell},
  doi          = {10.1090/btran/71},
  url          = {https://doi.org/10.1090/btran/71}
}

@article{ABR25,
  author  = {Abbott, Carolyn and Behrstock, Jason and Russell, Jacob},
  title   = {Structure invariant properties of the hierarchically hyperbolic boundary},
  journal = {Journal of Topology and Analysis},
  volume  = {17},
  number  = {3},
  pages   = {859--903},
  year    = {2025},
  doi     = {10.1142/S1793525323500516}
}

@article{Behrstock06,
  author    = {Jason A. Behrstock},
  title     = {{Asymptotic geometry of the mapping class group and Teichm{\"u}ller space}},
  journal   = {Geometry \& Topology},
  volume    = {10},
  year      = {2006},
  pages     = {1523--1578},
  doi       = {10.2140/gt.2006.10.1523}
}

@article{BF02,
  author  = {Bestvina, Mladen and Fujiwara, Koji},
  title   = {Bounded cohomology of subgroups of mapping class groups},
  journal = {Geometry \& Topology},
  volume  = {6},
  year    = {2002},
  pages   = {69--89},
  doi     = {10.2140/gt.2002.6.69}
}

@article{BHS17smallcancellation,
  title={Asymptotic dimension and small-cancellation for hierarchically hyperbolic spaces and groups},
  author={Behrstock, Jason A. and Hagen, Mark F. and Sisto, Alessandro},
  journal={Proceedings of the London Mathematical Society},
  volume={114},
  number={5},
  pages={890-926},
  year={2017},
  publisher={Wiley Online Library}
}

@article{Man24,
  title={Short hierarchically hyperbolic groups I: uncountably many coarse median structures},
  author={Mangioni, Giorgio},
  journal={arXiv preprint arXiv:2410.09232},
  year={2024}
}

@article{BHMS24,
author = {Behrstock, Jason A. and Hagen, Mark F. and Martin, Alexandre and Sisto, Alessandro},
title = {A combinatorial take on hierarchical hyperbolicity and applications to quotients of mapping class groups},
journal = {Journal of Topology},
volume = {17},
number = {3},
pages = {e12351},
doi = {https://doi.org/10.1112/topo.12351},
url = {https://londmathsoc.onlinelibrary.wiley.com/doi/abs/10.1112/topo.12351},
year = {2024}
}

@article{BHS21,
author = {Behrstock, Jason A. and Hagen, Mark F. and Sisto, Alessandro},
title = {{Quasiflats in hierarchically hyperbolic spaces}},
volume = {170},
journal = {Duke Mathematical Journal},
number = {5},
publisher = {Duke University Press},
pages = {909-996},
year = {2021},
doi = {10.1215/00127094-2020-0056},
URL = {https://doi.org/10.1215/00127094-2020-0056}
}

@article{BR20,
	title={{Hierarchical hyperbolicity of graph products}},
	author={Berlyne, Daniel and Russell, Jacob},
	journal={Groups, Geometry, and Dynamics},
	year={2022},
    volume={16},
    number={2},
    pages={523-580},
    doi={DOI 10.4171/GGD/652}
}

@article{BHS17I,
	title={{Hierarchically hyperbolic spaces, I: Curve complexes for cubical groups}},
	author={Behrstock, Jason A. and Hagen, Mark F. and Sisto, Alessandro},
	journal={Geometry \& Topology},
	volume={21},
	number={3},
	pages={1731-1804},
	year={2017},
	publisher={Mathematical Sciences Publishers}
}

@article{BHS19,
	title={{Hierarchically hyperbolic spaces II: Combination theorems and the distance formula}},
	author={Behrstock, Jason A. and Hagen, Mark F. and Sisto, Alessandro},
	journal={Pacific Journal of Mathematics},
	volume={299},
	number={2},
	pages={257-338},
	year={2019},
	publisher={Mathematical Sciences Publishers}
}

@book{BH99,
author = {Bridson, Martin and Haefliger, Andr\'e},
year = {1999},
title = {Metric Spaces of Non-Positive Curvature},
isbn = {978-3-642-08399-0},
publisher = {Grundlehren der mathematischen Wissenschaften 319, Springer},
doi = {10.1007/978-3-662-12494-9}
}

@article{Bow08,
    author = {Bowditch, Brian H.},
    title = {Tight geodesics in the curve complex},
    journal = {Inventiones mathematicae},
    year = {2008},
    volume = {171},
    number = {2},
    pages = {281-300}
}

@misc{Bow19,
  author = {Brian H. Bowditch},
  title = {{Quasiflats in coarse median spaces}},
  year = 2019,
howpublished = {Preprint available at {\href{https://bhbowditch.com/papers/quasiflats.pdf}{bhbowditch.com/papers/quasiflats.pdf}}}
}

@misc{CHK,
    title={Real cubings and asymptotic cones of hierarchically hyperbolic groups},
    author={Casals-Ruiz, Montserrat and Hagen, Mark F. and Kazachkov, Ilya},
howpublished = {Preprint available at {\href{https://www.wescac.net/cones.html}{www.wescac.net/cones.html}}}
}

@article{Cou16,
     author = {Coulon, R\'emi B.},
     title = {Partial periodic quotients of groups acting on a hyperbolic space},
     journal = {Annales de l'Institut Fourier},
     pages = {1773-1857},
     publisher = {Association des Annales de l{\textquoteright}institut Fourier},
     volume = {66},
     number = {5},
     year = {2016},
     doi = {10.5802/aif.3050},
     language = {en},
     url = {https://aif.centre-mersenne.org/articles/10.5802/aif.3050/}
}

@article {CSS08,
    AUTHOR = {Crisp, John and Sageev, Michah and Sapir, Mark},
     TITLE = {{Surface subgroups of right-angled Artin groups}},
   JOURNAL = {Internat. J. Algebra Comput.},
  FJOURNAL = {International Journal of Algebra and Computation},
    VOLUME = {18},
      YEAR = {2008},
    NUMBER = {3},
     PAGES = {443-491},
      ISSN = {0218-1967,1793-6500},
   MRCLASS = {20F36 (05C10 05C25 57M07 57M50)},
  MRNUMBER = {2422070},
MRREVIEWER = {Robert\ W.\ Bell},
       DOI = {10.1142/S0218196708004536},
       URL = {https://doi.org/10.1142/S0218196708004536},
}

@article{CP01,
author = {Crisp, John and Paris, Luis},
year = {2001},
pages = {19-36},
title = {{The solution to a conjecture of Tits on the subgroup generated by the squares of the generators of an Artin group}},
volume = {145},
issue = {1},
journal = {Inventiones mathematicae},
doi = {10.1007/s002220100138}
}

@article{CW,
author = {Crisp, John and Wiest, Bert},
year = {2004},
pages = {439-472},
title = {{Embeddings of graph braid and surface groups in right-angled Artin groups and braid groups}},
volume = {4},
journal = {Algebraic \& Geometric Topology},
doi = {10.2140/agt.2004.4.439}
}

@article{BLM83,
  title={Abelian and solvable subgroups of mapping class group},
  author={BIRMAN, Joan S. and Lubotzky, Alex and McCarthy, John},
  journal={Duke mathematical journal},
  volume={50},
  number={4},
  pages={1107--1120},
  year={1983}
}

@article{CLM12,
title={{The geometry of right-angled Artin subgroups of mapping class groups}},
author={Clay, Matt T. and Leininger, Christopher J. and Mangahas, Johanna},
journal={Groups, Geometry, and Dynamics},
volume={6},
number={2},
pages={249-278},
year={2012},
doi={10.4171/GGD/157}
}

@article{DSS89,
  title={{Surface subgroups of graph groups}},
  author={Droms, Carl and Servatius, Brigitte and Servatius, Herman},
  journal={Proceedings of the American Mathematical Society},
  volume={106},
  number={3},
  pages={573–578},
  year={1989},
  publisher={}
}

@article{DHS17,
  title={{Boundaries and automorphisms of hierarchically hyperbolic spaces}},
  author={Durham, Matthew and Hagen, Mark F. and Sisto, Alessandro},
  journal={Geometry \& Topology},
  volume={21},
  number={6},
  pages={3659-3758},
  year={2017},
  publisher={Mathematical Sciences Publishers}
}

@article{DHS20,
  title={{Correction to the article Boundaries and automorphisms of hierarchically hyperbolic spaces}},
  author={Durham, Matthew and Hagen, Mark F. and Sisto, Alessandro},
  journal={Geometry \& Topology},
  volume={24},
  number={2},
  pages={1051-1073},
  year={2020},
  publisher={Mathematical Sciences Publishers}
}

@article{DL24,
  title={{Right-angled Artin subgroups of right-angled Coxeter and Artin groups}},
  author={Dani, Pallavi and Levcovitz, Ivan},
  journal={Algebraic \& Geometric Topology},
  volume={24},
  number={2},
  pages={755-802},
  year={2024},
  doi={10.2140/agt.2024.24.755}
}

@article{DYER90,
title = {{Reflection subgroups of Coxeter systems}},
journal = {Journal of Algebra},
volume = {135},
number = {1},
pages = {57-73},
year = {1990},
issn = {0021-8693},
doi = {https://doi.org/10.1016/0021-8693(90)90149-I},
url = {https://www.sciencedirect.com/science/article/pii/002186939090149I},
author = {Dyer, Matthew}
}

@article{HM95,
title = {{Algorithms and Geometry for Graph Products of Groups}},
journal = {Journal of Algebra},
volume = {171},
number = {1},
pages = {230-257},
year = {1995},
issn = {0021-8693},
doi = {https://doi.org/10.1006/jabr.1995.1010},
url = {https://www.sciencedirect.com/science/article/pii/S0021869385710101},
author = {S. Hermiller and J. Meier}
}

@article{HHP23,
title={{Coarse injectivity, hierarchical hyperbolicity and semihyperbolicity}},
author={Haettel, Thomas and Hoda, Nima and Petyt, Harry},
journal={Geometry \& Topology},
volume={27},
issue={4},
pages={1587–1633},
year={2023}
}

@article{PS23,
	title={{Unbounded domains in hierarchically hyperbolic groups}},
	author={Petyt, Harry and Spriano, Davide},
	journal={Groups, Geometry, and Dynamics},
	volume={17},
	number={2},
	pages={479-500},
	year={2023}
}

@article{Ser89,
title = {Automorphisms of graph groups},
journal = {Journal of Algebra},
volume = {126},
number = {1},
pages = {34-60},
year = {1989},
issn = {0021-8693},
doi = {https://doi.org/10.1016/0021-8693(89)90319-0},
url = {https://www.sciencedirect.com/science/article/pii/0021869389903190},
author = {Herman Servatius}
}

@article{Seo21,
author = {Seo, Donggyun},
year = {2021},
pages = {1511-1533},
title = {{Powers of Dehn twists generating right-angled Artin groups}},
volume = {21},
issue = {3},
journal = {Algebraic \& Geometric Topology},
doi = {10.2140/agt.2021.21.1511}
}

@article{HMS26,
	title={{A combinatorial structure for many hierarchically hyperbolic spaces}},
	author={Hagen, Mark and Mangioni, Giorgio and Sisto, Alessandro},
	journal={Pacific Journal of Mathematics},
	volume={341},
	number={2},
	pages={305--377},
	year={2026},
	publisher={Mathematical Sciences Publishers}
}

@article{MM98,
  title={{Geometry of the complex of curves II: Hierarchical structure}},
  author={Masur, Howard A and Minsky, Yair N},
  volume = {10},
  journal = {Geometric and Functional Analysis},
  number = {4},
  pages = {902-974},
  year = {2000},
  doi = {10.1007/PL00001643}
}

@article{Run21,
  title={{Effective generation of right-angled Artin groups in mapping class groups}},
  author={Runnels, Ian},
  journal={Geometriae Dedicata},
  volume={214},
  pages={277-294},
  year={2021},
  publisher={Springer}
}

@article{ANSGP24,
  title={{Hierarchically hyperbolic groups and uniform exponential growth}},
  author={Abbott, Carolyn and Ng, Thomas and Spriano, Davide and Gupta, Radhika and Petyt, Harry},
  journal={Mathematische Zeitschrift},
  volume={306},
  number={1},
  eid={18},
  year={2024},
  doi={10.1007/s00209-023-03411-6}
}

@article{DMS23,
  title={{Stable cubulations, bicombings, and barycenters}},
  author={Durham, Matthew G and Minsky, Yair N and Sisto, Alessandro},
  journal={Geometry \& Topology},
  volume={27},
  number={6},
  pages={2383-2478},
  year={2023},
  publisher={Mathematical Sciences Publishers}
}

@article{KK13,
  title={{Embedability between right-angled Artin groups}},
  author={Kim, Sang-hyun and Koberda, Thomas},
  journal={Geometry \& Topology},
  volume={17},
  number={1},
  pages={493-530},
  year={2013},
  publisher={Mathematical Sciences Publishers}
}

@article{KK14,
    author = {Kim, Sang-Hyun and Koberda, Thomas},
    title = {{An Obstruction to Embedding Right-Angled Artin Groups in Mapping Class Groups}},
    journal = {International Mathematics Research Notices},
    volume = {2014},
    number = {14},
    pages = {3912-3918},
    year = {2014},
    month = {04},
    issn = {1073-7928},
    doi = {10.1093/imrn/rnt064},
    url = {https://doi.org/10.1093/imrn/rnt064},
    eprint = {https://academic.oup.com/imrn/article-pdf/2014/14/3912/18890184/rnt064.pdf},
}

@article{Kob12,
  title={{Right-angled Artin groups and a generalized isomorphism problem for finitely generated subgroups of mapping class groups}},
  author={Koberda, Thomas},
  journal={Geometric and Functional Analysis},
  volume={22},
  pages={1541-1590},
  year={2012},
  publisher={Springer}
}

@article{CR15,
  title={{Embeddability and universal theory of partially commutative groups}},
  author={Casals-Ruiz, Montserrat},
  journal={International Mathematics Research Notices},
  volume={2015},
  number={24},
  pages={13575-13622},
  year={2015},
  publisher={Oxford University Press}
}

@article{CRDK13,
title = {{Embedddings between partially commutative groups: Two counterexamples}},
journal = {Journal of Algebra},
volume = {390},
pages = {87-99},
year = {2013},
issn = {0021-8693},
doi = {https://doi.org/10.1016/j.jalgebra.2013.04.036},
url = {https://www.sciencedirect.com/science/article/pii/S0021869313002779},
author = {Casals-Ruiz, Montserrat and Duncan, Andrew and Kazachkov, Ilya}
}

@article{JS22,
  title={{Right-angled Artin subgroups of Artin groups}},
  author={Jankiewicz, Kasia and Schreve, Kevin},
  journal={Journal of the London Mathematical Society},
  volume={106},
  number={2},
  pages={818-854},
  year={2022},
  publisher={Wiley Online Library}
}

@article{BBM20, 
title={{Embedding right-angled Artin groups into Brin–Thompson groups}}, 
volume={169}, 
DOI={10.1017/S0305004119000112}, 
number={2}, 
journal={Mathematical Proceedings of the Cambridge Philosophical Society}, 
author={Belk, James and Bleak, Colin and Matucci, Francesco}, 
year={2020}, 
pages={225–229}
}

@article{Sal07,
  title={{Embedding right-angled Artin groups into graph braid groups}},
  author={Sabalka, Lucas},
  journal={Geometriae Dedicata},
  volume={124},
  pages={191-198},
  year={2007},
  publisher={Springer}
}

@article{KL08,
  title={Shadows of mapping class groups: capturing convex cocompactness},
  author={Kent, Autumn E and Leininger, Christopher J},
  journal={Geometric and Functional Analysis},
  volume={18},
  number={4},
  pages={1270-1325},
  year={2008},
  publisher={Springer}
}

@book{KS22,
  title={Trees of hyperbolic spaces},
  author={Kapovich, Michael and Sardar, Pranab},
  volume={282},
  year={2024},
  publisher={American Mathematical Society}
}

@article{WY25,
title = {Uniform exponential growth for groups with proper product actions on hyperbolic spaces},
journal = {Journal of Algebra},
volume = {676},
pages = {189-235},
year = {2025},
issn = {0021-8693},
doi = {https://doi.org/10.1016/j.jalgebra.2025.03.028},
url = {https://www.sciencedirect.com/science/article/pii/S0021869325001814},
author = {Wan, Renxing and Yang, Wenyuan}
}

@article{Mou18,
author = {Sarah C Mousley},
title = {{Nonexistence of boundary maps for some hierarchically hyperbolic spaces}},
volume = {18},
journal = {Algebraic \& Geometric Topology},
number = {1},
publisher = {MSP},
pages = {409 -- 439},
year = {2018},
doi = {10.2140/agt.2018.18.409},
URL = {https://doi.org/10.2140/agt.2018.18.409}
}

@article{DROMSRAAsubgroup,
title = {Subgroups of graph groups},
journal = {Journal of Algebra},
volume = {110},
number = {2},
pages = {519-522},
year = {1987},
issn = {0021-8693},
doi = {https://doi.org/10.1016/0021-8693(87)90063-9},
url = {https://www.sciencedirect.com/science/article/pii/0021869387900639},
author = {Carl Droms}
}

@misc{OPRAAG,
title={{From algebraic orthogonality to RAAG embedding obstructions in hierarchically hyperbolic groups}},
author={Oh, Sangrok and Park, Jihoon},
year={2026},
eprint={2608.09535},
archivePrefix={arXiv}
}
 
\end{document}